\documentclass[journal,onecolumn,draftclsnofoot]{IEEEtran}
\usepackage{cite}
\usepackage{amsmath,amssymb,amsfonts}
\usepackage{graphicx}
\usepackage{xcolor}

\makeatletter
\newcommand*\rel@kern[1]{\kern#1\dimexpr\macc@kerna}
\newcommand*\widebar[1]{%
	\begingroup
	\def\mathaccent##1##2{%
		\rel@kern{0.8}%
		\overline{\rel@kern{-0.8}\macc@nucleus\rel@kern{0.2}}%
		\rel@kern{-0.2}%
	}%
	\macc@depth\@ne
	\let\math@bgroup\@empty \let\math@egroup\macc@set@skewchar
	\mathsurround\z@ \frozen@everymath{\mathgroup\macc@group\relax}%
	\macc@set@skewchar\relax
	\let\mathaccentV\macc@nested@a
	\macc@nested@a\relax111{#1}%
	\endgroup
}

\DeclareMathOperator*{\argmin}{argmin}

\DeclareMathOperator*{\minimize}{minimize}
\DeclareMathOperator*{\maximize}{maximize}
\DeclareMathOperator*{\subjectto}{subject~to}

\DeclareMathOperator*{\otherwise}{otherwise}

\newcommand{\reals}[1]{\mathbb{R}^{#1}}
\newcommand{\cN}{\mathcal{N}}
\DeclareMathOperator{\dom}{\mathbf{dom}}

\DeclareMathOperator{\diag}{diag}

\DeclareMathOperator{\blkdiag}{blkdiag}
\DeclareMathOperator{\prox}{\mathbf{prox}}
\DeclareMathOperator{\sign}{sign}
\DeclareMathOperator{\dist}{\mathbf{dist}}

\DeclareMathOperator{\indicator}{I}

\newcommand{\norm}[1]{\| #1 \|}    
\newcommand{\inner}[2]{\langle #1,#2\rangle}    

\newcommand{\blue}[1]{#1}

\newcommand{\ds}{\displaystyle}
\newcommand{\set}[1]{\{ #1 \}}
\newcommand{\enma}[1]{\ensuremath{#1}}

\newcommand{\non}{\nonumber}

\newcommand{\beq}{\begin{equation}}
\newcommand{\eeq}{\end{equation}}
\newcommand{\ba}{\begin{array}}
\newcommand{\ea}{\end{array}}
\newcommand{\bseq}{\begin{subequations}}
\newcommand{\eseq}{\end{subequations}}

\newcommand{\DefinedAs}[0]{\mathrel{\mathop:}=}

\newcommand{\rmd}{\mathrm{d}}
\newcommand{\rme}{\mathrm{e}}

\newcommand{\matbegin}{\left[}
\newcommand{\matend}{\right]}

\newcommand{\tbo}[2]{
		\matbegin \begin{array}{c}
				#1 \\ #2
			\end{array} \matend }
\newcommand{\thbo}[3]{
		\matbegin \begin{array}{c}
				#1 \\ #2 \\ #3
			\end{array} \matend }

\newcommand{\obth}[3]{
	\matbegin \begin{array}{ccc}
		#1 & #2 & #3
	\end{array} \matend }

\newcommand{\thbt}[6]{
		\matbegin \begin{array}{cc}
				#1 & #2 \\ #3 & #4 \\ #5 & #6
			\end{array} \matend }
\newcommand{\tbth}[6]{
		\matbegin \begin{array}{ccc}
				#1 & #2 & #3\\ #4 & #5 & #6
			\end{array} \matend }
\newcommand{\thbth}[9]{
		\matbegin \begin{array}{ccc}
				#1 & #2 & #3 \\
				#4 & #5 & #6 \\
				#7 & #8 & #9
			\end{array}\matend}

\newcommand{\ww}{w}
\newcommand{\nn}{\nu}
\newcommand{\zz}{z}
\newcommand{\yy}{y}
\newcommand{\pp}{\psi}
\newcommand{\PP}{\Psi}
\newcommand{\lam}{\lambda}

\newcommand{\xs}{x^\star}
\newcommand{\ws}{\ww^\star}
\newcommand{\zs}{\zz^\star}
\newcommand{\ys}{\yy^\star}
\newcommand{\lams}{\lam^\star}
\newcommand{\pps}{\pp^\star}
\newcommand{\PPs}{\PP^\star}
\newcommand{\xss}{x^\ast}

\newcommand{\zss}{\zz^\ast}
\newcommand{\yss}{\yy^\ast}
\newcommand{\lamss}{\lam^\ast}

\newcommand{\wb}{\widebar{\ww}}
\newcommand{\xb}{\widebar{x}}
\newcommand{\zb}{\widebar{\zz}}

\newcommand{\tb}{\widebar{t}}
\newcommand{\ab}{\widebar{\alpha}}

\newcommand{\xd}{\dot{x}}
\newcommand{\zd}{\dot{\zz}}
\newcommand{\yd}{\dot{\yy}}
\newcommand{\lamd}{\dot{\lam}}

\newcommand{\xt}{\tilde{x}}
\newcommand{\yt}{\tilde{\yy}}
\newcommand{\zt}{\tilde{\zz}}
\newcommand{\lamt}{\tilde{\lam}}
\newcommand{\wt}{\widetilde{\ww}}
\newcommand{\ft}{\widetilde{f}}
\newcommand{\xbt}{\widetilde{\xb}}
\newcommand{\zbt}{\widetilde{\zb}}
\newcommand{\wbt}{\widetilde{\wb}}

\newcommand{\xbp}{\xb^\prime}
\newcommand{\zbp}{\zb^\prime}
\newcommand{\wbp}{\wb^\prime}

\newcommand{\xh}{\widehat{x}}
\newcommand{\zh}{\widehat{\zz}}
\newcommand{\wh}{\widehat{\ww}}

\newcommand{\xp}{x^\prime}
\newcommand{\zp}{\zz^\prime}
\newcommand{\yp}{\yy^\prime}
\newcommand{\lamp}{\lam^\prime}

\newcommand{\xas}{\vec{x}^\star}
\newcommand{\zas}{\vec{\zz}^\star}
\newcommand{\yas}{\vec{\yy}^\star}
\newcommand{\lamas}{\vec{\lam}^\star}
\newcommand{\pas}{\vec{\pp}^\star}

\newcommand{\cX}{\enma{\mathcal X}}

\newcommand{\cL}{\enma{\mathcal L}}
\newcommand{\cLm}{\enma{\mathcal L}_\mu}
\newcommand{\csLm}{\enma{\mathcal L}_\mu^s}
\newcommand{\more}[2]{M_{#1}(#2)}
\newcommand{\moreau}[1]{\more{\mu g}{#1}}
\newcommand{\dual}{d}
\newcommand{\dopt}{\dual^\star}
\newcommand{\cDs}{\enma{\mathcal D^\star}}
\newcommand{\cPs}{\enma{\mathcal P^\star}}
\newcommand{\cPsw}{\enma{\mathcal P^\star}_{\ww}}
\newcommand{\PPsw}{\PP^\star_\ww}
\newcommand{\cPw}{\enma{\mathcal P_w}}

\newcommand{\cE}{\enma{\mathcal E}}

\newcommand{\cG}{\enma{\mathcal G}}

\newcommand{\cM}{\enma{\mathcal M}}
\newcommand{\cC}{\enma{\mathcal C}}
\newcommand{\cI}{\enma{\mathcal I}}
\newcommand{\cD}{\enma{\mathcal D}}

\newcommand{\cP}{\enma{\mathcal P}}
\newcommand{\cR}{\enma{\mathcal R}}
\newcommand{\cS}{\enma{\mathcal S}}
\newcommand{\cJ}{\enma{\mathcal J}}
\newcommand{\sigu}{\widebar{\sigma}}
\newcommand{\sigl}{\underline{\sigma}}

\usepackage[colorlinks=true, allcolors=blue]{hyperref}

\usepackage{lipsum}

\newtheorem{theorem}{Theorem}

\newtheorem{lemma}{Lemma}
\newtheorem{remark}{Remark}
\newtheorem{assumption}{Assumption}

\newtheorem{corollary}{Corollary}
\newenvironment{mylemma}[1]
  {\count@\c@lemma
   \global\c@lemma#1 %
    \global\advance\c@lemma\m@ne
   \lemma}
  {\endlemma
   \global\c@lemma\count@}

\newcommand{\nvs}{0.1cm}
\newcommand{\vsp}{\vspace*{0.1cm}}
\newcommand{\asp}[1]{\!\!\! #1 \!\!\!}

\renewcommand{\baselinestretch}{.99} 

\IEEEoverridecommandlockouts                              
\begin{document}
\title{Stability of Primal-Dual Gradient Flow Dynamics for Multi-Block Convex Optimization Problems}
\author{Ibrahim K.\ Ozaslan, Panagiotis Patrinos, and  Mihailo R.\ Jovanovi\'c
	\thanks{I.\ K.\ Ozaslan and M.\ R.\ Jovanovi\'c are with the Ming Hsieh Department of Electrical and Computer Engineering, University of Southern California, Los Angeles, CA 90089, USA. P.\ Patrinos is with the Department of Electrical Engineering (ESAT-STADIUS), KU Leuven, 3001 Leuven, Belgium. E-mails: ozaslan@usc.edu, panos.patrinos@esat.kuleuven.be, mihailo@usc.edu.}}

\maketitle

	\vspace*{-4ex}
\begin{abstract}
We examine stability properties of primal-dual gradient flow dynamics for composite convex optimization problems with multiple, possibly nonsmooth, terms in the objective function under the generalized consensus constraint. The proposed dynamics are based on the proximal augmented Lagrangian and they provide a viable alternative to ADMM which faces significant challenges from both analysis and implementation viewpoints in large-scale multi-block scenarios. In contrast to  customized algorithms with individualized convergence guarantees, we develop a systematic approach for solving a broad class of challenging composite optimization problems. We leverage various structural properties to establish global (exponential) convergence guarantees for the proposed dynamics. Our assumptions are much weaker than those required to prove (exponential) stability of primal-dual dynamics as well as (linear) convergence of discrete-time methods such as standard two-block and multi-block ADMM and EXTRA algorithms. Finally, we show necessity of some of our structural assumptions for exponential stability and provide computational experiments to demonstrate the convenience of the proposed approach for parallel and distributed computing applications.
\end{abstract}
	\vspace*{-1ex}

\section{Introduction}\label{sec.intro}
We study the composite constrained optimization problems of the form
	\begin{subequations}
	\label{eq.intro}
\beq\label{eq.intro1}
	\ba{rcl}
	\minimize\limits_{x, \, \zz}
	& \asp{} &
	{f}(x)   \;+\;  g(\zz) 
	\\[\nvs]
	\subjectto 
	& \asp{} &
	Ex \,+\, F\zz \,-\, q  \,=\,  0
	\ea
\eeq
where $x\in\reals{m}$ and $\zz\in\reals{n}$ are the optimization variables, $E\in\reals{p\times m}$, $F\in \reals{p\times n}$, and $q\in \reals{p}$ are the problem data, and $f:$ $\reals{m}\to\reals{}$, $g:$ $\reals{n}\to\reals{}\cup\set{\pm\infty}$ are the separable convex functions given by
	\beq\label{eq.intro2}
		f(x) 
		 \,=\,  
		\sum\limits_{i \,=\, 1}^k f_i(x_i), 
		~~
		g(\zz) 
		 \,=\,  
		\sum\limits_{j\,=\,1}^\ell g_j(\zz_j).
	\eeq
	\end{subequations}
Depending on $f_i$'s and $g_j$'s, the optimization variables in~\eqref{eq.intro} may have arbitrary partitions $x = [x_1^T~\cdots~x_k^T]^T$ and $\zz = [\zz_1^T~\cdots~\zz_\ell^T]^T$, which induce conformable partitions of $E = [{E_1}{~\cdots~}{E_k}]$ and $F = [{F_1}{~\cdots~}{F_\ell}]$. We denote the set of solutions by~$\cPs$ and assume that it is nonempty. Furthermore, we let each $f_i$ be a convex function with a Lipschitz continuous gradient (i.e., smooth) and each $g_j$ be a closed proper convex (possibly nondifferentiable) function with efficiently computable proximal \blue{operator}. Examples of $g_j$ include indicator functions of convex sets, support functions, group-lasso, as well as $\ell_1$, $\ell_2$, $\ell_\infty$, and nuclear norms~\cite{parboy13}. While we do not assume existence of any smooth term in the objective function (i.e., we allow $k = 0$ in~\eqref{eq.intro2}), if a smooth term does exist, it should be included in the $x$-block rather than in the $z$-block. This separation between smooth and nonsmooth parts of the objective function plays an important role in identifying weakest set of assumptions that are required to establish our results; it also alleviates cumbersome notation resulting from the introduction of auxiliary variables.

\vsp

Since function $g$ is allowed to be nondifferentiable, a wide range of constraints can be included into problem~\eqref{eq.intro}. In particular, convex constraints $x_i \in \cX_i$ for some $i\in\set{1,\dots,k}$ can be easily incorporated into~\eqref{eq.intro} by augmenting the objective function with the indicator function of set $\cX_i$,
	\beq
		\indicator_{\cX_i}(x) \;\DefinedAs\; 
		\left\{
		\ba{ll}
		0, 
		& 
		x \, \in\,\cX_i
		\\[0.0cm] 
		\infty,
		& 
		 \otherwise\! .
		\ea
		\right.
		\non
	\eeq
For example, linear inequality constraint $G_ix_i \leq  h_i$ can be handled by introducing a slack variable $\zz_{\ell+1} \geq 0$, converting the inequality constraint to equality constraint $G_i x_i + \zz_{\ell + 1} = h_i$, and by adding the indicator function associated with the positive orthant, $g_{\ell+1}(\zz_{\ell+1}) = \indicator_{\reals{}_+} \! (\zz_{\ell+1})$, to the objective function in~\eqref{eq.intro}. Furthermore, even nondifferentiable convex inequality constraints can be included in~\eqref{eq.intro} as long as the projection operator associated with the constraint set is easily computable.

\vsp

Optimization problem~\eqref{eq.intro} arises in a host of applications ranging from signal processing and machine learning to statistics and control theory; see Section~\ref{sec.numeric} for detailed examples. A particularly important class of problems captured by~\eqref{eq.intro} is the regularized consensus problem~\cite{boyparchupeleck11} in which $k$ agents in a connected undirected network aim to cooperatively solve 
\bseq\label{ex.dist_opt}
\beq\label{ex.dist_opt.ori}
\ba{rcl}
\minimize\limits_{\xt}
&\asp{}&
\ds\sum\limits_{i  \,=\,  1}^{k} \left( f_i(\xt) \,+\, g_i(C_i\xt) \right)
\ea
\eeq
where the matrix $C_i\in\reals{\widetilde{m}\times\widetilde{n}}$ and convex functions $f_i$: $\reals{\widetilde{m}}\to\reals{}$ and $g_i$: $\reals{\widetilde{n}}\to\reals{}\cup\set{\pm\infty}$, with the former being smooth, are known only by agent $i$. Each node in the network represents an agent and each edge represents a communication channel between two agents. The information exchange between two agents occurs only if there exists an edge between the corresponding nodes. 
	\blue{To cast this problem into the general form~\eqref{eq.intro}, we follow a standard distributed optimization approach~\cite{mer94,boyparchupeleck11,bertsi15}. Specifically, we assign a local decision variable $x_i$ to each agent and enforce the consensus constraint $x_1=\cdots=x_k$ by imposing $Tx=0$ on the stacked variable $x\in\reals{k\widetilde{m}}$, where $T^T$ is the incidence matrix of the underlying communication network~\cite{mer94}. We further apply variable splitting of the form $C_i x_i = z_i$, resulting in}
\beq\label{ex.dist_opt.equ}
\ba{rcl}
\minimize\limits_{x,\, \zz}
&\asp{}&
f(x) \,+\, g(\zz)
\\[\nvs]
\subjectto 
&\asp{}&
\tbo{T}{C}x \,+\, \tbo{0}{-I}\zz \,=\, 0
\ea
\eeq
\eseq
\blue{where the stacked variable $z\in\reals{k\widetilde{n}}$, and $C$ is the block-diagonal matrix defined as $C \DefinedAs \blkdiag \, (C_1,\ldots,C_k)$.}

\vsp

\blue{Beyond constrained formulations, such as~\eqref{ex.dist_opt}, the unconstrained composite optimization problems of the form}
	\beq\label{eq.unc_problem}
		\minimize\limits_x~f(x) \;+\;  \sum\limits_{j \, = \, 1}^\ell g_j(T_jx) 
	\eeq      
	can also be brought into~\eqref{eq.intro} by setting $z_j    =   T_jx$, $q = 0$, and
	\beq
		E    =    \left[\!
				\begin{array}{c}
				{T_1}
				\\
				{T_2}
				\\[-1ex]
				{\vdots}
				\\
				{T_\ell}
				\end{array}
				\!\right]\!,
		~
		F_1   =    
		\left[\!
		\begin{array}{c}
		{-I}
		\\
		{0}
		\\[-1ex]
		{\vdots}
		\\
		{0}
		\end{array}
		\!\right]\!,
		\dots~
		F_\ell    =    
		\left[\!
				\begin{array}{c}
				{0}
				\\[-1ex]
				{\vdots}
				\\
				{0}
				\\
				{-I}
				\end{array}
				\!\right]\!, ~  j = 1,\dots,\ell.
		\non
	\eeq

\vsp

\blue{Splitting methods provide an effective means for solving the class problems that can be brought into the form of~\eqref{eq.intro} by facilitating separate treatment of different blocks.} If the problem is properly formulated, these methods are also convenient for distributed computations and parallelization. For example, the Alternating Direction of Method of Multipliers (ADMM), which represents a particular instance of more general splitting techniques~\cite{chapoc16, compes12,con13, vu13, latpat17, davyin17}, has attracted significant attention because of its straightforward and efficient implementation~\cite{boyparchupeleck11}.

\vsp

The multi-block ADMM for problem~\eqref{eq.intro} takes the form,
\bseq\label{eq.admm}
\begin{eqnarray}
	\hspace{-6ex}x_i^{t + 1}  
	&\asp{\!=\!}&  
	\argmin\limits_{x_i}~\csLm(x^{t+1|t}_{k|i}, \zz^t; \lam^t), 
	~~\,\quad i \,=\, 1, \ldots, k
	\\[-0.1cm]
	\hspace{-6ex}\zz_j^{t + 1}
	&\asp{\!=\!}&   
	\argmin\limits_{\zz_j}~\csLm(x^{t + 1},\zz_{\ell|j}^{t+1|t}; \lam^t),
	~~j \,=\, 1, \ldots, \ell
	\\[-0.1cm]
	\hspace{-6ex}\lam^{t + 1}
	&\asp{\!=\!}&  
	\lam^t \,+\, \rho\nabla_\lam\csLm(x^{t + 1}, \zz^{t + 1}; \lam^t)
\end{eqnarray}
\eseq
where $x_{k|i}^{t+1|t} \DefinedAs (x_1^{t  + 1}\dots,  x_{i  - 1}^{t + 1},  x_i, x_{i + 1}^t,  \ldots, x_{k}^t)$, $t$ is the iteration index, $\csLm$ is the augmented Lagrangian associated with~\eqref{eq.intro}, $\lam$ is the Lagrange multiplier, and $\rho$ is the step-size.

The convergence properties of~\eqref{eq.admm} are well-understood for two-block problems with $k = \ell = 1$ in~\eqref{eq.intro}~\cite{denyin16}. However, the multi-block case with either $k>1$ or $\ell>1$ is much more subtle and without imposing strong assumptions it is challenging to maintain convergence guarantees~\cite{cheheyeyua16,davyin17,latpat17} and computational convenience~\cite{ma16}. In particular, the multi-block ADMM~\eqref{eq.admm} is not (i) necessarily convergent unless additional strong convexity and rank assumptions are introduced~\cite{cheheyeyua16}; (ii) amenable to parallel implementation because the minimization in each block requires access to the solution of previous blocks~\cite{denlaipenyin17}. Although variable splitting can be used to bring the multi-block problem into the two-block setup~\cite{bertsi15}, the subproblems can become difficult to solve and the efficiency is compromised because of a significant increase in the number of variables and constraints~\cite{denlaipenyin17}. Although various modifications have been proposed for multi-block ADMM~\eqref{eq.admm} to circumvent strong assumptions that ensure convergence and computational convenience~\cite{denlaipenyin17, ma16, hehouyua15, hetaoyua12}, in contrast to standard two-block ADMM, the convergence properties of these variations remain unclear in certain scenarios. To the best of our knowledge, sufficient conditions ensuring linear convergence of these variations have not yet been established. Moreover, the empirical evidence suggests that these variations are much slower than the standard multi-block ADMM~\cite{linmazha16}.

\vsp

The primal-dual (PD) gradient flow dynamics offer a viable alternative to ADMM in terms of implementation: while ADMM requires explicit minimization, only the gradient of the Lagrangian is required to update iterates. Furthermore, in contrast to ADMM, the PD gradient flow dynamics are convenient for parallel and distributed computing even in multi-block problems without requiring any modifications relative to the two-block setup. They are thus appealing for large-scale applications and have attracted significant attention since their introduction as continuous-time dynamical systems in seminal work~\cite{arrhuruza58}.

\vsp

Recent effort centered on studying stability and convergence properties of PD gradient flow dynamics under various scenarios. Early results~\cite{feipag10,chemalcor16} focused on the asymptotic stability of the PD gradient flow dynamics that are based on the Lagrangian associated with differentiable constrained problems. Some of these results have also been extended to general saddle functions~\cite{cheghacor17, chemallowcor18, holles14} and, in a more recent effort, the focus started shifting toward proving the exponential stability~\cite{quli18, cornie19, cheli20, ozajovACC23, ozajovCDC23,dinhudhijovCDC18,dhikhojovTAC19, ozahasjovACC22, dinjovCDC20, ozajovCDC22, tanquli20} and the contraction~\cite{nguvuturslo18, cisjafbul22,davcengokrusbul23} properties. Also, advancements in Nesterov-type acceleration and design of second-order PD algorithms have been made in~\cite{attchbfadria22, zenleiche23,hehufan22a}~and~\cite{dhikhojovTAC22}, respectively. In particular,~\cite{dhikhojovTAC19} introduced a framework to bring the augmented Lagrangian associated with equality constrained convex problems into a smooth form even if the objective function contains nondifferentiable terms. This approach facilitates the use of the PD gradient flow dynamics for nonsmooth problems without resorting to the use of subgradients which complicate the analysis and substantially slow convergence. 

\vsp

In this paper, we utilize the proximal augmented Lagrangian associated with problem~\eqref{eq.intro} to introduce primal-dual gradient flow dynamics. We establish asymptotic and exponential stability of these dynamics under three sets of assumptions:
\begin{itemize}
\item Theorem~\ref{theorem.gas} requires only feasibility and convexity of the problem and establishes the global asymptotic stability.
	\item Theorem~\ref{theorem.sges} restricts the smooth blocks to convex functions satisfying Polyak-Lojasiewicz (PL) condition and the nonsmooth blocks to either polyhedral functions or group lasso penalties. These structural properties allow us to establish the semi-global exponential stability without strong convexity assumptions on the objective function or rank requirements on the constraint matrices $E$ and $F$.
	\item Theorem~\ref{theorem.ges}  removes all restrictions from the nonsmooth blocks at the expense of a range-space condition on constraint matrices and establishes the global exponential stability. While Theorem~\ref{theorem.ges} does not require presence of strongly convex terms in the objective function, the lack of strong convexity is compensated by additional column-rank conditions on the constraint matrices.
	\item Theorem~\ref{theorem.counter} proves that the range-space condition in Theorem~\ref{theorem.ges} is necessary for global exponential stability and that it cannot be relaxed without imposing additional restrictions on the nonsmooth block.  
\end{itemize}

\vsp

Under the aforementioned assumptions, the proposed dynamics have a continuous but nondifferentiable right-hand-side and admit a continuum of equilibria. This precludes the use of standard techniques including quadratic Lyapunov functions or linearization to analyze exponential stability. Even if these were applicable, the spectral analysis of a $4\times4$ block-Lyapunov or block-Jacobian matrix would be analytically challenging. To circumvent all these issues, we develop a novel Lyapunov function and establish exponential stability without imposing any regularity conditions on the equilibria.

\vsp

Theorems~\ref{theorem.sges} and \ref{theorem.ges} characterize the largest known class of problems for which the PD gradient flow dynamics exhibit exponentially fast convergence. Our assumptions are weaker than those typically required to prove linear convergence of discrete-time algorithms. Notably, unlike customized algorithms whose convergence guarantees are tailored to specific problem structure, our approach offers a systematic and broadly applicable framework for solving composite problems.

\vsp

Moreover, our approach facilitates parallel and distributed computation without requiring additional modifications. For example, application of proposed dynamics to consensus optimization problem~\eqref{ex.dist_opt} leads to a distributed implementation. In contrast to existing guarantees on the distributed algorithms, our analysis establishes exponential convergence of the distributed dynamics even in the presence of nonsmooth terms in the objective function; see Section~\ref{sec.discussion} for detailed comparison of our results with the related literature.

\vsp

The rest of the paper is organized as follows. In Section~\ref{sec.background}, we provide background material and introduce the primal-dual gradient flow dynamics. In Section~\ref{sec.summary}, we summarize our main results. In Section~\ref{sec.related_works}, we discuss related work, and compare our findings with the existing literature. In Section~\ref{sec.proofs}, we prove our main theorems; in Section~\ref{sec.numeric}, we utilize computational experiments to demonstrate the merits of our analyses; and in Section~\ref{sec.conc}, we conclude our presentation with remarks. 


\vsp

\textit{Notation:} We use $\norm{\cdot}$ and $\inner{\cdot}{\cdot}$ to denote the Euclidean norm and the standard inner product, $\sigu(A)$ and $\sigl(A)$ to denote the largest and smallest nonzero singular values of a matrix $A$, and $\cN(A)$ and $\cR(A)$ to denote the null and range spaces of $A$. We define the Euclidean distance between the vector $\pp$ and the set $\PP$ as $\dist(\pp,\PP) = \min_{\phi \, \in \, \PP}\norm{\pp - \phi}^2$.

\section{Background and motivation}\label{sec.background}

We start by providing background material and motivation for our study. In Section~\ref{sec.lag},  we introduce the Lagrangian associated with problem~\eqref{eq.intro} and derive the optimality conditions. In Section~\ref{sec.pal}, we derive a continuously differentiable saddle function and, in Section~\ref{sec.pd}, we utilize this saddle function to introduce the primal-dual (PD) gradient flow dynamics. Finally, in Section~\ref{sec.distributed}, we show that the PD gradient flow dynamics applied to the consensus optimization problem~\eqref{ex.dist_opt} results in a distributed algorithm.

	\vspace*{-1ex}
\subsection{Lagrange saddle function}
	\label{sec.lag}


Optimization problem~\eqref{eq.intro} can be lifted to a higher dimensional space by introducing auxiliary variables $\ww_i$ for each nonsmooth block associated with $\zz_i$,
	\beq\label{eq.lifted}
	\ba{rcl}
		\minimize\limits_{x,\, \zz,\, \ww}
		& \asp{} &
		f(x) \,+\, g(\ww)
		\\[\nvs]
		\subjectto 
		& \asp{} &
		Ex \,+\, F \zz  \,-\, q\,=\, 0
		\\[\nvs]
		& \asp{} &
		\zz  \,-\,  \ww\,=\, 0 
	\ea
	\eeq
where $\ww = [{\ww_1^T}{~\cdots~}{\ww_\ell^T}]^T\in\reals{n}$. The auxiliary variables isolate each nonsmooth block in the objective function and facilitate the derivation of a continuously-differentiable saddle function in Section~\ref{sec.pal}. We denote the set of all solutions to~\eqref{eq.lifted} by $\cPsw$; clearly, \mbox{$\cPsw = \set{(x,z,z)|(x,z)\in\cPs}$}. Throughout the manuscript, we use the subscript $( \, \cdot \, )_\ww$ to highlight that the solution set is associated with the lifted problem. 

\vsp

The Lagrangian associated with problem~\eqref{eq.lifted} is given by,
	\beq\label{eq.lagrange}
	\ba{r}
	\!\!\!\cL(x,\zz, \ww; \yy, \lam)  
	\, = \,
	f(x) \,+\, g(\ww) \,+\, 
	\lam^T(Ex \,+\, F\zz \,-\, q)  \,+\, \yy^T(\zz \,-\, \ww)
		\ea
	\eeq
where 
	$
	\yy  =   [{\yy_1^T}{~\cdots~}{\yy_\ell^T}]^T\in\reals{n}
	$ 
and 
	$
	\lam\in\reals{p}
	$ 
are the dual variables. Throughout the paper, we assume that there exists $(x,\zz)\in\reals{m}\times\mathbf{ri}\dom g$ such that $Ex+F\zz = q$, where $\mathbf{ri}\dom g$ denotes the relative interior of the domain of $g$~\cite[Sec.\ 6.2]{baucom11}. \blue{These assumptions ensure} that the strong duality holds~\cite[Thm.~15.23 and Prop.~15.24(x)]{baucom11}. Consequently, the necessary and sufficient conditions for $(\xs,\zs,\ws,\ys,\lams)$ to be an optimal primal-dual pair of problem~\eqref{eq.lifted} are given by the Karush-Kuhn-Tucker (KKT) conditions,
\bseq
	\label{eq.kkt}
		\begin{align}
			\nabla f(\xs) 
			& \,=\, 
			- E^T\lams
			\label{eq.opt1}
			\\
			\ys 
			& \,=\, 
			- F^T\lams
			\label{eq.opt4}
			\\
			\partial g(\ws) 
			&\,\ni\,
			\ys
			\label{eq.opt2}
			\\
			\zs 
			& \,=\, 
			\ws
			\label{eq.opt5}
			\\
			q 
			& \,=\, 
			E\xs \,+\, F\zs.
			\label{eq.opt3}
			\end{align}
	\eseq
Let $\PPsw$ denote the set of all points satisfying optimality conditions~\eqref{eq.kkt}. Since the KKT system~\eqref{eq.kkt} is challenging to solve because of nonlinear inclusions~\eqref{eq.opt1} and~\eqref{eq.opt2}, we utilize the fact that every solution $(\xs,\zs,\ws,\ys,\lams)\in\PPsw$ is a saddle point of the Lagrangian that satisfies,
	\beq\label{eq.saddle}
	\blue{\cL(\xs, \zs, \ws; \yy, \lam) \,\leq\,   \cL(\xs, \zs, \ws; \ys, \lams) \,\leq\,   \cL(x, \zz, \ww; \ys, \lams), \quad \forall x, \zz, \ww, \yy, \lam}.
	\eeq
Based on this characterization, a solution to~\eqref{eq.lifted} can be computed by simultaneous minimization and maximization of the Lagrangian over primal variables $(x,\zz,\ww)$ and dual variables $(\yy,\lam)$, respectively. In what follows, we describe how to obtain a continuously-differentiable Lagrange saddle function. We also develop primal-dual algorithms with superior performance relative to the first-order methods that utilize subgradients, which suffer from slow convergence rate even for strongly convex problems; e.g., see~\cite[Sec.~3.4]{bub15}.

	\vspace*{-1ex}
\subsection{Proximal augmented Lagrangian}\label{sec.pal}
Computation of saddle points that satisfy~\eqref{eq.saddle} is, in general, a challenging task because of the presence of nondifferentiable terms. We can alleviate these difficulties by exploiting the structure of the associated proximal operator that yields the manifold on which the augmented Lagrangian is minimized with respect to the auxiliary variable $\ww$. The augmented Lagrangian, which has the same saddle points as~\eqref{eq.lagrange}, is obtained by adding a quadratic penalty to~\eqref{eq.lagrange} for each equality constraint in~\eqref{eq.lifted} with  a penalty parameter $\mu>0$,
 	\beq\non
 		\cLm(x, \zz, \ww; \yy, \lam)  \,=\,  f(x) \,+\, g(\ww) \,+\, \lam^T(Ex \,+\, F\zz \,-\, q)  \,+\,   
 				\yy^T(\zz \,-\, \ww) \,+\, \tfrac{1}{2\mu}\norm{Ex \,+\, F\zz \,-\, q}^2
 				\,+\,   \tfrac{1}{2\mu}\norm{\zz \,-\, \ww}^2. 
 	\eeq
Completion of squares yields
 	\beq\label{eq.aug_lagrangian}
 		\cLm(x, \zz, \ww; \yy, \lam) 
 		\,=\, 
 		f(x) \,+\, g(\ww)  \,+\,  \tfrac{1}{2\mu}\norm{\ww  \,-\,  (\zz \,+\, \mu\yy) }^2 \,+\,
 		\tfrac{1}{2\mu}\norm{Ex \,+\, F\zz \,-\, q \,+\, \mu\lam}^2 \,-\, \tfrac{\mu}{2}\norm{\yy}^2 \,-\, \tfrac{\mu}{2}\norm{\lam}^2
 	\eeq
and the explicit minimizer of $\cLm$ with respect to $\ww$ is determined by the proximal operator of the function $g$,
	\beq\label{eq.wopt}
		\wb(\zz; \yy) 
		\,\DefinedAs\, 
		\argmin\limits_\ww~\cLm(x,  \zz,  \ww; \yy,\lam)
		\,=\,
		\prox_{\mu g}(z  \, + \,  \mu y).
	\eeq

For a closed proper convex function $g$ and a positive parameter $\mu$, \blue{let $L_{\mu g}(s,v) \DefinedAs g(s) + \tfrac{1}{2\mu}\norm{s - v}^2$. The Moreau envelope and the proximal operator associated with $g$ are defined, respectively, as $\cM_{\mu g}(v) = \minimize_s L_{\mu g}(s,v)$ and $\prox_{\mu g}(v) = \argmin_s L_{\mu g}(s,v)$.} Moreau envelope allows us to perform the explicit minimization of the augmented Lagrangian over $\ww$ and obtain the saddle function that is referred to as~{\em the proximal augmented Lagrangian\/}~\cite{dhikhojovTAC19},
	\beq\label{eq.pal}
	\begin{split}
		\cLm(x, \zz; \yy, \lam)  
		\,\DefinedAs\,& 
		\minimize\limits_\ww~\cLm(x, \zz, \ww; \yy, \lam)  \,=\, \cLm(x, \zz, \wb(\zz; \yy); \yy, \lam)
		\\[\nvs]
	 	 \,=\,& 
		f(x) \,+\, \cM_{\mu g}{(\zz \, + \, \mu\yy)} \,+\, \tfrac{1}{2\mu}\norm{
		Ex \, + \, F\zz \, - \, q \, + \, \mu\lam}^2 \,-\, \tfrac{\mu}{2}\norm{\yy}^2 \, - \, \tfrac{\mu}{2}\norm{\lam}^2.
	\end{split}
	\eeq
In contrast to the augmented Lagrangian which is a nonsmooth function of $\ww$, the proximal augmented Lagrangian has Lipschitz continuous gradients with respect to both primal ($x,\zz$) and dual ($\yy,\lam$) variables. This follows from the fact that the Moreau envelope is a continuously differentiable function with \blue{Lipschitz continuous gradient, $\nabla \cM_{\mu g}(v) = \tfrac{1}{\mu}(v - \prox_{\mu g}(v))$}~\cite[Prop.~12.30]{baucom11}.

\vspace*{-1ex}
\subsection{Primal-dual gradient flow dynamics}\label{sec.pd}

Since the proximal augmented Lagrangian is a continuously differentiable saddle function, first-order algorithms can be used to compute its saddle points. In particular, we utilize primal-descent dual-ascent gradient flow dynamics,\nobreak
	\bseq\label{eq.dyn}
	\begin{align}
		\xd  
		& \,=\,   
		- \nabla_x\cLm(x, \zz; \yy, \lam)  
		\label{eq.dynx}
		\\[\nvs]
		\zd
		& \,=\, 
		 - \nabla_\zz\cLm(x, \zz; \yy, \lam)  
		  \label{eq.dynz}
		 \\[\nvs]
		 \yd
		 & \,=\, 
		 \alpha\nabla_\yy\cLm(x, \zz; \yy, \lam)
		  \label{eq.dyny}
		 \\[\nvs]
		 \lamd
		 & \,=\, 
		 \alpha\nabla_\lam\cLm(x, \zz; \yy, \lam)
		 \label{eq.dynlam}
	\end{align}
	\eseq 
where $x$: $[0, \infty) \to \reals{m}$, $\zz$: $[0, \infty) \to \reals{n}$, $\yy$: $[0, \infty) \to \reals{n}$, $\lam$: $[0, \infty) \to \reals{p}$, and $\alpha$ is a positive parameter that determines the time constant of the dual dynamics. We denote the state vector in~\eqref{eq.dyn} by $\pp = (x,\zz,\yy,\lam)$. 

\vsp
	
By construction, the equilibrium points of primal-dual gradient flow dynamics~\eqref{eq.dyn} are the saddle points of the proximal augmented Lagrangian which, in conjunction with~\eqref{eq.wopt}, satisfy KKT conditions~\eqref{eq.kkt}. To show this, we set the right-hand-side of~\eqref{eq.dyn} to zero. Equation~\eqref{eq.dynlam} gives condition~\eqref{eq.opt3}. Equation~\eqref{eq.dyny} yields $\zz = \prox_{\mu g}(\zz+\mu\yy)$ which together with~\eqref{eq.wopt} implies~\eqref{eq.opt5}. Furthermore, by the definition of proximal operator, $\zz = \prox_{\mu g}(\zz+\mu\yy)$ is equivalent to $\yy\in\partial g(\zz)$ which together with~\eqref{eq.opt5} results in~\eqref{eq.opt2}. Equation~\eqref{eq.dynz} together with~\eqref{eq.dyny} and~\eqref{eq.opt3} yields~\eqref{eq.opt4}. Finally, equation~\eqref{eq.dynx} combined with~\eqref{eq.dynlam}  provides~\eqref{eq.opt1}.


\vsp

Separable structure in~\eqref{eq.intro2}  allows us to recast~\eqref{eq.dyn} in terms of individual blocks,
	\beq\non
	\ba{rclr}
	\lamd
	&\asp{=}& 
	\alpha \, \blue{h_\lambda (x,z)}&
	\\[\nvs]
	\yd_j
	&\asp{=}&
	\blue{\alpha \, h_{\yy_j} \! (y_j,z_j)}
	&~
	j \,=\, 1,\dots,\ell\phantom{.\,}
	\\[\nvs]
	\zd_j
	&\asp{=}&
	-\,\big(y_j \,+\, \frac{1}{\mu} \blue{ h_{\yy_j} \! (y_j,z_j)}\big) 
	\,-\, 
	F_j^T\big(\lam \,+\, \frac{1}{\mu} \blue{h_\lambda (x,z)}\big)
	&~
	j \,=\, 1,\dots,\ell\phantom{.\,}
	\\[\nvs]
		\xd_i  
		&\asp{=}&  
		- \, \nabla f(x_i) \,-\, E_i^T\big(\lam \,+\, \frac{1}{\mu} \blue{h_\lambda (x,z)}\big)
		&
		~
		i \,=\, 1,\dots,k
	\ea
	\eeq
where
	\beq
	\blue{\ba{rcl}
	h_\lambda (x,z)
	& \! \asp{\DefinedAs} & 
	\ds\sum_{i \,=\, 1}^{k}E_ix_i \,+\, \sum_{j \,=\, 1}^{\ell}F_j\zz_j \,-\, q
	\\[\nvs]
	\!\!
	h_{\yy_j} \! (y_j,z_j)
	& \! \asp{\DefinedAs} & 
	\zz_j  \,-\, \prox_{\mu g_j}(\zz_j \,+\, \mu\yy_j).
	\ea}
	\non
	\eeq	
Here, once $\lam$-block is computed, each $x_i$-block can be implemented in parallel to other $x$-blocks as well as to the other $\yy$- and $\zz$-blocks. Similarly, each $\yy_j$- and $\zz_j$-block can be computed in parallel to other $\yy$- and $\zz$-blocks.

	\vspace*{-1ex}
\subsection{Distributed implementation}\label{sec.distributed}
\blue{The primal-dual gradient flow dynamics are also well suited for solving distributed optimization problems.} In what follows, we show that by defining $\lam_1 \DefinedAs [T^T~0]\lam$ and $\lam_2 \DefinedAs [0~I]\lam$, we can express dynamics~\eqref{eq.dyn} associated with problem~\eqref{ex.dist_opt.equ} in a way that each update requires only local information available to the agents. Substituting these expressions for $\lam_1$ and $\lam_2$ into~\eqref{eq.dyn} and using the fact that the constraint matrices in~\eqref{ex.dist_opt.equ} are given by $E = [T^T \; C^T]^T$ and $F = [0 \; -I]^T$ yields
\beq\non
\ba{rcl}
\lamd_1
&\asp{=}&
\alpha T^TTx
\\[\nvs]
\lamd_2
&\asp{=}&
\alpha(Cx - \zz)
\\[\nvs]
\yd
&\asp{=}&
\alpha(\zz - \prox_{\mu g}(\zz  + \mu\yy))
\\[\nvs]
\zd&\asp{=}&
-\,\frac{1}{\mu}(\zz + \mu\yy - \prox_{\mu g}(\zz + \mu\yy)) + \lam_2 + \frac{1}{\mu}(Cx - \zz)
\\[\nvs]
\xd 
&\asp{=}&
-\nabla f(x)  -  \lam_1 - C^T\lam_2 - \frac{1}{\mu}(T^TTx + C^TCx - C^T\zz)
\ea
\eeq
where $T^TT$ is \blue{the Laplacian matrix of the network~\cite{mer94}. Let $N_i$ be the set} of neighbors of agent $i$. Owing to the structure of the Laplacian matrix and block diagonal form of $C$, each agent $i = 1,\ldots,k$ in the network needs to compute
\bseq\label{eq.decentralized_dyn}
\begin{align}
\lamd_{1,i}
& \,=\, 
\alpha\, h_{\lambda_{1,i}}\big(x_i, \set{x_j}_{j\in N_i}\big)
\\
\lamd_{2,i}
& \,=\, 
\alpha\, h_{\lambda_{2,i}}(x_i,z_i)
\\
\yd_i
& \,=\, 
\alpha\, h_{y_{i}}(y_i,z_i)
\\
\zd_i
& \,=\, 
-\,y_i  -  \tfrac{1}{\mu}h_{y_{i}}(y_i,z_i)  \,+\,  \lam_{2,i}  \,+\,  \tfrac{1}{\mu}h_{\lambda_{2,i}}(x_i,z_i)
\\
\xd_i 
& \,=\, 
-\nabla f_i(x_i)  \,-\,  \lam_{1,i} \,-\, C_i^T\lam_{2,i}  \,-\,\tfrac{1}{\mu}\big(h_{\lambda_{1,i}}\big(x_i, \set{x_j}_{j\in N_i}\big) + C_i^Th_{\lambda_{2,i}}(x_i,z_i)\big)
\end{align}
where 
\begin{align}
h_{\lambda_{1,i}}(x_i, \set{x_j}_{j\in N_i}) & \,\DefinedAs\, |N_i|x_i  \,-\,  \sum_{j\in N_i} x_j
\\
h_{\lambda_{2,i}}(x_i,z_i) & \,\DefinedAs\, C_ix_i  \,-\,  z_i
\\
h_{y_{i}}(y_i,z_i) & \,\DefinedAs\, z_i  \,-\,  \prox_{\mu g_i}(z_i + \mu y_i).
\end{align}
\eseq

\vsp

The forward Euler discretization of~\eqref{eq.decentralized_dyn} gives the following distributed {\em discrete-time\/} algorithm
\begin{align*}
&\lam_{1,i}^{t+1}
\,=\,
\lam_{1,i}^t \,+\, \eta\alpha\Big( |N_i|x_i^t \,-\, \ds\sum\limits_{j\,\in\,N_i}x_j^t\Big) 
\\[\nvs]
&\lam_{2,i}^{t + 1}
\,=\,
\lam_{2,i}^{t} \,+\, \eta\alpha(C_ix_i^t \,-\, \zz_i^t)\phantom{- \tfrac{1}{\alpha\mu}\left( (\yy_i^{t+1} - \yy_i^t) -(\lam_{2,i}^{t+1} - \lam_{2,i}^t)\right) }
\end{align*}
\begin{align*}
&\yy_i^{t + 1}
\,=\,
\yy_i^t \,+\, \eta\alpha(\zz_i^t \,-\, \prox_{\mu g_{i}}(\zz_i^t \,+\, \mu\yy_i^t))
\\[\nvs]
&\zz_i^{t + 1}
\,=\,
\zz_i^t \,-\, \eta(\yy_i^{t} - \lam_{2,i}^{t}) \,-\,  \tfrac{1}{\alpha\mu}\left( (\yy_i^{t+1} - \yy_i^t) -(\lam_{2,i}^{t+1} - \lam_{2,i}^t)\right) 
\\[\nvs]
&x_i^{t+1}
\,=\,
x_i^t  \,-\, \eta\nabla f_i(x_i^t)  \,-\, \eta(\lam_{1,i}^{t}  \,+\, C^T_i\lam_{2,i}^{t }) \,-\,  \tfrac{1}{\alpha\mu} \left( (\lam_{1,i}^{t + 1}  \,-\, \lam_{1,i}^{t }) \,+\, C_i^T(\lam_{2,i}^{t + 1} \,-\, \lam_{2,i}^{t })\right) 
\end{align*}
where $t\geq0$ is the iteration index, $\eta$ is a step-size, and $|N_i|$ is the cardinality of the set $N_i$. \blue{Each agent $i$ in the network requires only the optimization} variables of its neighbors, i.e., $\set{x_j}_{j\in N_i}$, to compute $\lam_{1i}^{t+1}$; furthermore, with access to $\set{x_j}_{j\in N_i}$, each node can update its own state independently of all other nodes in the network.

	\vspace*{-1ex}
\section{Main results}\label{sec.summary}
In this section, we summarize our stability results for PD gradient flow dynamics~\eqref{eq.dyn}. Our first theorem establishes Global Asymptotic Stability (GAS) under Assumption~\ref{ass.zero}, which only requires feasibility and convexity \mbox{of the problem.}

\vsp 

\begin{assumption}[constraint qualification]\label{ass.zero}
	There exists $(x,\zz)\in\reals{m}\times \mathbf{ri}\dom g$ such that $Ex+ F\zz = q$; function $f$ in problem~\eqref{eq.intro} is convex with an $L_f$-Lipschitz continuous gradient $\nabla f$; and function $g$ is proper, \mbox{closed, and convex.} 
	\end{assumption}

\vsp



\begin{theorem}[GAS]\label{theorem.gas}
Let Assumption~\ref{ass.zero} hold. The set of equilibrium points $\PPs$ of PD gradient flow dynamics~\eqref{eq.dyn}, characterized by KKT conditions~\eqref{eq.kkt}, is globally asymptotically stable and the solution to~\eqref{eq.dyn} converges to a \mbox{point in this set.}
\end{theorem}

\begin{proof}
See Section~\ref{sec.proof.theorem.gas}.
\end{proof}

\vsp

In Theorem~\ref{theorem.sges}, we establish Local Exponential Stability (LES) for a continuum of equilibria. This is done by restricting the class of functions allowed in both the smooth (Assumption~\ref{ass.eb_smooth}) and nonsmooth (Assumption~\ref{ass.eb_nonsmooth}) blocks without introducing any assumption on the constraint matrices \mbox{$E$ and $F$.}

\vsp

\begin{assumption}[Relaxation of strong convexity]\label{ass.eb_smooth}
Each smooth component in~\eqref{eq.intro2} is given by $f_i(x_i) = h_i(A_ix_i)$ for all $i=1,\dots,k$ where $x_i\in\reals{m_i}$, $h_i$: $\reals{\widebar{m}_i}\to\reals{}$ is a strongly convex function with a Lipschitz continuous gradient, and $A_i\in\reals{\widebar{m}_i\times m_i}$ is a (possibly zero) matrix.
\end{assumption}

\vsp

\begin{remark}\label{remark.str_cvx}

Since $A_i$'s in Assumption~\ref{ass.eb_smooth} are not assumed to be full-column rank, smooth block $f$ in~\eqref{eq.intro} is not necessarily a strongly convex function, but it satisfies the Polyak-Lojasiewicz (PL) condition. It is even allowed to have $f = 0$.

\end{remark}

\vsp

\begin{assumption}[Restriction of nonsmooth functions]\label{ass.eb_nonsmooth}
Each nonsmooth component in~\eqref{eq.intro2} is either (i) a polyhedral function, i.e., their epigraph can be represented as intersection of finitely many halfspaces, or (ii) a group lasso penalization, i.e., $g_j(z_j) = \eta\norm{z_j}_1 + \sum_\cJ\omega_{\cJ}\norm{z_\cJ}$, where $\zz_j\in\reals{n_j}$, $\omega_\cJ>0$, and $\cJ$ is an index \mbox{partition of $\set{1,\ldots,n_j}$.}
\end{assumption}

\vsp

\begin{remark}\label{remark.nuclear}
Functions that are frequently used in practice that satisfy Assumption~\ref{ass.eb_nonsmooth} include, but are not limited to, hinge loss, piecewise affine functions (e.g., $\ell_1$ and $\ell_\infty$ norms), indicator functions of polyhedral sets (i.e., sets associated with linear equality and inequality constraints), and $\ell_{1,2}$-norm regularization. Under additional complementary-type constraint qualifications on $\PPs$, even nuclear norm $\norm{\cdot}_\star$  can be included \mbox{in this list~\cite[Prop.~12]{zhoso17};} \blue{see Section~\ref{sec.numeric} for practical applications in which these functions arise.}
\end{remark}
\vsp

\begin{theorem}[LES]\label{theorem.sges}
Let Assumptions~\ref{ass.zero},~\ref{ass.eb_smooth}, and~\ref{ass.eb_nonsmooth} hold. There exists a time $\tb\in(0,\infty)$ such that for $t \geq \tb$ and $\alpha \in (0,\ab_1)$, any solution $\psi(t)$ to PD gradient flow dynamics~\eqref{eq.dyn} satisfy 
\beq\label{eq.theorem.sges}
\dist(\psi(t),\PPs)  \,\leq\, M_1\dist(\psi(\tb),\PPs)\rme^{-\rho_1(t \,-\, \tb)}.
\eeq
\blue{The positive constants $\tb$, $\ab_1$, $M_1$, and $\rho_1$ are defined in Section~\ref{sec.nonquad_lyap} (see Lemma~\ref{lemma.lyap2} and~\eqref{eq.M-rho1}).}
Among these constants, only $\tb$ depends on the initial distance $\dist(\psi(0),\PPs) $. 
\end{theorem}

\begin{proof}
See Section~\ref{sec.proof.semiges}.
\end{proof}

\vsp

\begin{remark}
In Theorem~\ref{theorem.sges}, the domain in which trajectories decay exponentially is characterized by the constants arising in the PL and Hoffman error-bound inequalities; see Section~\ref{sec.proof.semiges} for details. Even though these constants are independent of the initial conditions, the time $\tb$ required for the trajectory to enter this domain depends on the initial distance.
\end{remark}

\vsp

Theorem~\ref{theorem.sges} in conjunction with Theorem~\ref{theorem.gas} implies existence of an exponentially decaying global upper bound on the distance to the equilibrium points. Since this upper bound depends on the initial distance to the equilibrium points, it implies semi-global exponential stability; see~\cite[Sec.~5.10] {sas13} for the definition.

\vsp

\begin{corollary}[Semi-GES]\label{cor.sGES}
Let Assumptions~\ref{ass.zero},~\ref{ass.eb_smooth}, and~\ref{ass.eb_nonsmooth} hold. For $\alpha\in(0,\ab_1)$ where $\ab_1$ is independent of the initial condition (see~\eqref{eq.M-rho1}), any solution to PD gradient flow dynamics~\eqref{eq.dyn} is semi-globally exponentially stable, i.e., there are constants $M_\psi$ and $\rho_\psi$ depending on $\dist(\psi(0),\PPs)$ such that
\beq\label{eq.cor.sGES}
\dist(\psi(t),\PPs)  \,\leq\, M_\psi\dist(\psi(0),\PPs)\rme^{-\rho_\psi t}, ~~ t \,\geq\, 0.
\eeq
\end{corollary}

\begin{proof}
See~Section~\ref{sec.proof.cor.sges}.
\end{proof}

\vsp

In Theorem~\ref{theorem.ges}, we prove the Global Exponential Stability (GES) of PD gradient flow dynamics~\eqref{eq.dyn}.  While strong convexity of the objective function along with the invertibility of matrices $E E^T$ and $F$ in~\eqref{eq.intro} is typically required to establish GES of PD gradient flow dynamics~\cite{dhikhojovTAC19}, we identify different structural properties that allow us to relax these assumptions. Our requirements, summarized in Assumptions~\ref{ass.str_cvx}~and~\ref{ass.lip}, ensure strong convexity of the proximal augmented Lagrangian with respect to the primal variables and allow for the rows of the matrix $E$ to be linearly dependent as long as the range space of the (possibly singular) matrix $F$ is contained in the range space of~$E$. In Theorem~\ref{theorem.counter}, we further show that Assumption~\ref{ass.lip} indeed provides a necessary condition for GES.

\vsp

\begin{assumption}[str.\ convexity of aug.\ Lag.\ wrt primal var.]\label{ass.str_cvx}
	Let $I \subseteq \set{1,\dots,k}$ and  $J \subseteq \set{1, \dots, \ell}$ be the sets of indices such that for $i \in I$ and $j \in J$, functions $f_i$ and $g_j$ in~\eqref{eq.intro2} are not strongly convex. Let $E_{I}$ and $F_J$ contain the columns of matrices $E$ and $F$ associated with the blocks indexed by $I$ and $J$, respectively, and let $[{E_I}~{F_J}]$ be a full-column rank matrix.
\end{assumption}

\vsp

\begin{assumption}[range condition on $E$]\label{ass.lip}
	Constraint matrices $E$ and $F$ in~\eqref{eq.intro} satisfy
	$
	\cR(F)  \subseteq \cR(E). 
	\non
	$
\end{assumption}

\vsp

\blue{In Theorem~\ref{theorem.ges}, $m_g$ denotes} the strong convexity constant of the sum of strongly convex nonsmooth components of the objective function in~\eqref{eq.intro} and it is allowed to be zero in the absence of strongly convex \mbox{nonsmooth terms in~\eqref{eq.intro}.}

\vsp

\begin{theorem}[GES]\label{theorem.ges}
Let Assumptions~\ref{ass.zero},~\ref{ass.str_cvx}, and~\ref{ass.lip} hold and let $\alpha\in(0,\ab_2)$ and $\mu m_g \leq 1$. Any solution $\pp(t)$ to PD gradient flow dynamics~\eqref{eq.dyn} is globally exponentially stable, i.e.,
\beq\label{eq.theorem.ges}
	\norm{\pp(t) \,-\, \pas}^2  \,\leq\, M_2\norm{\pp(0) \,-\, \pas}^2 \,\rme^{-\rho_2t}, ~~ t \,\geq\, 0
\eeq
where the limit point of the trajectory, $\pas = \lim_{t\to\infty}\pp(t)$, is the orthogonal projection of $\pp(0)$ onto $\PPs$, i.e., $\pas = \argmin_{\phi\in\PPs}\norm{\pp(0) - \phi}$. \blue{The positive constants $\ab_2$, $M_2$, and $\rho_2$ are defined in Section~\ref{sec.sum_lyap} (see Lemma~\ref{lemma.conseq_ass_str} and~\eqref{eq.M-rho2}). These constants do not depend on the initial distance $\dist(\psi(0),\PPs)$.}
\end{theorem}

\begin{proof}
See Section~\ref{sec.proof.ges}.
\end{proof}

\vsp

\begin{remark}\label{remark.smooth}
Unlike~Theorem~\ref{theorem.ges}, Theorems~\ref{theorem.gas} and~\ref{theorem.sges} do not require existence of smooth components in problem~\eqref{eq.intro}, i.e., both $f$ and $E$ are allowed to be identically zero.    
\end{remark}

\vsp

\begin{remark}
In contrast to Theorem~\ref{theorem.sges}, in Theorem~\ref{theorem.ges}, we prove exponential stability for the equilibrium points that form an affine set without introducing any restrictions on the nonsmooth blocks at the expense of additional range space requirements on the constraint matrices.    
\end{remark}

\vsp

\begin{remark}
In the absence of nonsmooth blocks in problem~\eqref{eq.intro}, i.e., when both $g$ and $F$ in~\eqref{eq.intro} are absent/zero, Theorem~\ref{theorem.ges} proves the global exponential stability of PD gradient flow dynamics~\eqref{eq.dyn} for strongly convex $f$ without any additional rank assumptions on the matrix~$E$. This relaxation is especially useful in consensus problems; see Section~\ref{sec.decentralized} for details.
\end{remark}

\vsp

\begin{remark}\label{remark.alpha}
The upper bound on time constant $\alpha$ in Theorems~\ref{theorem.sges} and~\ref{theorem.ges} reflects conservatism of Lyapunov-based analysis that we utilize in our proofs. As shown in Section~\ref{sec.numeric}, dynamics~\eqref{eq.dyn} exhibit exponential convergence even for $\alpha=1$.
\end{remark}

\vsp

\begin{theorem}[Necessary cond. for GES]\label{theorem.counter}
Let Assumption~\ref{ass.zero} hold. Assumption~\ref{ass.lip} represents a necessary condition for global exponentially stability of PD gradient flow dynamics~\eqref{eq.dyn} applied to problem~\eqref{eq.intro}; this assumption cannot be relaxed without introducing additional restrictions on nonsmooth blocks.
\end{theorem}

\begin{proof}
See Section~\ref{sec.proof.counter}.
\end{proof}

	\vspace*{-1ex}
\section{Related works and discussion}\label{sec.related_works}\label{sec.discussion}
In this section, we compare and contrast our results with the existing literature; a summary table highlighting selected comparisons is provided in \hyperref[table.comp]{Table~1}.

\vsp

\subsubsection{Primal-dual gradient flow dynamics}
All problem instances studied in~\cite{quli18, cornie19, cheli20, ozajovACC23, ozajovCDC23,dinhudhijovCDC18,dhikhojovTAC19, ozahasjovACC22,dinjovCDC20, ozajovCDC22} can be cast as~\eqref{eq.intro}. We note that, Assumptions~\ref{ass.zero},~\ref{ass.str_cvx}, and~\ref{ass.lip} are much weaker than those required in these references to prove global exponential stability. In particular, the stability analyses in~\cite{quli18, cornie19, cheli20, ozajovACC23, ozajovCDC23} are limited to smooth problems with linear constraints and the main focus in~\cite{dinhudhijovCDC18,dhikhojovTAC19, ozahasjovACC22, ozajovCDC22, dinjovCDC20} is on unconstrained problems of form~\eqref{eq.unc_problem}. We note that while~\eqref{eq.unc_problem} can be brought into the form of~\eqref{eq.intro}, the converse is not possible unless $F$ is an invertible matrix. Thus, the additional challenges arising from the consensus constraint $Ex + Fz = q$ in the presence of nonsmooth components are not addressed in these existing works. 

\vsp

Moreover, Theorem~\ref{theorem.sges} exploits structural properties that allow us to establish exponentially fast convergence (i.e., semi-global exponential stability) for a class of problems for which the global exponential stability is not feasible. To the best of our knowledge, no other studies provide exponential convergence guarantees for the primal-dual gradient flow dynamics for applications studied in~\cite{taoyua11,penganwrixuma12,zarchejovgeoTAC17,simfrihastib13,pilerg20}; see Section~\ref{sec.numeric} for detailed examples.

\vsp
		 
\subsubsection{ADMM}\label{sec.diss_admm}
Our primary goal is to identify the minimal structural assumptions required for exponential convergence of PD gradient flow dynamics. In this context, existing conditions that guarantee linear convergence of ADMM variants provide valuable insights for understanding the merits of our results. However, since these methods represent distinct classes of algorithms, our intention is not to directly compare convergence of PD gradient flow dynamics to that of ADMM: the former is a gradient-based continuous-time method, whereas the latter is a discrete-time algorithm that requires explicit minimization (with respect to primal variables) of the augmented Lagrangian at every iteration.


\vsp

In~\cite[Table 1]{denyin16}, four different scenarios were provided for linear convergence of the standard two-block ADMM (i.e.,~\eqref{eq.admm} with $k=\ell=1$). The analyses in~\cite{liomer79, nislesrecpacjor15, gisboy16} fall into one of these scenarios. While Assumptions~\ref{ass.str_cvx} and~\ref{ass.lip} are satisfied in all of these four scenarios\footnote{In~\cite[Scenario 1]{denyin16}, we assume $Q\succ 0$; otherwise it is not clear how to obtain an exact solution to the nondifferentiable problem in $\yy$-update.}, our results are not restricted to the two-block case. In~\cite[Table 2]{linmazha15}, three of four scenarios considered in~\cite{denyin16} were generalized to multi-block ADMM~\eqref{eq.admm}, but the resulting conditions are much more restrictive than those introduced in Assumptions~\ref{ass.str_cvx} and~\ref{ass.lip}. For example, \cite{linmazha15} requires all $g_j$'s in problem~\eqref{eq.unc_problem} to be strongly convex and $E$ to be a full-row rank matrix; in contrast, Theorem~\ref{theorem.ges} does not impose any requirements on nondifferentiable terms for full-row rank~$E$. Lastly, in~\cite{honluo17}, the Hoffman error bound~\cite{hof52} was utilized to prove the existence of a linear convergence rate for~\eqref{eq.admm} without imposing any assumptions on the constraint matrices when the nondifferentiable components are polyhedral. In Corollary~\ref{cor.sGES}, we obtain similar results under weaker assumptions. For example, in contrast to~\cite{honluo17}, Assumptions~\ref{ass.eb_smooth} and~\ref{ass.eb_nonsmooth} do not require constraints to be compact sets and do not impose any restriction on the constraint matrices.

\vsp

Existing studies of other splitting methods~\cite{chapoc16,compes12,con13,vu13,davyin17,latpat17} target a more general class of problems than~\eqref{eq.intro} and, consequently, impose stronger assumptions when applied to~\eqref{eq.intro}. For instance, results analogous to Theorem~\ref{theorem.ges} for the Condat--Vu algorithm require smooth and strongly convex objective functions~\cite{driehrschtan24}. Similarly, guarantees comparable to Corollary~\ref{cor.sGES} for Condat--Vu~\cite{jiawucaizha22} or AFBA~\cite{latfrepat19} restrict $f$ and $g$ to be piecewise linear--quadratic functions~\cite[Ch.~10.E]{rocwet09}, which constitute a smaller function class than those satisfying Assumptions~\ref{ass.eb_smooth} and~\ref{ass.eb_nonsmooth}.

\vsp
 
\subsubsection{Comparison of related works on a consensus optimization problem}\label{sec.decentralized}
Consensus problem~\eqref{ex.dist_opt} can be used to demonstrate the utility of our analysis in the multi-block setup. For this purpose, we first examine the smooth version of the problem, 
\beq\label{eq.cons2}
\ba{rcl}
\minimize\limits_{x}
&\asp{}&
f(x)
\\[\nvs]
\subjectto 
&\asp{}&
Tx= 0
\ea
\eeq
where  $T^T$ is the incidence \blue{matrix~\cite{mer94}} of a connected undirected network, $f(x) = \sum_{i \, = \, 1}^k f_i(x_i)$, and $x = [x_1^T~\cdots~x_k^T]^T$. 

\vsp

For GES of primal-dual gradient flow dynamics, the previous results~\cite{quli18, cornie19, cheli20, ozajovACC23, ozajovCDC23,dinhudhijovCDC18,dhikhojovTAC19, ozahasjovACC22, ozajovCDC22, dinjovCDC20, tanquli20} require strong convexity of each $f_i$ and (except~\cite{ozajovACC23} and~\cite{ozajovCDC23}) assume that $T$ is a full-row rank matrix; this rank assumption on~$T$ is rarely met by the incidence matrices in practice. Moreover, although the separable structure in the multi-block problems is not exploited in the scenarios considered in~\cite{denyin16} for ADMM, a decentralized ADMM that utilizes this structure is proposed in~\cite{shilinyuawuyin14}. However, the linear convergence of the decentralized variants also requires each $f_i$ to be strongly convex~\cite{makozd17}. Finally, the decentralized gradient method EXTRA~\cite{shilinwuyin15} provides linear convergence for problem~\eqref{eq.cons2} in which at least one $f_i$ is strongly convex. It is worth noting that EXTRA can be obtained via forward Euler discretization of~\eqref{eq.dyn}; see~\cite[Sec.~IV-C]{dhikhojovTAC19}. 

\vsp

Theorem~\ref{theorem.ges} establishes GES of distributed dynamics~\eqref{eq.decentralized_dyn} applied to problem~\eqref{eq.cons2} by assuming strong convexity of only one $f_i$ without making any assumption on~$T$. This is because (i) Assumption~\ref{ass.str_cvx} is satisfied if only one $f_i$ is strongly convex as the incidence matrix of any connected undirected network becomes full-row rank when one of the rows is removed; (ii) Assumption~\ref{ass.lip} trivially holds as the condition $\cR(F) = \set{0}\subseteq\cR(E)$ is satisfied \mbox{for any~$E$.}

\vsp

Now, let us remove the smoothness assumption and study original problem~\eqref{ex.dist_opt}. None of the aforementioned works offer a convergence analysis for the associated algorithms because of the presence of both nonsmooth terms \blue{in the objective function} and consensus constraint. In~\cite{shilinwuyin15b}, a proximal variant of EXTRA that can also handle nonsmooth $g_i$'s in~\eqref{ex.dist_opt.ori} was proposed but $C_i$'s are taken to be identity matrices and only a sublinear convergence is established. In~\cite{latfrepat19}, the restriction on $C_i$'s is removed and linear convergence is obtained assuming that both $f$ and $g$ are \blue{piecewise linear--quadratic functions}. On the other hand, Corollary~\ref{cor.sGES} establishes Semi-GES of distributed dynamics~\eqref{eq.decentralized_dyn} for problem~\eqref{ex.dist_opt.equ} for a wider class of smooth and nonsmooth functions without making any rank or structural assumptions on the incidence matrix $T$ or $C_i$'s.

\begin{center}
\vspace{0.3cm}
\begin{minipage}{\linewidth}
{Table 1:} Comparison of the assumptions in the existing literature to guarantee exponential stability (con\-tin\-u\-ous\--time) or linear convergence (discrete-time) in problem~\eqref{eq.intro}. 
\label{table.comp}
\end{minipage}
\\
\vspace{0.2cm}

\renewcommand{\arraystretch}{1.4}
\begin{tabular}{|p{2.9cm}|p{2.8cm}|p{2.8cm}|p{2.9cm}|l|}
\hline 
{\it Methods} & {\it Smooth part} & {\it Nonsmooth part} & {\it Matrices} & {\it Rate} \\ 
\hline
{\textbf{Proposed} \phantom{xxxx} (Theorem~\ref{theorem.sges})}  & {PL condition} & {polyhedral/ \phantom{xxx} group lasso} & {arbitrary} & {semi-GES} \\
\hline
{\textbf{Proposed} \phantom{xxxx}  (Theorem~\ref{theorem.ges})} & {strongly convex \phantom{xx} (subset of indices \phantom{x} $I^c\subseteq\set{1,\ldots,k}$)}  & {strongly convex \phantom{xx} (subset of indices \phantom{x} $J^c\subseteq\set{1,\ldots,\ell}$)} & full column-rank\phantom{xxb} $[E_I\,\,F_J]$, \phantom{xxxxx} $\cR(F) \subseteq \cR(E)$ & {GES} \\
\hline
PD flow~\cite{quli18,cornie19} & strongly convex & N/A (smooth) & full row-rank $E$ & GES \\
\hline
PD flow \cite{ozajovACC23, ozajovCDC23} & strongly convex & N/A (smooth) & arbitrary & GES \\
\hline
PD flow~\cite{dinhudhijovCDC18,dhikhojovTAC19,dinjovCDC20,ozahasjovACC22} & strongly convex & arbitrary & invertable $F$, \phantom{xxxxx} full row-rank $E$ & GES \\
\hline
PD flow \cite{ ozajovCDC22} & strongly convex & {polyhedral/ \phantom{xxx} group lasso} & invertable $F$ \phantom{x} & semi-GES \\
\hline
2-ADMM \cite{denyin16} & strongly convex & arbitrary & full column-rank $F$, full row-rank $E$& global linear \\
\hline
Multi-ADMM~\cite{linmazha15} & arbitrary & strongly convex & full row-rank $E$ & global linear \\
\hline
Multi-ADMM \cite{honluo17} & PL condition & {polyhedral/ \phantom{xxx}  group lasso} & full column-rank \phantom{xx} $E_i$ and $F_j$ $\forall i,j$ & Q-linear \\
\hline
\mbox{Condat-Vu~\cite{jiawucaizha22}} AFBA~\cite{latfrepat19} & Piecewise-Lin\-ear\-- Quadratic functions &  Piecewise-Lin\-ear\-- Quadratic functions& arbitrary & Q-linear\\
\hline
\end{tabular}
\end{center}

\vsp

\section{Proof of main results}\label{sec.proofs}
We next prove our main results. While we use a quadratic Lyapunov function $V_1$ and LaSalle's invariance principle to prove GAS (Theorem~\ref{theorem.gas}), we employ a nonquadratic Lyapunov function $V_2$ under additional Assumptions~\ref{ass.eb_smooth} and~\ref{ass.eb_nonsmooth} to establish LES (Theorem~\ref{theorem.sges}). We show that GES can be established by augmenting $V_1$ with $V_2$ under Assumptions~\ref{ass.str_cvx} and~\ref{ass.lip} (Theorem~\ref{theorem.ges}). Finally, we provide a counter example to demonstrate that Assumption~\ref{ass.lip} is necessary for GES unless the nonsmooth blocks are restricted to a subclass of convex functions (Theorem~\ref{theorem.counter}). The proofs of the lemmas presented in this section as well as a table summarizing the notation and constants used throughout the manuscript are given in the Appendix.

	\vspace*{-1ex}
\subsection{Proof of Theorem~\ref{theorem.gas}: Global asymptotic stability} \label{sec.proof.theorem.gas}

We use LaSalle's invariance principle in conjunction with the following quadratic Lyapunov function to establish the global asymptotic stability of PD gradient flow dynamics~\eqref{eq.dyn},
	\beq\label{eq.lyap1}
	\ba{l}
		V_1(\pp)   \,=\, V_1(x, \zz; \yy, \lam)  \,=\, \tfrac{1}{2}(\alpha\norm{x - \xs}^2  \,+\, \alpha\norm{\zz - \zs}^2 \,+\, \norm{\yy - \ys}^2 \,+\, \norm{\lam - \lams}^2)
	\ea
	\eeq
where $\pps \DefinedAs (\xs, \zs, \ys, \lams) $ is an arbitrary but fixed point in the solution set $\PPs$. Clearly, $V_1$ is positive definite and radially unbounded. Lemma~\ref{lemma.lyap1} establishes a negative semi-definite upper bound on the time derivative of~$V_1$.

\vsp

\begin{lemma}\label{lemma.lyap1}
Let Assumption~\ref{ass.zero} hold. The time derivative of $V_1$ in~\eqref{eq.lyap1} along the solutions of primal-dual gradient flow dynamics~\eqref{eq.dyn} with $\alpha > 0$ satisfies 
	\beq\non
	\ba{l}
		\dot{V}_1(t)   
		\,\leq\,   
	 	- c_0
		(\norm{\nabla f(x(t))  -  \nabla f(\xs)}^2  +  \norm{\nabla_{\yy,\lam}\cLm(\pp(t)) }^2)
	 	  \\[\nvs]
	\ea
	\eeq
where $c_0 = \alpha/\max(L_f,\, \mu)$.
\end{lemma}

\vsp

\begin{proof}
	See Appendix~\ref{proof.lemma.lyap1}.
\end{proof}

\vsp

Since the upper bound on $\dot{V}_1$ can possibly be zero outside the set of equilibrium points $\PPs$, based on Lemma~\ref{lemma.lyap1} we can only certify that $\dot{V}$ is a negative semi-definite function. This implies that the set of equilibrium points $\PPs$ is stable in the sense of Lyapunov, i.e., the trajectories of~\eqref{eq.dyn} always remain bounded. We next utilize LaSalle's Invariance Principle~\cite[Thm.\ 3.4]{kha02} to establish the global asymptotic stability of the set of equilibrium points $\PPs$.

\vsp


On the set of points where the upper bound on $\dot{V}_1$ in~Lemma~\ref{lemma.lyap1} is equal to zero, primal-dual dynamics~\eqref{eq.dyn} simplify to
	\beq
	\xd  \,=\, -\, E^T\lamt, ~~ \zd  \,=\, -\,(\yt \, + \, F^T\lamt ), ~~ \yd  \,=\,  0,~~\lamd \,=\, 0
	\non
	\eeq
	where $\xt \DefinedAs x  -  \xs$, $\zt  \DefinedAs \zz - \zs$, $\yt \DefinedAs \yy - \ys$, and $\lamt \DefinedAs \lam - \lams$. Hence, the time derivative of $V_1$ for these points becomes
	\beq \non
	\dot{V}_1(x,\zz,\yy,\lam)   \,=\,  -\alpha(\inner{\lamt}{E\xt  \,+\,  F\zt}  \,-\,  \inner{\yt}{\zt}).
	\eeq
	Let $V_1(t)$ denote the value of the Lyapunov function along the solution of~\eqref{eq.dyn} at time $t$ and let
	\beq\non
	\cD \,\DefinedAs\, \set{(x,\zz,\yy,\lam)\ | \ V_1(x,\zz,\yy,\lam)  \,\leq\,  V_1 (0)}.
	\eeq
	Since $V_1$ is a radially unbounded function, its sublevel sets, and hence $\cD$, are compact. Let $\cC\subseteq \cD$ denote the set in which $\dot{V}_1(t) = 0$, i.e.,
	\beq
	\ba{l}
	\cC  \,\DefinedAs\,\big\{(x, \zz, \yy, \lam) \in\cD \,|\, \nabla f(x)  \,=\,  \nabla f(\xs),\, E\xt + F\zt \,=\, 0,\, \zt \,=\, \wt,\, \inner{\lamt}{E\xt  \,+\,  F\zt}  \,+\,\inner{\yt}{\zt} \,=\, 0\big\}
	\ea
	\non
	\eeq
and let $\Omega$ denote the largest invariant set inside $\cC$. LaSalle's Invariance Principle combined with stability of~\eqref{eq.dyn} implies the global asymptotic stability of $\Omega$. Moreover, since the proximal augmented Lagrangian has a Lipschitz continuous gradient~\cite{dhikhojovTAC19}, we can use~\cite[Lem.\ A.3]{cheghacor17} together with the stability of dynamics to conclude that the solutions converge to a point in $\Omega$. In what follows, we show that $\Omega\subseteq \PPs$.
	
\vsp	
	
	Since $\Omega$ is invariant under dynamics~\eqref{eq.dyn}, $\dot{V}_1$ remains zero in $\Omega$. Hence, we have
	\beq\non
	\ba{rcl}
	0 
	&\asp{=}&
	\frac{\rmd}{\rmd t}( \inner{\lamt}{E\xt  \,+\,  F\zt}  \,+\,\inner{\yt}{\zt})
	\\[\nvs]
	&\asp{=}& 
	\inner{\lamt}{E\xd  +  F\zd}  \,+\,\inner{\yt}{\zd}  \,=\,  -\, \norm{E^T\lamt}^2  -  \norm{\yt + F^T\lamt}^2
	\ea
	\eeq
	which implies that $E^T\lamt =  0$ and $\yt + F^T\lamt = 0$. Thus, every point $\pp = (x, \zz, \yy, \lam)  \in \Omega$ satisfies the following conditions
	\bseq
	\label{eq.set_omega}
	\begin{eqnarray}
		\nabla f(x) 
		&\asp{=}&
		\nabla f(\xs)
		\label{eq.omega1}
		\\[\nvs]
		Ex \,+\, F\zz 
		&\asp{=}&
		q
		\label{eq.omega2}
		\\[\nvs]
		\zz
		&\asp{=}& 
		\prox_{\mu g}(\zz \,+\, \mu\yy)
		\label{eq.omega3}
		\\[\nvs]
		E^T\lam
		&\asp{=}&
		E^T\lams
		\label{eq.omega4}
		\\[\nvs]
		\yy
		&\asp{=}&
		-\,F^T\lam .
		\label{eq.omega5}
	\end{eqnarray}
	\eseq
	Summing~\eqref{eq.omega1}~and~\eqref{eq.omega4}~gives~\eqref{eq.opt1}. Definition of $\prox_{\mu g}$ together with~\eqref{eq.wopt}~and~\eqref{eq.omega3} implies $y \in \partial g(\ww)$ which combined with~\eqref{eq.omega5} yields~\eqref{eq.opt2}. Equations~\eqref{eq.omega2},~\eqref{eq.omega3},~and~\eqref{eq.omega5} are the same as~\eqref{eq.opt3},~\eqref{eq.opt5},~and~\eqref{eq.opt4}, respectively. Hence, $\Omega$ is a subset of equilibrium points characterized by KKT conditions~\eqref{eq.kkt}. The globally asymptotically stability of the equilibrium points follows from the fact that $\pps = (\xs, \zs, \ys, \lams)$ in~\eqref{eq.lyap1} is an arbitrary point in $\PPs$.

	\vspace*{-1ex}
\subsection{Proof of Theorem~\ref{theorem.sges}: Local exponential stability}\label{sec.proof.semiges}
The quadratic Lyapunov function $V_1$ employed for establishing GAS does not provide any convergence rate for dynamics~\eqref{eq.dyn} since its derivative is a negative semidefinite function outside the equilibria. We introduce a novel Lyapunov function $V_2$ that is based on the associated Lagrange dual problem to obtain an exponential convergence rate. We also restrict the class of functions in~\eqref{eq.intro} and exploiting structural properties expressed in terms of local error bounds. In Section~\ref{sec.dual}, we examine  properties of Lagrange dual problem associated with~\eqref{eq.intro} and the consequences of Assumptions~\ref{ass.eb_smooth} and~\ref{ass.eb_nonsmooth}. Then, in Section~\ref{sec.nonquad_lyap}, we introduce our Lyapunov function candidate and complete the stability analysis based on structural properties obtained in Section~\ref{sec.dual}. Finally, in Section~\ref{sec.proof.cor.sges}, we show how the global asymptotic stability can be incorporated into the results obtained in Section~\ref{sec.nonquad_lyap} to establish the semi-global exponential stability of the dynamics. 

\vsp

\subsubsection{Lagrange dual problem}\label{sec.dual}
Minimizing the proximal augmented Lagrangian over primal variables $(x, \zz)$ yields the Lagrange dual function associated with the lifted problem~\eqref{eq.lifted}
	\beq\label{eq.dual}
		\dual(\yy, \lam) \! \DefinedAs\! \minimize\limits_{x,\, \zz}~\cLm(x, \zz; \yy, \lam)  \!=\! \cLm(\xb (\yy, \lam), \zb (\yy, \lam);  \yy, \lam) 
	\eeq
where $(\xb (\yy, \lam), \zb (\yy, \lam))$ denotes a solution to the following system of nonlinear equations
	\bseq\label{eq.bar_system}
	\begin{eqnarray}
		\nabla f(\xb) 
		  \,+\,  E^T\big( \lam  \,+\, \tfrac{1}{\mu}(E\xb + F\zb - q)\big)  &\!\asp{=}\!& 0
		\label{eq.bar_system1}
		\\
		\hspace{-20pt}
		\nabla_\zz\cM_{\mu g}(\zb  +  \mu\yy)
		 \,+\,  F^T\big( \lam  \,+\,  \tfrac{1}{\mu}(E\xb  +  F\zb -  q)\big) &\!\asp{=}\!& 0.
		\label{eq.bar_system2}
	\end{eqnarray}
	\eseq
We denote set of all solutions to~\eqref{eq.bar_system} at $(\yy,\lam)$ by $\cP(\yy,\lam)$. Lemma~\ref{lemma.gradient_dual} shows that even if $\cP(\yy,\lam)$ for a given $(\yy,\lam)$ is not a singleton, the dual function $d(\yy,\lam)$ has a Lipschitz continuous gradient.

\vsp
\begin{lemma}\label{lemma.gradient_dual}
	The gradient of the dual function $\dual (\yy, \lam) $,
	\beq\label{eq.lemma.grad_dual}
		\nabla\dual(\yy, \lam)  \!=\! \tbo{\!\!\!\nabla_\yy\dual(\yy, \lam)\!\!\!}{\!\!\!\nabla_\lam\dual(\yy, \lam)\!\!\!} \!=\! \tbo{\!\!\!\zb (\yy, \lam) - \prox_{\mu g}(\zb (\yy, \lam) + \mu\yy)\!\!\!\!}{\!\!\!E\xb (\yy, \lam) + F\zb (\yy, \lam) - q\!\!\!\!}  
	\eeq
	is Lipschitz continuous with modulus $\mu$, where $(\xb, \zb)$ denotes a $(\yy, \lam)$-parameterized solution to~\eqref{eq.bar_system}.
\end{lemma}

\vsp

\begin{proof}
	See Appendix~\ref{proof.lemma.gradient_dual}.
\end{proof}

\vsp

The set of optimal dual variables, denoted by $\cDs$, is determined by the set of points where $\nabla \dual(\yy,\lam) = 0$. Due to the strong duality, the set of (primal) solutions to the original problem~\eqref{eq.intro} is given by $\cPs = \cup_{(\ys,\lams)\in\cDs}\cP(\ys,\lams) = \cP(\ys,\lams)$ for any $(\ys,\lams)\in\cDs$ where the second equality follows from~\cite[Thm.\ 11.50]{rocwet09}. Moreover, the optimal value of problem~\eqref{eq.intro} is equal to the maximum value of the dual function, $\dopt \DefinedAs \maximize_{\yy,\, \lam} \dual(\yy, \lam)$. 

\vsp

In Lemma~\ref{lemma.eb_primal}, we exploit the relation between the generalized gradient map associated with the augmented Lagrangian~\eqref{eq.aug_lagrangian} and the gradient of the proximal augmented Lagrangian~\eqref{eq.pal} to establish an upper bound on the distance between the solutions of dynamics~\eqref{eq.dyn} and the manifold on which the proximal augmented Lagrangian evaluates to the dual function. To achieve this goal, we utilize a PL-type inequality~\cite{karnutsch16} for minimization of the proximal augmented Lagrangian with respect to the primal variables,  which necessitates additional Assumptions~\ref{ass.eb_smooth} and~\ref{ass.eb_nonsmooth} on the objection function.  

\vsp

\begin{lemma}\label{lemma.eb_primal}
Let Assumptions~\ref{ass.zero},~\ref{ass.eb_smooth}, and~\ref{ass.eb_nonsmooth} hold. There exist positive constants $\kappa_p$ and $\delta_p$ such that the following inequalities hold when $\norm{\nabla_{x,\zz}\cLm(x,\zz;\yy,\lam)}\leq \delta_p$,
	\bseq\label{eq.lemma.eb_primal}
	\begin{eqnarray}
\hspace{-4ex}\kappa_p\dist((x,\zz),\cP(\yy,\lam))  &\asp{\leq}&  \norm{\nabla_{x,\zz}\cLm(x,\zz;\yy,\lam)} \label{eq.lemma.eb_primal.eb}
\\
\hspace{-4ex}\tfrac{\kappa_p}{2}\dist^2((x,\zz),\cP(\yy,\lam)) &\asp{\leq}& \cLm(x,\zz;\yy,\lam) \,-\, \dual(\yy,\lam)\label{eq.lemma.eb_primal.qg}
\\
\hspace{-4ex}\cLm(x,\zz;\yy,\lam) \,-\, \dual(\yy,\lam) &\!\asp{\leq}\!& \tfrac{L_{xz}}{2\kappa_p}\norm{\nabla_{x,\zz}\cLm(x,\zz;\yy,\lam)} ^2 \label{eq.lemma.eb_primal.pl}
\end{eqnarray}
\eseq
where $L_{x\zz}$ is the Lipschitz constant of~$\nabla_{x,\zz}\cLm$.
\end{lemma}

\vsp

\begin{proof}
	See Appendix~\ref{proof.lemma.eb_primal}.
\end{proof}

\vsp

In Lemma~\ref{lemma.eb_dual}, we obtain an upper bound on the distance between the manifold on which the proximal augmented Lagrangian is equal to the dual function and the set of optimal dual variables.

\vsp

\begin{lemma}\label{lemma.eb_dual}
Let Assumptions~\ref{ass.zero},~\ref{ass.eb_smooth}, and~\ref{ass.eb_nonsmooth} hold. There exist positive constants $\kappa_d$ and $\delta_d$ such that the following inequality holds when $\norm{\nabla\dual(\yy,\lam)}\leq \delta_d$,
	\bseq\label{eq.lemma.eb_dual}
	\begin{eqnarray}
\kappa_d\dist((\yy,\lam),\cDs)  &\asp{\leq}&\norm{\nabla\dual(\yy,\lam)}  \label{eq.lemma.eb_dual.eb}
\\
\tfrac{\kappa_d}{2}\dist^2((\yy,\lam),\cDs) &\asp{\leq}&  \dopt \,-\, \dual(\yy,\lam) \label{eq.lemma.eb_dual.qg}
\\
\dopt \,-\, \dual(\yy,\lam) &\asp{\leq}& \tfrac{\mu}{2\kappa_d}\norm{\nabla\dual(\yy,\lam)} ^2 \label{eq.lemma.eb_dual.pl}.
\end{eqnarray}
\eseq
\end{lemma}

\begin{proof}
	See Appendix~\ref{proof.lemma.eb_dual}.
\end{proof}

\vsp

Lemmas~\ref{lemma.eb_primal} and~\ref{lemma.eb_dual} suggest that the sum of functions on the right-hand-side of~\eqref{eq.lemma.eb_primal.qg} and~\eqref{eq.lemma.eb_dual.qg} quantifies the distance to the equilibrium points of dynamics~\eqref{eq.dyn}. In the next section, based on this observation, we propose a nonquadratic Lyapunov function candidate.

\vsp

\subsubsection{A nonquadratic Lyapunov function}\label{sec.nonquad_lyap} \blue{We now introduce our Lyapunov function candidate, which is constructed by exploiting the error bound conditions obtained in Lemmas~\ref{lemma.eb_primal} and~\ref{lemma.eb_dual}; it serves as a distance metric to the optimal \mbox{solution set},}
\beq\label{eq.lyap2}
	V_2(\pp)   \,=\,  \cLm(x, \zz; \yy, \lam) \,-\, \dual(\yy, \lam)  \,+\,  \dopt  \,-\,  \dual(\yy, \lam) .
\eeq
Here, $\cLm(x,\zz; \yy, \lam) - \dual(\yy, \lam)$ denotes the primal gap, i.e., the distance from the trajectories to the \blue{manifold $\cP(x,y)$ on which the proximal augmented Lagrangian coincides with the Lagrange dual function, while $\dopt - \dual(\yy, \lam)$ quantifies the dual gap, i.e., the distance between this manifold and the set of optimal dual variables $\cDs$.} Since either primal or dual gap is positive out\-side the equilibria, $V_2$ is a positive definite function. \mbox{To the best} of our knowledge, apart from our recent work~\cite{ozajovCDC22}, $V_2$ has not been utilized for a Lyapunov-based analysis. One key property of $V_2$ is that it is differentiable owing to the proximal augmented Lagrangian unlike many other quantities used in the analysis of similar  optimization algorithms such as ADMM~\cite{honluo17}. 

\vsp

We start our Lyapunov-based analysis by showing that $V_2$ is upper bounded by the distance to the equilibrium points.

\vsp

\begin{lemma}\label{lemma.lyap2_upper}
	Lyapunov function $V_2$ in \eqref{eq.lyap2} satisfies
	\beq\label{eq.lemma.lyap2_upper}
	V_2(x, \zz; \yy, \lam)   \,\leq\, c_1\dist^2((x,\zz,\yy,\lam),\PPs) 
	\eeq 
	where $c_1 = (L_{xz}/2+1)\max(1,\mu)$ and $L_{x\zz}$ is the Lipschitz constant of~$\nabla_{x,\zz}\cLm$.
\end{lemma} 

\vsp

\begin{proof}
	See Appendix~\ref{proof.lemma.lyap2_upper}.
\end{proof}

\vsp

Moreover, Theorem~\ref{theorem.gas} implies that both $\norm{\nabla \cLm}$ and $\norm{\nabla\dual}$ along the solutions of dynamics~\eqref{eq.dyn} decay to zero, thus guarantees the existence of a finite time $\tb$ after which the proximity conditions in Lemmas~\ref{lemma.eb_primal} and~\ref{lemma.eb_dual} are satisfied. Lemma~\ref{lemma.lyap2} establishes a strictly negative upper bound on the time derivative of $V_2$ for $t\geq \tb$. 

\vsp

\begin{lemma}\label{lemma.lyap2}
	Let Assumptions~\ref{ass.zero},~\ref{ass.eb_smooth}, and~\ref{ass.eb_nonsmooth} hold and let $\tb\geq 0$ be such that $\norm{\nabla_{x,\zz}\cLm(x(\tb),\zz(\tb);\yy(\tb),\lam(\tb))}\leq \delta_p$ and $\norm{\nabla\dual(\yy(\tb),\lam(\tb))}\leq \delta_d$ for constants $\delta_p$ and $\delta_d$ given in Lemmas~\ref{lemma.eb_primal} and~\ref{lemma.eb_dual}, respectively. The time derivative of $V_2$ along the solutions of~\eqref{eq.dyn} with a time scale $\alpha \in (0,\ab_1)$ satisfies
	\beq\label{eq.lemma.lyap2}
		\dot{V}_2(t)   \,\leq\,   -\, \rho_1 V_2(t), \quad t\geq \tb
	\eeq
	where $L_{xz}$ is the Lipschitz constant of $\nabla_{x,\zz}\cLm$, 
	\beq\non
	\rho_1 \,=\, \min(1,2\alpha)/\max( L_{xz}/\kappa_p, \mu/\kappa_d), ~~\ab_1 \,=\, 0.5\kappa_p^2/(\sigu^2([E\ F]) +  4).
	\eeq
\end{lemma}

 \vsp

\begin{proof}
See Appendix~\ref{proof.lemma.lyap2}.
\end{proof}

\vsp

Lemma~\ref{lemma.lyap2} in conjunction with the Gronwall's inequality and Lemma~\ref{lemma.lyap2_upper} implies that for $t\geq \tb$, we have
\beq\label{eq.expo_decay}
V_2(t)  \leq V_2(\tb)\rme^{-\rho_1(t - \tb)} \leq c_1\dist^2(\pp(\tb),\PPs) \rme^{-\rho_1(t - \tb)}. 
\eeq
Substituting quadratic growth condition~\eqref{eq.lemma.eb_dual.qg} into~\eqref{eq.expo_decay} yields an exponentially decaying upper bound on the distance to the optimal dual variables for all $t\geq \tb$,
\beq\label{eq.expo_decay_dual}
\!\!\!\dist^2((\yy(t),\lam(t)),\cDs)  \leq  \tfrac{2c_1}{\kappa_d}\dist^2(\pp(\tb),\PPs) \rme^{-\rho_1(t  - \tb)}.
\eeq
Furthermore, for the distance to the optimal primal variables, Theorem~\ref{theorem.gas} implies that $(\xas,\zas)\DefinedAs\lim_{t\to\infty}(x(t),\zz(t)) \in\cPs$. Hence, using the fundamental theorem of calculus, we obtain
\begin{align}
\dist^2((x(t), \zz(t)), \cPs)  
 &\,\leq\,  \norm{(x(t), \zz(t)) \,-\,  (\xas,\zas)}^2   
\non
\\
& \,\leq\, \ds  \int_{t}^{\infty} \norm{(\xd(\tau), \zd(\tau))}^2\rmd\tau \,=\,    \int_{t}^{\infty}\norm{\nabla_{x,z}\cLm(\pp(\tau))}^2\rmd\tau
\non
\\
& \,\leq\,  \ds \int_{t}^{\infty}2L_{xz}(\cLm(\pp(\tau)) \,-\, \dual(\yy(\tau),\lam(\tau)))\rmd\tau
\non
\\
& \,\leq\,  \ds \int_{t}^{\infty}2L_{xz}V_2(\tb)\rme^{-\rho_1(\tau \,-\, \tb)}\rmd\tau 
\;=\; \tfrac{2L_{xz}}{\rho_1}V_2(\tb)\rme^{-\rho_1(t \,-\, \tb)}
\non
\\
& \,\leq\,  \tfrac{2L_{xz}}{\rho_1} \dist^2(\pp(\tb),\PPs) \rme^{-\rho_1(t \,-\, \tb)}\label{eq.expo_decay_primal}
\end{align}
where the third inequality follows from the cocoercivity of the $\nabla_{x,\zz}\cLm$~\cite[Cor.\ 18.17]{baucom11}, the fourth inequality from the fact that the dual gap is nonnegative, and the last inequality  from~\eqref{eq.expo_decay}. Combining~\eqref{eq.expo_decay_dual} with~\eqref{eq.expo_decay_primal} completes the proof with the following constants: $\ab_1 = 0.5\kappa_p^2/(\sigu^2([E\ F]) +  4)$, and
\beq\label{eq.M-rho1}
  \!\!\!M_1  =  \tfrac{2L_{xz}\max(1,L_{xz})\max(1,\mu)}{\min(\kappa_d, \rho_1)}, ~~ \rho_1  = \tfrac{\min(1,2\alpha)}{\max\left( L_{xz}/\kappa_p, \mu/\kappa_d\right) }.
\eeq

\vsp

\subsubsection{Proof of Corollary~\ref{cor.sGES}}\label{sec.proof.cor.sges}
The global asymptotic stability implies that the trajectories of dynamics~\eqref{eq.dyn} remain in the compact sublevel set $\set{\pp~|~V_1(\pp)\leq V_1(\pp(t_0))}$ where $V_1$ is a quadratic Lyapunov function used in the proof of Theorem~\ref{theorem.gas}. Using the compactness of this set, Lemmas~\ref{lemma.eb_primal} and ~\ref{lemma.eb_dual} can be improved in such a way that the local error bounds~\eqref{eq.lemma.eb_primal} and~\eqref{eq.lemma.eb_dual} hold for any time $t\geq t_0$, i.e. $\delta_p = \delta_d = \infty$, while the error constants $\kappa_p$ and $\kappa_d$ are parameterized by the initial distance $\dist(\pp(t_0), \PPs)$. In what follows, we prove this only for Lemma~\ref{lemma.eb_primal}, but the same arguments can be employed for Lemma~\ref{lemma.eb_dual}.

\vsp

Let $\cC = \set{\pp~|~V_1(\pp)\leq V_1(\pp(t_0))}$. Theorem~\ref{theorem.gas} proves that set $\cC$ is invariant under dynamics~\eqref{eq.dyn}. Moreover, from Theorem~\ref{theorem.sges}, we know that there exists a time $\tb\geq t_0$ such that the inequality $\norm{\nabla_{x,\zz}\cLm(\pp(t))}\leq\delta_p$ holds for $t\geq \tb$. However, for $t\leq \tb$, the ratio 
\beq\non
\dist((x(t), \zz(t)), \cP(\yy(t), \lam(t)))/\norm{\nabla_{x,z}\cLm(\pp(t))}
\eeq
is a continuous function~\cite[Proof of Lemma 2.3-(b)]{honluo17} and well-defined over the compact set $$\cC\cap \set{\pp~|~\norm{\nabla_{x,\zz}\cLm(\pp(t))}\geq \delta_p}.$$ Hence, it can be upper bounded by a constant $\kappa_p^\prime$ which depends on set $\cC$ and thus the initial distance $\dist(\pp(t_0), \PPs)$. 

\vspace*{-1ex}
\subsection{Proof of Theorem~\ref{theorem.ges}: Global exponential stability}\label{sec.proof.ges}
In isolation, Lyapunov functions~\eqref{eq.lyap1} and~\eqref{eq.lyap2} cannot be used to establish GES. Specifically, bounding the distance to the set of optimal dual variables is the main difficulty for establishing an exponential convergence rate. In the proof of Semi-GES, we obtain this bound in Lemmas~\ref{lemma.eb_primal} and~\ref{lemma.eb_dual} by exploiting local error bound conditions, but these conditions cannot be promoted to global guarantees unless the dual function is strongly concave. In this proof, we utilize a different set of assumptions and a pathway to obtain global results. We show that GES can be established by augmenting $V_1$ with $V_2$ under Assumptions~\ref{ass.str_cvx} and~\ref{ass.lip}. In Section~\ref{sec.dual_distance}, we use Assumption~\ref{ass.lip} and~\eqref{eq.bar_system} to substitute Lemma~\ref{lemma.eb_dual} with some global guarantees. In Section~\ref{sec.implications_ass2}, we use Assumption~\ref{ass.str_cvx} to promote local error bounds in Lemma~\ref{lemma.eb_primal} into global certificates and improve the upper bounds on $\dot{V}_1$ and $\dot{V}_2$. Finally, we integrate all findings and complete the proof in Section~\ref{sec.sum_lyap}.

\vsp

\subsubsection{Implications of Assumption~\ref{ass.lip}}\label{sec.dual_distance}
In the absence of Assumption~\ref{ass.eb_nonsmooth}, without having any additional restrictions on the nonsmooth components in problem~\eqref{eq.intro}, we cannot expect the dual function to have a particular structural property amenable for deriving bounds on the distance to the solutions. However, the connection between the dual and primal variables established in~\eqref{eq.bar_system} can be used for this purpose under additional assumptions on the constraint matrices. In particular, for arbitrary point $(\ys,\lams)\in\cDs$, we can use~\eqref{eq.bar_system1} to derive an upper bound on $\norm{\lam-\lams}$ and~\eqref{eq.bar_system2} for $\norm{\yy-\ys}$. However, we need to ensure that the difference $\lam(t) - \lams$ belongs to the range space of $E$ for all times since we can only observe  the multiplication $E^T(\lam(t) - \lams)$ through~\eqref{eq.bar_system1}. In Lemma~\ref{lemma.upper_lam}, we utilize Assumption~\ref{ass.lip} to satisfy this condition.

\vsp

\begin{lemma}\label{lemma.upper_lam}
	Let Assumptions~\ref{ass.zero} and~\ref{ass.lip} hold and let $\pas = (\xas, \zas, \yas, \lamas) \DefinedAs \lim_{t\to\infty}\pp(t)$. Then, for any $ t \geq 0$,
	\beq\non\label{eq.lemma.upper_lam}
		\norm{\lam(t) \,-\, \lamas}^2 \,\leq\, c_2(\norm{\xb(t) \,-\, \xas}^2 \,+\, \norm{\nabla_\lam\dual(\yy(t), \lam(t))}^2)
	\eeq
	where $(\xb(t), \zb(t))$ is an arbitrary point in $\cP(\yy(t), \lam(t))$ and $c_2 = \max(2L_f^2/\sigl^2(E),\,  1/\mu^{2})$. 
\end{lemma}

\vsp

\begin{proof}
	See Appendix~\ref{proof.lemma.upper_lam}.
\end{proof}

\vsp

We can derive an upper bound on $\norm{\yy - \ys}$ without needing any additional assumption or the limit argument used in Lemma~\ref{lemma.upper_lam}, as follows. For arbitrary points $(\ys, \lams)\in\cDs$ and $(\yy,\lam)$, let $(\xb, \zb)\in\cP(\yy,\lam)$. Equation~\eqref{eq.bar_system2} together with KKT condition~\eqref{eq.opt4} (i.e., $\ys =- F^T\lams$) and Lemma~\ref{lemma.gradient_dual} yields
	\begin{align}
	\norm{y - \ys}^2
	&=
	\tfrac{1}{\mu^2}\norm{\nabla_\yy\dual(\yy, \lam) + F^T(\mu(\lam - \lams) + \nabla_\lam\dual(\yy, \lam))}^2
	\non
	\\
	&\leq
	c_3(\norm{\lam - \lams}^2 + \norm{\nabla\dual(\yy, \lam)}^2)	
	\label{eq.upper_y}
	\end{align}
where $c_3 = (2/\mu^2)\max(1, \sigu^2(F), \mu^2\sigu^2(F))$. Combining Lemma~\ref{lemma.upper_lam} with~\eqref{eq.upper_y} and using the triangle inequality, $\norm{\xb-\xs}^2\leq2(\norm{\xb-x}^2+\norm{x-\xs}^2)$, we obtain that for $t\geq 0$,
\beq\label{eq.upper_dual_dist}
\ba{l}
\norm{\yy(t) \,-\, \yas}^2 \,+\, \norm{\lam(t) \,-\, \lamas}^2 
 \,\leq\, c_4( \norm{x(t) \,-\, \xb(t)}^2 \,+\,  \norm{x(t) \,-\, \xas}^2  \,+\, \norm{\nabla\dual(\yy(t),\lam(t))}^2)
\ea
\eeq
where $c_4 = 2(c_3 + 1)(c_2 + 1)$. The upper bound obtained in~\eqref{eq.upper_dual_dist} together with Lemmas~\ref{lemma.lyap1} and~\ref{lemma.lyap2} suggests that $\dot{V}_1$ and $\dot{V}_2$ can be used together to upper bound the distance to the set of optimal dual variables.

\vsp

\subsubsection{Implications of Assumption~\ref{ass.str_cvx}}\label{sec.implications_ass2}
Assumption~\ref{ass.str_cvx} has two benefits: (i) it provides a sufficient condition under which the proximal augmented Lagrangian $\cLm$ is strongly convex in primal variables $(x,\zz)$; (ii) it improves the inequality derived in Lemma~\ref{lemma.lyap1} in such a way that the left-hand-side includes a $\zz$-dependent quadratic term. Furthermore, strong convexity of $\cLm$ allows us to replace error bound condition~\eqref{eq.lemma.eb_primal.eb} in the proof of Lemma~\ref{lemma.lyap2} (upper bound on $\dot{V}_2$) and show that the time derivative of $V_2$ along the solutions of dynamics~\eqref{eq.dyn} is a negative definite function outside the equilibria, not just in certain neighborhood around equilibria. All these results are summarized in Lemma~\ref{lemma.conseq_ass_str} under an additional technical assumption that $\mu m_g\leq 1$ where $m_g$ denotes the strong convexity constant of the strongly convex nonsmooth blocks. Notably, in most problems, $m_g = 0$, hence this condition imposes no additional restriction on the selection of $\mu$.

\vsp

\begin{lemma}\label{lemma.conseq_ass_str}
Let Assumptions~\ref{ass.zero} and~\ref{ass.str_cvx} hold and let $\mu m_g \leq 1$.
\begin{enumerate}
\item[(a)] The proximal augmented Lagrangian~\eqref{eq.pal} is strongly convex in primal variables $(x,z)$ with modulus $m_{xz}$; see~\eqref{eq.mxz} for an explicit expression of  $m_{xz}$.
\item[(b)] There is a unique solution to problem~\eqref{eq.intro}, i.e., $\cPs = \set{(\xs,\zs)}$, while $\cDs$ may not be a singleton.
\item[(c)] The time derivative of quadratic Lyapunov function $V_1$ in~\eqref{eq.lyap1} along the solutions of~\eqref{eq.dyn} for any $t\geq 0$ and $\alpha>0$ satisfies
\beq\label{eq.lemma.conseq_ass_str.lyap1}
\dot{V}_1(t) \,\leq\, -\,\alpha m_{xz}\norm{(x(t), \zz(t)) \,-\, (\xs,\zs) }^2.
\eeq 
\item[(d)] The time derivative of nonquadratic Lyapunov function $V_2$ in~\eqref{eq.lyap2} along the solutions of dynamics~\eqref{eq.dyn} for any $t\geq 0$ and $\alpha\in(0,\ab_2]$ satisfies
\beq\label{eq.lemma.conseq_ass_str.lyap2}
\ba{l}
	\dot{V}_2(t)    \,\leq\,   - \min(0.5,\alpha)\big(\norm{\nabla\dual(\yy(t), \lam(t))}^2 \,+\, \norm{(x(t),\zz(t)) \,-\, (\xb(t),\zb(t))}^2\big)
\ea
\eeq
where $\{(\xb(t), \zb(t))\} = \cP(\yy(t), \lam(t))$ and $\ab_2 = 0.5m_{x\zz}^2/(\sigu^2([E~F]) +  4)$.
\end{enumerate}
\end{lemma}

\vsp

\begin{proof}
See Appendix~\ref{proof.lemma.conseq_ass_str}.
\end{proof}


\vsp

\subsubsection{Sum of two Lyapunov functions}\label{sec.sum_lyap}
To prove the global exponential stability using nonquadratic Lyapunov function $V_2$, we need to find an upper bound on~\eqref{eq.lemma.lyap2_upper} in terms of~\eqref{eq.lemma.conseq_ass_str.lyap2}. We could use the upper bound~\eqref{eq.upper_dual_dist} on the dual gap if there was not a $\norm{x-\xs}$ term in~\eqref{eq.upper_dual_dist}. \blue{This shortcoming can be remedied by augmenting $V_2$ with the quadratic Lyapunov function $V_1$, thereby creating an energy-like function that captures the coupling between primal and dual convergence. Hence, we employ the sum of $V_1$ and $V_2$, $V_3 \DefinedAs V_1 + V_2$, to establish the global exponential stability of the dynamics~\eqref{eq.dyn} as follows.}
\vsp

Since Lemma~\ref{lemma.upper_lam} provides guarantees with respect to the limit point of the trajectories $\pas \DefinedAs \lim_{t\to\infty}\pp(t)$, we fix the arbitrary reference point in the quadratic Lyapunov function to $\pas$. Theorem~\ref{theorem.gas} guarantees $\pas \in \PPs$. While Lemma~\ref{lemma.lyap2_upper} provides a quadratic upper bound on $V_3$, a quadratic lower bound is given by $V_1$ itself,
\beq\non
\xi_1\norm{\pp(t)  \,-\,  \pas}^2   \,\leq\,   V_3 (t)   \,\leq\,  \xi_2\norm{\pp(t)  \,-\,   \pas}^2
\eeq
where $\xi_1 = \alpha/2$, $\xi_2 = c_1 + 0.5$, and $c_1$ is given in Lemma~\ref{lemma.lyap2_upper}. Furthermore, combining~\eqref{eq.lemma.conseq_ass_str.lyap1} and~\eqref{eq.lemma.conseq_ass_str.lyap2} \mbox{in Lemma~\ref{lemma.conseq_ass_str} yields}
\beq\label{eq.upper_v3}
\ba{l}
\hspace{-3ex}\dot{V}_3 (t) 
  \,\leq\, 
- \min(0.5, \alpha, \alpha m_{xz})\big(\norm{\nabla\dual(\yy(t),\lam(t))}^2 \,+\,
\norm{x(t)  \,-\,  \xb(t)}^2 
\,+\, 
\norm{(x(t),\zz(t))  \,-\,  (\xas,\zas)}^2\big) .
\ea
\eeq
Moreover, upper bound~\eqref{eq.upper_dual_dist} on the distance to~$\cDs$ leads to
\beq\label{eq.upper_phi}
\ba{l}
\norm{\pp(t)  \,-\,  \pas}^2   \,\leq\,  2(c_2  \,+\,  1)(c_3  \,+\,  1)\big(\norm{x(t)  \,-\,  \xb(t)}^2 \,+\, \norm{(x(t),\zz(t)) \,-\,  (\xas,\zas)}^2    \,+\,  \norm{\nabla\dual(\yy(t),\lam(t))}^2\big)
\ea
\eeq 
where  $c_2$ and $c_3$ are given in Lemma~\ref{lemma.upper_lam} and~\eqref{eq.upper_y}, respectively. Combining~\eqref{eq.upper_v3} with~\eqref{eq.upper_phi} results in
\beq\non
\dot{V}_3 (t)  \leq  -\,\xi_3 \norm{\pp(t)   \,-\,   \pas}^2
\eeq 
where $\xi_3 \DefinedAs \min(0.5, \alpha, \alpha m_{xz})/(2(c_2 +1)(c_3+1))$. Thus, by \cite[Thm.\ 4.10]{kha02}, for any $t\geq 0$,  we have
\beq\non
\norm{\pp(t)  \,-\,  \pas}^2  \,\leq\,  \tfrac{\xi_2}{\xi_1}\norm{\pp(0)  \,-\,  \pas}^2\rme^{-\,\tfrac{\xi_3}{\xi_2}t}
\eeq
which leads to the constants: $\ab_2 =  m_{x\zz}^2/(\sigu^2([E~F]+4)$, and
\beq\label{eq.M-rho2}
M_2   =   \tfrac{\xi_2}{\xi_1}   =  \tfrac{2c_1  \,+\,  1}{\alpha}, ~~ \rho_2   =   \tfrac{\xi_3}{\xi_2}  =  \tfrac{\min(0.5, \alpha, \alpha m_{xz})}{(2c_1  \,+\,  1)(c_2 \,+\, 1)(c_3 \,+\, 1)}.	
\eeq
In what follows, we prove that $\pas$ is the orthogonal projection of $\pp(0)$ onto~$\PPs$. In Lemma~\ref{lemma.conseq_ass_str}, we show that there exists a unique solution to problem~\eqref{eq.intro} under Assumption~\ref{ass.str_cvx}. Moreover, under Assumption~\ref{ass.lip}, $\cR(F) \subseteq \cR(E)$ which together with KKT conditions~\eqref{eq.opt1} and~\eqref{eq.opt2} implies that the set of optimal dual variables $\yy$ is singleton. Thus, $\PPs$ is an affine set
\beq\non
\ba{l}
\PPs   \,=\,  \set{\xs}\times\set{\zs}\times\set{\ys} \times\set{\lam \,\in\,\reals{p} ~|~\lam   \,=\,   \lams_0  \,+\, \lam_\perp,~~\lams_\perp\,\in\,\cN([E~F]^T)}

\ea
\eeq
where $(\xs,\zs,\ys,\lams_0)$ is the unique solution in $\PPs$ with $\lams_0\in\cR([E~F])$. Moreover, derivative $\lamd(t)$ is always perpendicular to the null space component of $\PPs$, i.e., $\lamd(t)\in\cR([E~F])$. Consequently, solution $\pp(t)$ converges to the orthogonal projection of $\pp(0)$ onto $\PPs$, i.e., $\pas = \argmin_{\phi\in\PPs}\norm{\pp(0)-\phi}^2$, see~\cite{ozajovACC23} for additional discussions. 

\vsp

\subsection{Proof of Theorem~\ref{theorem.counter}: A necessary condition for GES}\label{sec.proof.counter}
We use the following academic example to prove that Assumption~\ref{ass.lip} is a necessary condition for the global exponential stability of dynamics~\eqref{eq.dyn} applied to~\eqref{eq.intro} under Assumption~\ref{ass.zero},
\beq \non\label{eq.counter1}
\ba{rl}
\minimize\limits_x &~\tfrac{1}{2} \, x^2 
\\
\subjectto  &~x  \, \leq \, 2 \text{ and } -x    \,\leq\,2.
\ea
\eeq
This problem can be brought into the form of~\eqref{eq.intro} as
\beq\label{eq.counter2}
\ba{rl}
\minimize\limits_{x,\zz} &~\tfrac{1}{2} x^2 \,+\, I_-(\zz) 
\\
\subjectto  & \tbo{-1}{1}x   \,-\, \zz  \,=\, \tbo{2}{2}
\ea
\eeq
where $f(x) = (1/2)x^2$, $g(z) = I_-(\zz)$ (indicator function of negative orthant), $E = [-1~1]^T$, $F = -I$, and $q = [2~2]^T$. Unlike Assumptions~\ref{ass.zero} and~\ref{ass.str_cvx}, Assumption~\ref{ass.lip} is not satisfied in~\eqref{eq.counter1} since $\cR(F) \subseteq\cR(E)$ holds if and only if $E$ is a full-row rank matrix. In what follows, we show that the primal-dual gradient flow dynamics~\eqref{eq.dyn} applied to~\eqref{eq.counter2} does not have global exponential stability, which implies that Assumption~\ref{ass.lip} is a necessary condition for the global exponential stability of dynamics~\eqref{eq.dyn} unless the nonsmooth block are restricted to a subclass of convex functions characterized by Assumption~\ref{ass.zero}.

\vsp 

The dynamics~\eqref{eq.dyn} applied to \eqref{eq.counter2} take the following form with a unique equilibrium point at the origin,
\beq\label{eq.dyn_counter1}
\ba{rcl}
\xd &\asp{=}& -(x  \,+\, (1/\mu)(\mu\lam \,+\,  Ex \,+\, F\zz \,-\, q))
\\[\nvs]
\zd &\asp{=}& -(1/\mu)(\zz \,+\, \mu\yy -  [\zz \,+\, \mu\yy]_-  \,+\, F^T(\mu\lam  \,+\,  Ex  \,+\,  F\zz  \,-\,  q))
\\[\nvs]
\yd &\asp{=}& \alpha(\zz  \,-\, [\zz  \,+\,  \mu\yy]_- )
\\[\nvs]
\lamd &\asp{=}& \alpha(Ex \,+\, F\zz \,-\, q).
\ea
\eeq
Our proof is based on the analytical expression of the solution to the dynamics~\eqref{eq.dyn_counter1}, but the complexity of the resulting expressions obscures the clarity of the presentation. However, since $F = -I$ in~\eqref{eq.counter2}, the proximal augmented Lagrangian~\eqref{eq.pal} contains redundant variables which can be eliminated by setting $z = Ex-q$ in the lifted problem~\eqref{eq.lifted}. This removes $\zz$ and $\lam$ variables and simplify the proximal augmented Lagrangian to
	\beq\label{eq.pal_counter}
		\cLm(x; \yy)  
		 \,=\, 
		\ds f(x) \,+\, \moreau{Ex \,-\, q \, + \, \mu\yy}  \,-\, \tfrac{\mu}{2}\norm{\yy}^2.
	\eeq
Moreover, the primal-dual gradient flow dynamics based on~\eqref{eq.pal_counter} takes the following reduced form, 
\beq\label{eq.counter_dyn2}
\ba{rcl}
\xd(t) &\asp{=}& -(x \,+\, E^T\yy  \,+\, (1/\mu)E^T(Ex \,-\, q)  \,+\,  (1/\mu)E^T[Ex \,-\, q \,+\, \mu\yy]_-)
\\[\nvs]
\yd(t) &\asp{=}& Ex \,-\, q \,-\, [Ex \,-\, q \,+\, \mu\yy]_-
\ea
\eeq
where $[\cdot]_-$ denotes the orthogonal projection onto the negative orthant; see~\cite{dhikhojovTAC19} for details. Therefore, we present our analysis for dynamics~\eqref{eq.counter_dyn2} here and note that the identical steps lead the same conclusions for~\eqref{eq.dyn_counter1}.    

\vsp

We proceed by defining the measurements $\nn = Ex-q+\mu\yy$, i.e.,
\beq\non
\nn_1  \,=\,  -x  \,-\, 2 \,+\, \mu \yy_1 \qquad \nn_2  \,=\, x  \,-\, 2 \,+\, \mu \yy_2 
\eeq
and the associated affine set 
\beq\non
\cC   \,=\,  \set{(x,\yy)\in\reals{3}\ | \ \nn_1(x,\yy)  \,\geq\, 0, \ \nn_2(x,\yy)  \,\geq\, 0 }.
\eeq
In $\cC$, the dynamics~\eqref{eq.counter_dyn2} take the following form
\beq\non\label{eq.sim_dyn}
\ba{rcl}
\xd &\asp{=}& -\left( (1 \,+\, (2/\mu))x \,-\, \yy_1 \,+\, \yy_2 \right) 
\\[\nvs]
\yd_1 &\asp{=}&  -\alpha(x \,+\, 2)  
\\[\nvs]
\yd_2 &\asp{=}& \alpha(x \,-\, 2) 
\ea
\eeq
which can be cast as a linear time invariant dynamical system
\beq\non
\ba{rcl}
\dot{\pp}(t) &\asp{=}& A\pp(t) \,+\, b
\\[\nvs]
\nn(t)v &\asp{=}& C\pp(t) \,+\, d
\ea
\eeq
where $\pp = [x^T~\yy^T]^T$ and the system matrices are given by
\beq\label{eq.counter_matrices}
A  \,=\, \thbth{-(1 \,+\, \tfrac{2}{\mu})}{1}{-1}{-\alpha}{0}{0}{\alpha}{0}{0}, \quad b  \,=\, \thbo{0}{-2\alpha}{-2\alpha}, \quad C \,=\, \tbth{-1}{\mu}{0}{1}{0}{\mu}, \quad d \,=\, \tbo{-2}{-2}.
\eeq  
We apply the coordinate transformation $\phi(t) = V^{-1}\pp(t)$ based on the eigenvalue decomposition of the state matrix, 
\beq\label{eq.counter_eig}
A = V\Lambda V^{-1}, 
~~ 
\Lambda  =  \diag(\sigma,\tfrac{2\alpha}{\sigma},0), 
~~ 
V  =  \thbth{\!\!\sigma}{-\tfrac{2\alpha}{\sigma}}{0\!\!}{\!\!-\alpha}{\alpha}{1\!\!}{\!\!\alpha}{-\alpha}{1\!\!}, 
~~  
V^{-1}  = \thbth{\!\!\tfrac{\sigma}{\sigma^2-2\alpha}}{\tfrac{1}{\sigma^2-2\alpha}}{\tfrac{-1}{\sigma^2-2\alpha}\!\!}{\!\!\tfrac{\sigma}{\sigma^2-2\alpha}}{\tfrac{\sigma^2}{2\alpha(\sigma^2-2\alpha)}}{\tfrac{-\sigma^2}{2\alpha(\sigma^2-2\alpha)}\!\!}{\!\!0}{\tfrac{1}{2}}{\tfrac{1}{2}\!\!}
\eeq
where $\sigma$ is the root of the following polynomial $1 + (2/\mu) \,+\, 2\alpha/\sigma \,+\, \sigma = 0.$
The closed form solution of linear dynamical systems yields
\beq\label{eq.counter_closed}
\phi(t)   \,=\, \rme^{\Lambda t}\phi(0) \,+\, \ds\int_{0}^{t}e^{\Lambda \tau}\rmd \tau\ V^{-1}b, \quad \nn(t) \,=\, CV\phi(t) \,+\, d.
\eeq
Substitution of~\eqref{eq.counter_matrices} and~\eqref{eq.counter_eig} into~\eqref{eq.counter_closed} gives
\beq\non
\phi(t)  \,=\, \obth{\rme^{\sigma t}\phi_1(0)}{\rme^{(2\alpha/\sigma) t}\phi_2(0)}{\phi_3(0)  \,-\, 2t}^T
\eeq
and
\beq\label{eq.counter_mea}
\ba{rcl}
\nn_1(t) &\asp{=}& -(\sigma \,+\, \alpha\mu)\rme^{\sigma t}\phi_1(0)  \,+\, \alpha(\mu \,+\, (2/\sigma))\rme^{(2\alpha/\sigma) t}\phi_2(0)  \,+\,  \mu(\phi_3(0)  \,-\, 2t) \,-\, 2
\\[\nvs]
\nn_2(t) &\asp{=}& (\sigma \,+\, \alpha\mu)\rme^{\sigma t}\phi_1(0)  \,-\, \alpha(\mu \,+\, (2/\sigma))\rme^{(2\alpha/\sigma) t}\phi_2(0)  \,+\,  \mu(\phi_3(0)  \,-\, 2t) \,-\, 2.
\ea
\eeq
The trajectory of the measurements~\eqref{eq.counter_mea} implies that if $\phi_1(0)=0$, $\phi_2(0)=0$, and $\phi_3(0) > 2/\mu$, then the dynamical system~\eqref{eq.dyn} starts in the closed set~$\cC$ and leaves it exactly at
\beq\non
t^\star  \,=\,  \tfrac{1}{2}(\phi_3(0) \,-\, \tfrac{2}{\mu}).
\eeq 
While the equilibrium point, i.e., the origin, is not in $\cC$, the time duration to leave $\cC$ depends linearly on the magnitude of the initial conditions. This contradicts with the global exponential stability. To see this, let $\mu =1$ and, let the initial conditions parameterized by a constant $\beta$  be chosen as
\beq\non
x(0)  \,=\, 0, \quad y_1(0) \,=\, 2\beta \,+\, 2, \quad y_2(0)  \,=\, 2\beta \,+\, 2.
\eeq
For $\beta >0$, $(x(0), \yy(0))\in\cC$ and the coordinate transformation~\eqref{eq.counter_eig} for any $\alpha>0$ yields
\beq\non
\phi(0)  \,=\,  \obth{0}{0}{2\beta \,+\, 2}^T, \quad \phi(t)  \,=\,  \obth{0}{0}{2\beta \,+\, 2 \,-\, 2t}^T.
\eeq
Hence, the trajectory leaves set $\cC$ exactly at $t^* = \beta$ and $\norm{\phi(t^*)} = 2$. Now assume that the system, i.e., the unique equilibrium point at the origin, is globally exponentially stable, i.e., there exist $M>0$ and $\rho>0$ (independent of $\norm{\phi(0)}$) such that
\beq\non
\norm{\phi(t)} \,\leq\, M\rme^{-\rho t}\norm{\phi(0)}.
\eeq
The exponential stability implies that $\norm{\phi(t)}\leq 2$ for $\forall t\geq \bar{t}$ where
\beq\non
 \bar{t} \,=\, \tfrac{1}{\rho}\log(\tfrac{M}{2}\norm{\phi(0)})  \,=\, \tfrac{1}{\rho}\log(M(\beta \,+\, 1)).
\eeq
However, this contradicts the fact that $t^*\leq \bar{t}$ since for any positive values of $M$ and $\rho$, the inequality 
\beq\non
\beta \,=\, t^\star  \,\leq\, \bar{t} \,=\,\tfrac{1}{\rho} \log\left( M(\beta+1)\right)    
\eeq
is violated for large values of $\beta$.

{
	\vspace*{-1ex}
\section{Examples and computational experiments}\label{sec.numeric}
In this section, we provide several examples that arise in applications to demonstrate the merits and the effectiveness of our approach for multi-block optimization problems. In our computational experiments, we set $\alpha =1$ in~\eqref{eq.dyn} to illustrate conservatism of the stability analysis on the upper bound of the time constant (see Remark~\ref{remark.alpha}). We conduct our computations on small-scale problems using Matlab's function \textsf{ode45} with relative and absolute tolerances $10^{-9}$ and $10^{-12}$, respectively.

\vsp

\vspace*{-1ex}
\subsection*{Example 1: Decentralized lasso over a network}
	\blue{To demonstrate convenience of primal-dual gradient flow dynamics~\eqref{eq.dyn} for distributed computing, let us consider consensus problem~\eqref{ex.dist_opt} over the communication network given in Fig.~\ref{fig.ex3}.} Each agent $i \in \set{1,\ldots,10}$ in the network minimizes the sum of the following functions
	\beq\label{eq.lasso}
	f_i(\xt)  \,=\,  \tfrac{1}{2}\norm{G_i\xt \,-\, h_i}^2, 
	~~
	g_i(\xt)  \,=\, \tau_i\norm{\xt}_1.
	\eeq
	Following the problem setup given in~\cite{shilinwuyin15b}, the data is generated as follows. Each measurement $h_i\in\reals{3}$ (known only by agent $i$) is constructed as $y_i = M_i\xt + \omega_i$ where entries of $M_i\in\reals{3\times100}$ and $\omega_i\in\reals{3}$ are sampled from standard normal distribution and $M_i$ is normalized to $\norm{M_i}_2 = 1$. Sparse signal $\xt\in\reals{100}$ has $5$ randomly chosen non-zero entries each of which is randomly drawn from $\set{1,\ldots,5}$. Regularization parameters $\set{\tau_i}_{i=1}^{10}$ are determined randomly to satisfy $\sum_{i=1}^{10}\tau_i = 1.15$. 
		
\vsp
	 
	Dynamics~\eqref{eq.dyn} applied to~\eqref{ex.dist_opt.equ} take the distributed form given in~\eqref{eq.decentralized_dyn} where $\nabla f_i(x_i)  \,=\, G_i^T(G_ix_i \,+\, h_i)$ and $\prox_{\mu g_{i}}(z)$ is the shrinkage operator $\cS_{1, \tau_i\mu}(z)$ whose $i$th entry is given by $[\cS_{1,\mu}(z)]_{i} =  \max\left(|z_{i}|-\mu, \ 0\right)\sign(z_{i})$. The parameter $\mu$ is set to $1$ and the initial conditions are set to $0$. 

\vsp
	
The plots of both relative state and objective function errors are given in Fig.~\ref{fig.ex3}. Since nonsmooth block $g$ satisfies Assumption~\ref{ass.eb_nonsmooth}(i), Corollary~\ref{cor.sGES} can be used to establish Semi-GES of distributed dynamics even when smooth component $f$ is not strongly convex; see Remark~\ref{remark.str_cvx}. In contrast, distributed PG-EXTRA~\cite{shilinwuyin15b} enjoys only sublinear convergence guarantees in the same setting.
  
	\begin{figure*}[h]
		\centering
		\begin{tabular}{c@{\hspace{0.1 cm}}c@{\hspace{-0.4 cm}}c@{\hspace{-0.2 cm}} c@{\hspace{-0.4 cm}}c}
			
			\begin{tabular}{c}
				\includegraphics[width=0.2\textwidth]{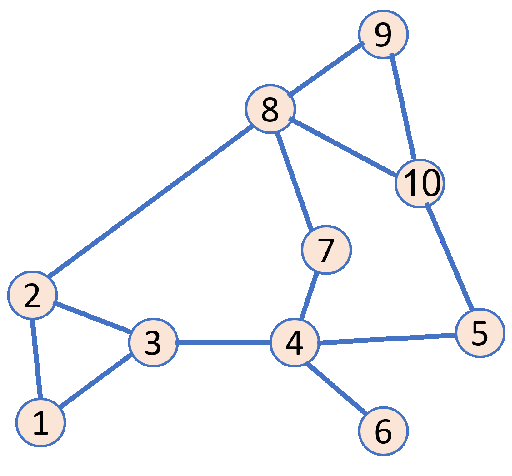}
				\\[0.2 cm]  {network topology}
			\end{tabular}
			&
			\begin{tabular}{c}
				\vspace{.25cm}
				\normalsize{\rotatebox{90}{ relative state error}}
			\end{tabular}
			&
			\begin{tabular}{c}
				\includegraphics[width=0.35\textwidth]{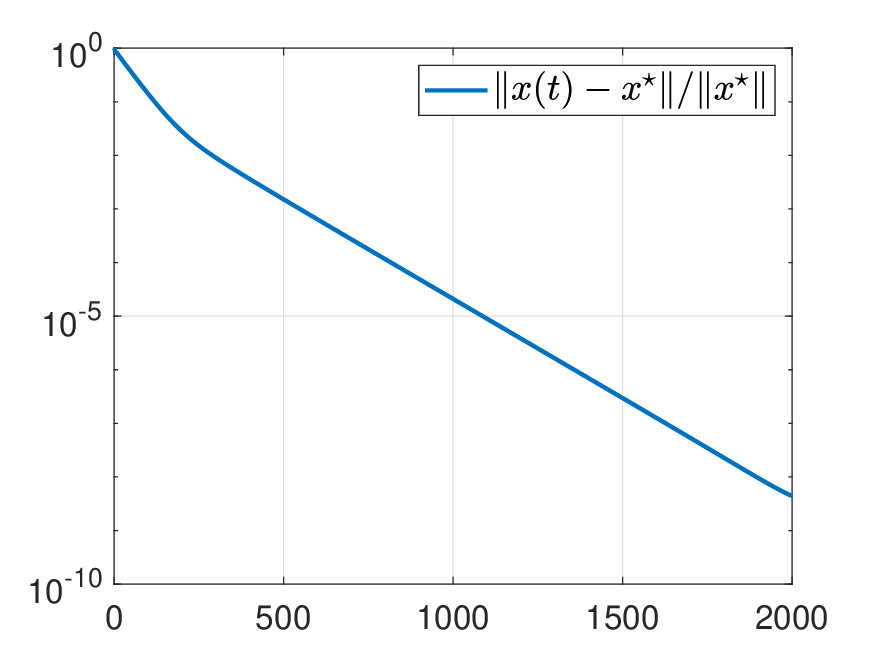}
				\\[-0.2 cm]  {time $t$}
			\end{tabular}
			&
			\begin{tabular}{c}
				\vspace{.25cm}
				\normalsize{\rotatebox{90}{ relative function error}}
			\end{tabular}
			&
			\begin{tabular}{c}
				\includegraphics[width=0.35\textwidth]{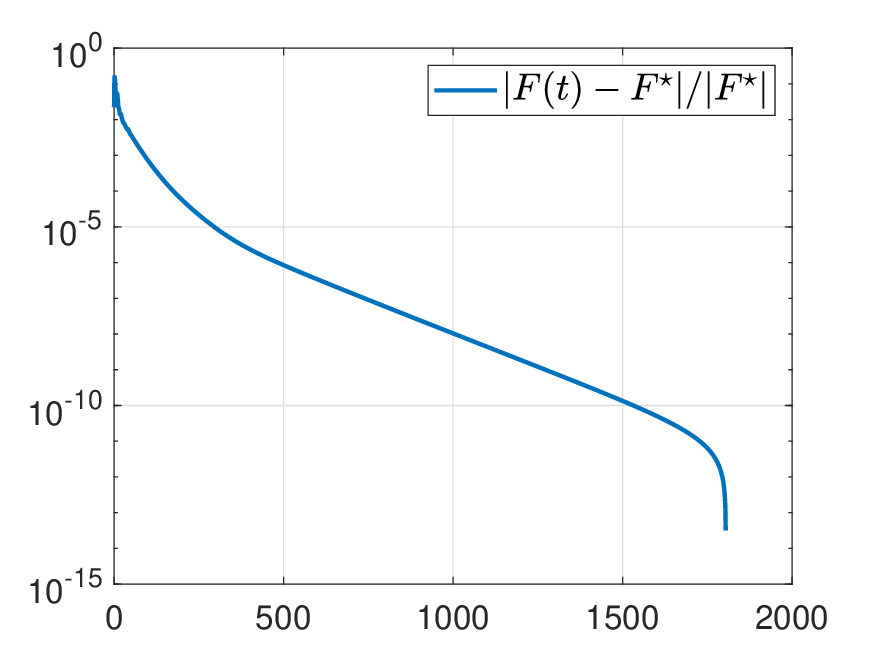}
				\\[-0.2 cm]  {time $t$}
			\end{tabular}	
		\end{tabular}
		\vspace{0cm}
		\caption{
			Topology of the underlying communication network in distributed lasso problem~\eqref{ex.dist_opt} and the Semi-GES of the distributed dynamics~\eqref{eq.decentralized_dyn}. $F(t)$ denotes the objective value of~\eqref{ex.dist_opt.equ} at time $t$. The reference solution is obtained by using CVX.}
		\label{fig.ex3}
	\end{figure*}

	\vspace*{-1ex}
	\subsection*{Example 2: Principal Component Pursuit}
	The following optimization problem arises in the recovery of low rank matrices from noisy incomplete observations~\cite{taoyua11},
\beq\label{ex.pcp}
\ba{rcl}
\minimize\limits_{Z_1,\, Z_2, \, Z_3} 
&\asp{}&
g_1(Z_1) \,+\, g_2(Z_2) \,+\, g_3(Z_3)
\\[\nvs]
\subjectto
&\asp{}&
\obth{I}{I}{I}\thbo{Z_1}{Z_2}{Z_3}  \,-~ Q\,=\, 0
\ea
\eeq 
where $g_1(Z_1) = \norm{Z_1}_\star$ is the nuclear norm, $g_2(Z_2) = \gamma_1\norm{Z_2}_1$ is $\ell_1$-norm, $g_3(Z_3) = \indicator_{\set{Z_3\,|\,\norm{P_\Omega(Z_3)}_F \leq \gamma_2}}(Z_3)$, and $P_\Omega(\cdot)$ is a binary mask that sets entries in the set $\Omega$ to zero. Similar models also arise in the estimation of sparse inverse covariance matrices~\cite{chaparwil10} and in the alignment of linearly correlated images under corruption~\cite{penganwrixuma12}.

\vsp

We generate the data for problem~\eqref{ex.pcp} using the setup given in~\cite{taoyua11} as follows. Constraint matrix $Q\in\reals{200\times 200}$ is given by $Q = Q_1 + Q_2 + Q_3$ with $Q_1 = R_1R_2^T$, where $R_1$ and $R_2$ are independent $200\times 10$ dimensional matrices whose entries are sampled from the standard normal distribution. The nonzero entries of binary mask $\Omega$ are determined at random and $80\%$ of all entries are set to $1$. The support of $Q_2$ is sampled uniformly among nonzero entries of $\Omega$ with sparsity ratio $5\%$ and the nonzero entries of $Q_2$ are uniformly sampled from interval $[-500 \ 500]$. Lastly, $Q_3$ is modeled as a white Gaussian noise with standard deviation $\sigma = 10^{-3}$. The remaining parameters are set to $\tau = 1/\sqrt{200}$ and $\delta = \sqrt{200 + \sqrt{1600}}\sigma$. 
	
\vsp	
	
	Dynamics~\eqref{eq.dyn} applied to~\eqref{ex.pcp} take the following form
	\bseq\label{eq.dyn_ex2}
	\begin{align}
			\dot{\Lambda} & \,=\,  \alpha h_{\Lambda}(Z)  \,\DefinedAs\,  \alpha(Z_1  \,+\, Z_2 \,+\, Z_3 \,-\, Q)
			\\
			\dot{Y}_j & \,=\,  \alpha h_{Y_j}(Y_j, Z_j)  \,\DefinedAs\, \alpha(Z_j  \,-\, \cS_{j,\mu_j}(Z_j  \,+\, \mu Y_j)),  	 &j \,=\, 1,2,3
			\\
			\dot{Z}_j & \,=\,  -\big(Y_j  \,+\,  \tfrac{1}{\mu} h_{Y_j}(Y_j, Z_j)\big)  \,-\,  \big(\Lambda  \,+\,  \tfrac{1}{\mu}h_{\Lambda}(Z)\big), &j \,=\, 1,2,3
		\end{align}
	\eseq
	where $(\mu_1, \mu_2, \mu_3) = (\tau\mu, \mu, \delta)$. Here, $(k,\ell)$-entry of the output of shrinkage operator $\cS_{1,\mu_1}$ is given by
	\beq\label{eq.shrinkage}
	\left[\cS_{1,\mu}(X)\right]_{k\ell}  \,=\,   \max\left(|X_{k\ell}| \,-\, \mu, 0\right)\sign(X_{k\ell}).
	\eeq
	The proximal operator $\cS_2$ of nuclear norm amounts to applying the shrinkage operator $\cS_1$ to the singular values, i.e.,
	\beq\label{eq.prox_nuclear}
	\cS_{2,\mu}(X) \,=\, U\cS_{1,\mu}(\Sigma)V^T
	\eeq
	where $X = U\Sigma V^T$ is the singular value decomposition of $X$. Lastly, the proximal operator $\cS_3$ of indicator function $g_3$ is the projection operator
	\beq
	\cS_{3,\mu}(X)  \,=\, X\circ\Omega^c  \,+\, \min(1, \, \mu/\norm{X\circ\Omega}_F) X\circ\Omega  
	\eeq
	where $\Omega_{ij}^c = 1 - \Omega_{ij}$ and $\circ$ denotes the Hadamard product. We choose zero initial conditions and, based on the formula given in~\cite{taoyua11}, we set $\mu = 1.75$. 

\vsp
	
	The plots of relative state and function errors are given in Figure~\ref{fig.ex1}. As proven in Theorem~\ref{theorem.gas}, PD gradient flow dynamics~\eqref{eq.dyn} converge globally even when the objective function does not include any smooth terms (see Remark~\ref{remark.smooth}). Furthermore, as the trajectories approach the equilibria, the convergence becomes exponential, which aligns well with Corollary~\ref{cor.sGES} (see Remark~\ref{remark.nuclear}). We note that proximal gradient methods~\cite{linganwriwuchema09} cannot be used to solve this problem because of the additional constraint on $Z_3$. The existing results on PD gradient flows such as~\cite{chemalcor16,dhikhojovTAC19,tanquli20} are not applicable in this setting because of the presence of multiple nonsmooth terms and singular constraint matrix $F = [I~I~I]$. Finally, while our approach~\eqref{eq.dyn} globally converges for arbitrary number of blocks, ADMM-based techniques, such as VASALM~\cite{taoyua11}, have convergence guarantees (without an explicit rate) only for three blocks.    
	
	\begin{figure*}[t]
		\centering
		\begin{tabular}{c@{\hspace{-0.4 cm}}c@{\hspace{0.2 cm}}c@{\hspace{-0.4 cm}}c}
			\begin{tabular}{c}
				\vspace{.25cm}
				\normalsize{\rotatebox{90}{ relative state error}}
			\end{tabular}
			&
			\begin{tabular}{c}
				\includegraphics[width=0.35\textwidth]{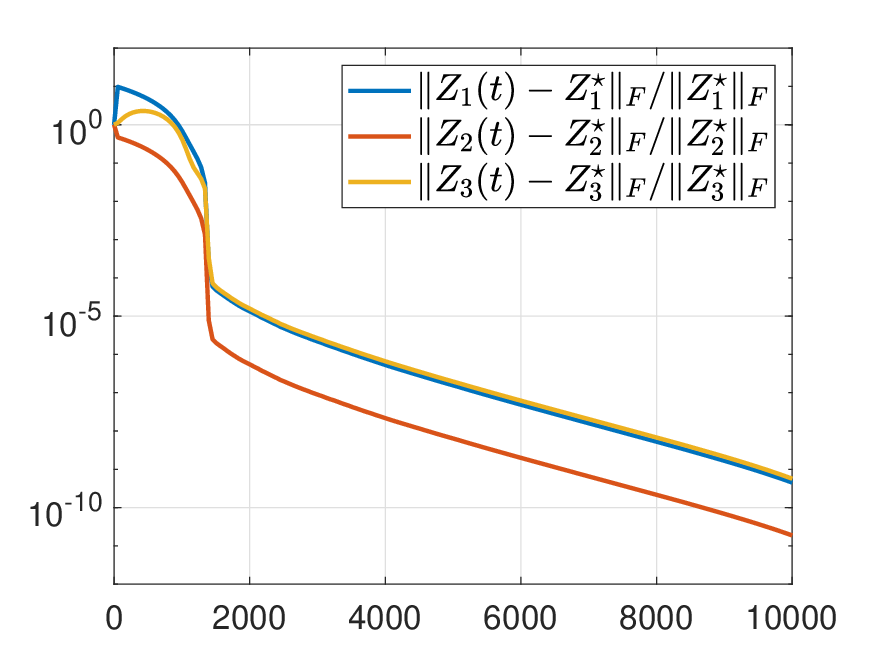}
				\\[-0.2 cm]  {time $t$}
			\end{tabular}
			&
			\begin{tabular}{c}
				\vspace{.25cm}
				\normalsize{\rotatebox{90}{ relative function error}}
			\end{tabular}
			&
			\begin{tabular}{c}
				\includegraphics[width=0.35\textwidth]{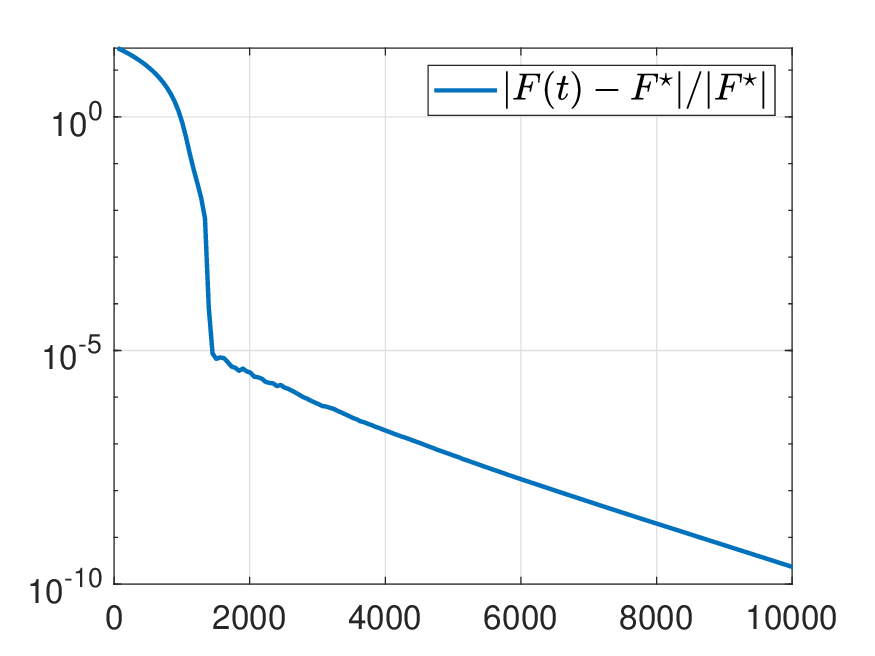}
				\\[-0.2 cm]  {time $t$}
			\end{tabular}	
		\end{tabular}
		\vspace{0cm}
		\caption{
		Semi-GES of dynamics~\eqref{eq.dyn_ex2} for principle component pursuit problem~\eqref{ex.pcp}. $F(t)$ denotes the objective value of~\eqref{ex.pcp} at time $t$. The reference solution is obtained by performing $10^4$ iterations of VASALM algorithm~\cite{taoyua11}.   }
		\label{fig.ex1}
	\end{figure*}

	\subsection*{Example 3: Covariance completion} 
	To demonstrate that our results carry over to setups where constraint matrices are replaced by bounded linear operators, we consider the optimization problem which arises in identification of statistics of disturbances to dynamical models~\cite{zarchejovgeoTAC17}, 
\beq\label{ex.cc}
\ba{rcl}
\minimize\limits_{X,\, Z} 
&\asp{}&
f(X) \,+\, g(Z)
\\[\nvs]
\subjectto
&\asp{}&
\tbo{\cE_1}{\cE_2}X   \,+\,  \tbo{I}{0}Z \,=\, \tbo{0}{Q}.
\ea
\eeq
Here, $f(X) = -\log\det(X+\delta I)$, $g(Z) = \gamma\norm{Z}_\star$ is the nuclear norm, and the linear operators are given by $\cE_1(X) = AX+XA^T$ and $\cE_2(X) = (BXB^T)\circ C$ where $\circ$ denotes the Hadamard product. We use additional regularization parameter $\delta$ to ensure that $f$ is a smooth convex function.

\vsp

	We use the mass-spring-damper example in~\cite{zarchejovgeoTAC17} to generate problem data for model~\eqref{ex.cc}. The parameters $A,B,C,Q$, and $\gamma$ are generated using the script\footnote{\url{https://www.ece.umn.edu/users/mihailo/software/ccama/run_ccama.html}} provided in~\cite{zarchejovgeoTAC17} for  $N= 40$ masses. The dynamics~\eqref{eq.dyn} applied to~\eqref{ex.cc} take the form,
	\bseq\label{eq.dyn_ex3}
	\begin{align}
	\hspace{-2ex}\dot{\Lambda}_1  & \,=\,  \alpha h_{\Lambda_1}(X,Z)  \,\DefinedAs\, \alpha(\cE_1(X)  \,+\, Z)
	\\
	\hspace{-2ex}\dot{\Lambda}_2  & \,=\,  \alpha h_{\Lambda_2}(X)  \,\DefinedAs\, \alpha(\cE_2(X)  \,-\, Q) \phantom{\hspace{15ex}} 
	\\
	\hspace{-3ex}\dot{Y} & \,=\,  \alpha h_{Y}(Y,Z)  \,\DefinedAs\, \alpha\big(Z \,-\, \cS_{2,\gamma\mu}(Z \,+\, \mu Y)\big)\label{eq.dyn_ex3y}
	\\
	\hspace{-3ex}\dot{Z} & \,=\,  -\big((Y \,+\, \tfrac{1}{\mu} h_{Y}(Y,Z)\big)  \,-\, \big(\Lambda_1 + \tfrac{1}{\mu}h_{\Lambda_1}(X,Z)\big)\label{eq.dyn_ex3z}
	\\
	\hspace{-3ex}\dot{X}  & \,=\,  (X + \delta I)^{-1}   \,-\,  \cE_1^\star\big(\Lambda_1  + \tfrac{1}{\mu}h_{\Lambda_1}(X,Z)\big)   \,-\, \cE_2^\star\big(\Lambda_2  +  \tfrac{1}{\mu}h_{\Lambda_2}(X)\big) 
	\label{eq.dyn_ex3x}
	\end{align}
	\eseq
	where the proximal operator of nuclear norm $\cS_{2,\mu}$ is given in~\eqref{eq.prox_nuclear}. We set $\mu = 1$, $\delta = 10^{-12}$ and the initial conditions are determined by the aforementioned script as follows. For $Z(0) = Y(0) = \Lambda_2(0) = I\in\reals{N\times N}$, $X(0)$ is the solution of $\cE_1(X) =- Z(0)$, and $\Lambda_1(0) = 10 \Lambda_1/\norm{\Lambda_1}_2$ where $\Lambda_1$ is the solution of the Lyapunov equation $A^T\Lambda_1+\Lambda_1A =- X(0)$. 
	
	\vsp
	
	Although existing works such as~\cite{chapoc16,compes12,con13,vu13,latpat17} provide asymptotic convergence guarantees for general splitting methods applied to this problem, they do not explicitly characterize a rate of convergence. An alternating minimization algorithm tailored for this specific problem was proposed in~\cite{zarchejovgeoTAC17}, but only sublinear convergence was established. As shown in Figure~\ref{fig.ex5}, the proposed dynamics~\eqref{eq.dyn} also exhibit sublinear convergence initially; however, as the trajectories approach the equilibria, convergence becomes exponentially fast. This behavior aligns well with Corollary~\ref{cor.sGES} (see Remark~\ref{remark.nuclear}). Convergence guarantees similar to those discussed in Remark~\ref{remark.nuclear} could also be obtained for ADMM in this setting if linear operators $\cE_1$ and $\cE_2$ had trivial null space, which holds for $\cE_1$ but not for $\cE_2$. Moreover, unlike the algorithm in~\cite{zarchejovgeoTAC17} or ADMM~\cite{honluo17}  which require minimization of Lagrangian at every iteration via incremental updates, dynamics~\eqref{eq.dyn_ex3} rely on simpler gradient-based updates and allow parallel computation of nonlinear blocks~\eqref{eq.dyn_ex3y}, \eqref{eq.dyn_ex3z}, and~\eqref{eq.dyn_ex3x}.
	
		\begin{figure*}[t]
			\centering
			\begin{tabular}{c@{\hspace{-0.4 cm}}c@{\hspace{0.2 cm}}c@{\hspace{-0.4 cm}}c}
				\begin{tabular}{c}
					\vspace{.25cm}
					\normalsize{\rotatebox{90}{ relative state error}}
				\end{tabular}
				&
				\begin{tabular}{c}
					\includegraphics[width=0.35\textwidth]{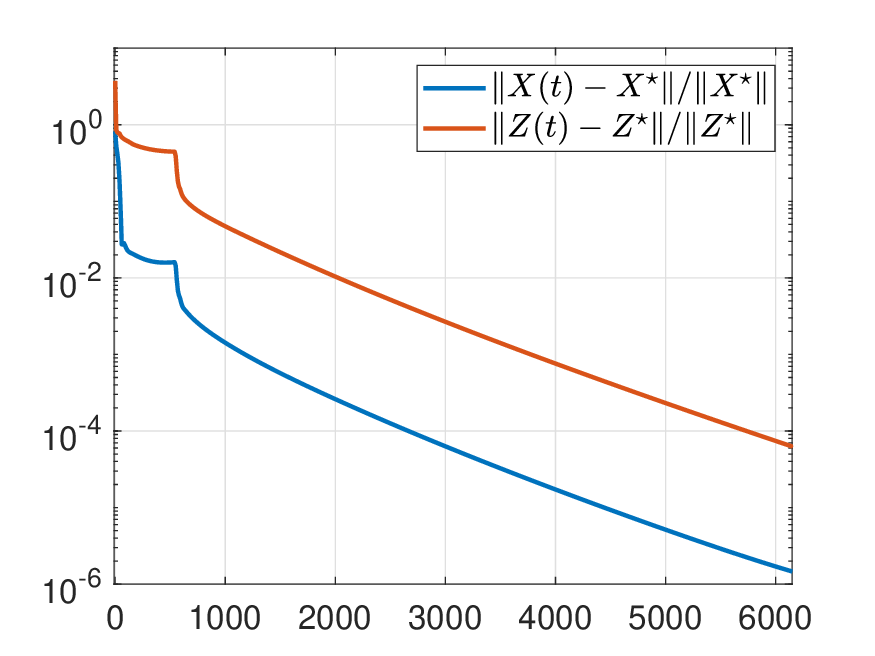}
					\\[-0.2 cm]  {time $t$}
				\end{tabular}
				&
				\begin{tabular}{c}
					\vspace{.25cm}
					\normalsize{\rotatebox{90}{ relative function error}}
				\end{tabular}
				&
				\begin{tabular}{c}
					\includegraphics[width=0.35\textwidth]{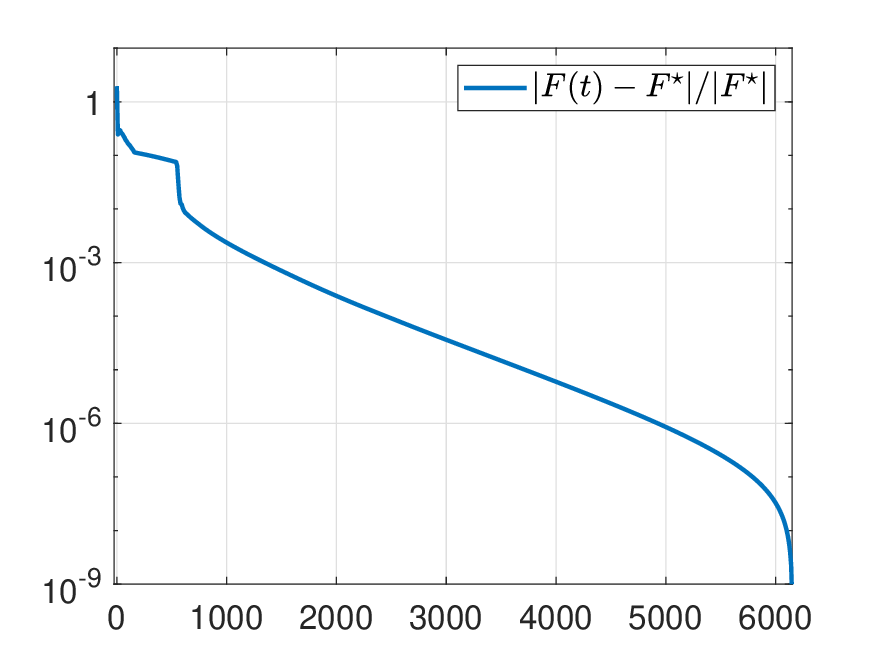}
					\\[-0.2 cm]  {time $t$}
				\end{tabular}	
			\end{tabular}
			\vspace{0cm}
			\caption{
			Semi-GES of dynamics~\eqref{eq.dyn_ex3} for covariance completion problem~\eqref{ex.cc}. $F(t)$ denotes the objective value of~\eqref{ex.cc} at time $t$. The reference solution is obtained by using CVX.   }
			\label{fig.ex5}
		\end{figure*}
		
	\vspace*{-1ex}
	
	\subsection*{Example 4: Sparse group lasso} 
	\blue{We consider the following problem to demonstrate Semi-GES of dynamics~\eqref{eq.dyn} when the nonsmooth block is not polyhedral but satisfies Assumption~\ref{ass.eb_nonsmooth}-(ii),
	\beq\label{ex.sgl}
	\ba{rcl}
	\minimize\limits_{x,\, z} &\!\asp{}\!& \frac{1}{2}\norm{x_1}_2^2  \,+\, \tau_1\norm{\zz_1}_1 \,+\, \tau_{2}\norm{\zz_2}_{1,2}
	\\[\nvs]
	\subjectto &\!\asp{}\!& \thbt{\!\!I\!}{T\!\!}{\!\!0\!}{I\!\!}{\!\!0\!}{I\!\!}\tbo{\!\!x_1\!\!}{\!\!x_2\!\!} \,+\, \thbt{\!\!0\!}{0\!\!}{\!\!-I\!}{0\!\!}{\!\!0\!}{-I\!\!}\tbo{\!\!\zz_1\!\!}{\!\!\zz_2\!\!}   \,=\,  
	\thbo{\!\!q\!\!}{\!\!0\!\!}{\!\!0\!\!}.
	\ea
	\eeq}
	We generate the data using the setup given in~\cite{simfrihastib13} as follows. The entries of $T\in\reals{60\times2000}$ are sampled from standard normal distribution and $q$ is constructed as
	$
	q = (T_1 + T_2 +  T_3)\xb + \sigma\omega
	$
	where $T = [T_1\cdots T_{50}]$ is a partition of columns, $\xb = [1\, 2\, 3\, 4\, 5\, 0\, \cdots\, 0]^T\in\reals{30}$, noise vector $\omega$ is sampled from standard normal distribution, and $\sigma$ is set so that the signal to noise ratio is $2$. Using the explicit formulas in~\cite{simfrihastib13}, the remaining parameters are chosen as $\tau_1 = 114$ and $\tau_2 = 37.94$.

	\vsp
	
	Dynamics~\eqref{eq.dyn} applied to above problem takes the following form
	\bseq\label{eq.dyn_ex4}
	\begin{eqnarray}
	\lamd &\asp{=}& \alpha h_{\lambda}(x,z) \,\DefinedAs\, \alpha\left( Ex  \,+\, F\zz  \,-\, q\right) 
	\\
	\tbo{\yd_1}{\yd_2} &\asp{=}& \alpha h_{y}(y,z) \,\DefinedAs\, \alpha \big( \tbo{\zz_1}{\zz_2} - \tbo{\cS_{1,\mu_1}(\zz_1 \,+\, \mu\yy_1)}{\cS_{4,\mu_2}(\zz_2 \,+\, \mu\yy_2)}\big) 
	\\
	\tbo{\zd_1}{\zd_2} &\asp{=}&  \,-\, \tbo{\yy_1}{\yy_2}  \,-\, \tfrac{1}{\mu_2}h_{y}(y,z) \,-\, F^T\big(\lam \,+\, \tfrac{1}{\mu} h_{\lambda}(x,z)\big)
	\\
	\tbo{\xd_1}{\xd_2} &\asp{=}&  \,-\, \tbo{x_1}{0} \,-\, E^T(\lam \,+\, \tfrac{1}{\mu} h_{\lambda}(x,z))
	\end{eqnarray}
	\eseq
	where $(\mu_1, \mu_2) = (\tau_1\mu, \tau_2\mu)$ and $\cS_1$ is the shrinkage operator given in~\eqref{eq.shrinkage}. Proximal operator $\cS_4$ is the block-shrinkage defined as
	\beq\non
	\left[\cS_{4,\mu}(\zz)\right) ]_h  \,=\, \max\left(0,\, 1 - \mu/\norm{[\zz]_h} \right) [\zz]_h
	\eeq 
	for all $h\in\cI$ where $\zz\in\reals{m}$ and $\cI$ is the partition of $\set{1,\ldots,m}$ encoded in the $\ell_{1,2}$-norm. In our case, $\cI$ is just uniform partition of $\set{1,\ldots,2000}$ to $50$ intervals each containing $40$ indices. Penalty parameter $\mu$ is taken $1$ and the initial conditions are chosen zero. The plots of relative state and function errors are given in Figure~\ref{fig.ex2}.
	
	\begin{figure*}[t]
		\centering
		\begin{tabular}{c@{\hspace{-0.4 cm}}c@{\hspace{0.2 cm}}c@{\hspace{-0.4 cm}}c}
			\begin{tabular}{c}
				\vspace{.25cm}
				\normalsize{\rotatebox{90}{ relative state error}}
			\end{tabular}
			&
			\begin{tabular}{c}
				\includegraphics[width=0.35\textwidth]{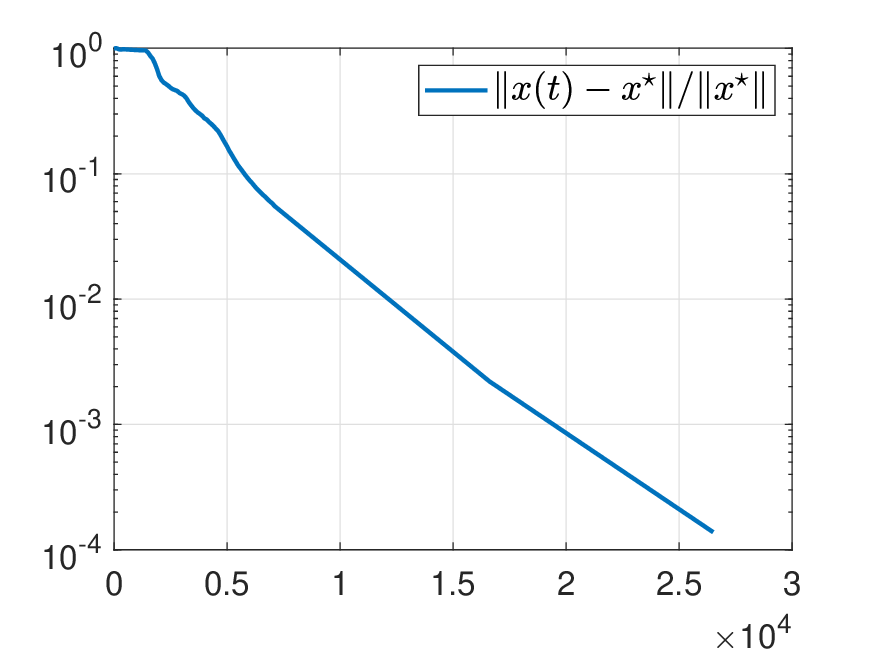}
				\\[-0.2 cm]  {time $t$}
			\end{tabular}
			&
			\begin{tabular}{c}
				\vspace{.25cm}
				\normalsize{\rotatebox{90}{ relative function error}}
			\end{tabular}
			&
			\begin{tabular}{c}
				\includegraphics[width=0.35\textwidth]{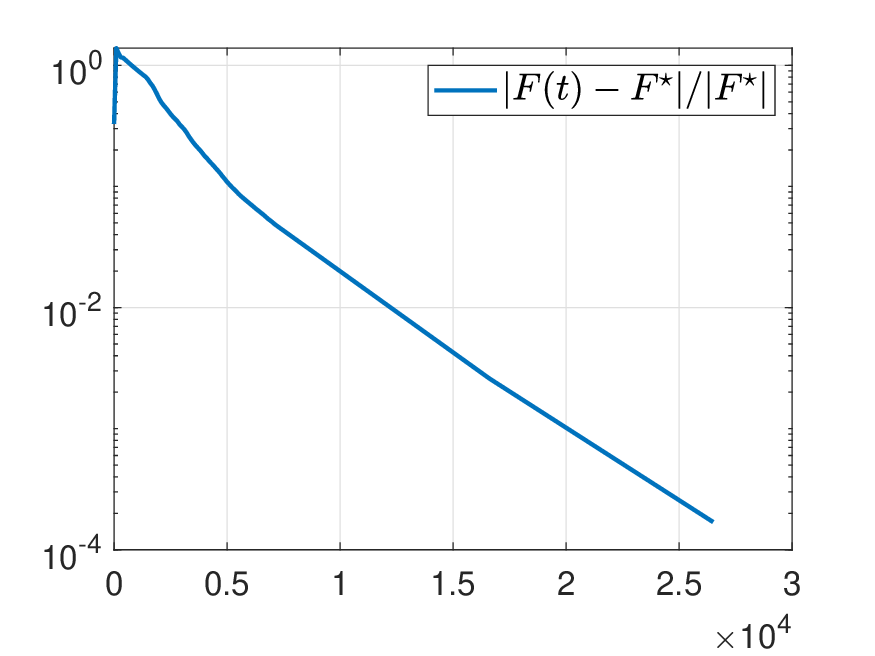}
				\\[-0.2 cm]  {time $t$}
			\end{tabular}	
		\end{tabular}
		\vspace{0cm}
		\caption{
	Semi-GES of dynamics~\eqref{eq.dyn_ex4} for sparse group lasso problem~\eqref{ex.sgl}. $F(t)$ denotes the objective value of~\eqref{ex.sgl} at time $t$. The reference solutions are obtained by using CVX.   }
		\label{fig.ex2}
	\end{figure*}

}
   
   \vspace*{-1ex}
\section{Concluding remarks}\label{sec.conc}

We have demonstrated the utility of primal-dual gradient flow dynamics for solving composite optimization problems in which a convex objective function is given by a sum of multiple, possibly nonsmooth, terms subject to the generalized consensus constraint. 
\blue{Our continuous-time framework provides a unified treatment of asymptotic and exponential convergence, leads to simple update rules, and is well suited for parallel and distributed implementation.} 
\blue{Theoretical bounds on algorithmic parameters may be conservative or difficult to compute precisely in practice, as also suggested by our numerical experiments; however, the proposed guarantees rely on assumptions that are weaker than those commonly imposed in existing works.}
\blue{Several illustrative examples, including distributed optimization problems, demonstrate the effectiveness of the proposed approach.}
\blue{Future work will focus on systematic discretizations of the proposed dynamics, tighter and more easily computable bounds, and data-driven parameter selection strategies.}

\vspace*{-1ex}

\newpage
\appendix
For clarity and brevity, we omit function arguments in the proofs whenever they can be inferred from context. Additionally, \hyperref[table.notation]{Table~2} summarizes the notation and constants used throughout the manuscript.
\begin{center}
\begin{minipage}{\linewidth}
{Table 2: Summary of the notation and the constants used thoroughout the manuscript.}
\label{table.notation}
\end{minipage}
\\
\vspace{0.2cm}
\begin{tabular}{l p{0.825\linewidth}}
\hline \hline
\textbf{Symbol} & \textbf{Description} \\ 
\hline
\multicolumn{2}{l}{\textit{Variable Modifiers}} \\
$\dot{x}$ & Time derivative of the variable $x$ \\
$x^\star$ & Optimal value or equilibrium point of variable $x$ \\
$\tilde{x}$ & Deviation from the equilibrium point, defined as $\tilde{x} = x - x^\star$ \\
$\bar{x}$ & Solution to the nonlinear system as defined in Lemma~\ref{lemma.gradient_dual} \\
$\vec{x}^\star$ & Limit point of the trajectory $x(t)$ as $t \to \infty$, i.e., $\vec{x}^\star = \lim_{t\to\infty}x(t)$ 
\\[0.3cm]
\hline
\multicolumn{2}{l}{\textit{General Operators and Spaces}} \\
$\sigu(A), \sigl(A)$ & Largest and smallest nonzero singular values of matrix $A$ \\
$\cN(A), \cR(A)$ & Null space and Range space of matrix $A$ \\
$\dist(p, \cS)$ & Euclidean distance from point $p$ to set $\cS$ \\
$\indicator_{\cS}(\cdot)$ & Indicator function of set $\cS$ \\
$\nabla f, \partial g$ & Gradient and subdifferential operators \\
$\prox_{\mu g}(\cdot)$ & Proximal operator of function $g$ with parameter $\mu$ \\
$\cM_{\mu g}(\cdot)$ & Moreau envelope of function $g$ with parameter $\mu$
\\[0.3cm]
\hline
\multicolumn{2}{l}{\textit{Problem Data and Variables}} \\
$x \in \reals{m}, \zz \in \reals{n}$ & Primal optimization variables in problem~\eqref{eq.intro}\\
$y \in \reals{n}, \lam \in \reals{p}$ & Dual optimization variables  in problem~\eqref{eq.intro}\\
$\ww \in \reals{n}$ & Auxiliary primal variable used in lifting  in problem~\eqref{eq.lifted}\\
$\pp = (x, \zz, \yy, \lam)$ & Stacked state vector of the dynamics~\eqref{eq.dyn}\\
$f_i, g_j$ & Smooth and nonsmooth function components in problem~\eqref{eq.intro}\\
$k, \ell$ & number of smooth and nonsmooth components in problem~\eqref{eq.intro} \\
$E, F, q$ & Constraint matrices and vector in problem~\eqref{eq.intro}\\
$C_i, T$ & Matrices associated with consensus problem~\eqref{ex.dist_opt}
\\[0.3cm]
\hline
\multicolumn{2}{l}{\textit{Algorithm and Analysis}} \\
$\alpha$ & Time constant for dual dynamics~\eqref{eq.dyn} \\
$\mu$ & Penalty parameter for the augmented Lagrangian~\eqref{eq.aug_lagrangian} \\
$\cL(x, \zz, \ww; \yy, \lam)$ & Standard Lagrangian of the lifted problem in~\eqref{eq.lagrange} \\
$\cLm(x, \zz; \yy, \lam)$ & Proximal augmented Lagrangian~(\eqref{eq.pal} \\
$\dual(\yy, \lam), \dopt$ & Lagrange dual function in~\eqref{eq.dual} and its maximum over $(y,\lam)$\\
$\cP(\yy, \lam)$ & Set of solutions to nonlinear equations~\eqref{eq.bar_system} for a given pair $(y,\lam)$ \\
$\cPs, \cPs_w, \cDs$ & Sets of primal, lifted primal, and dual solutions in Sections~\ref{sec.lag} and~\ref{sec.dual}\\
$\PPs$, $\PPs_w$ & Set of equilibrium points (KKT points) in problem~~\eqref{eq.intro} and~\eqref{eq.lifted}, respectively \\
$V_1, V_2, V_3$ & Lyapunov functions in~\eqref{eq.lyap1},~\eqref{eq.lyap2}, and~\eqref{eq.upper_v3}, respectively
\\[0.3cm]
\hline
\multicolumn{2}{l}{\textit{Constants}} \\
$\tb$ & Time threshold after which the trajectory enters the region of exponential attraction (Thm~\ref{theorem.sges}) \\
$\kappa_p, \delta_p$ & Error bound modulus and neighborhood radius for the primal gap (Lemma~\ref{lemma.eb_primal}) \\
$\kappa_d, \delta_d$ & Error bound modulus and neighborhood radius for the dual gap (Lemma~\ref{lemma.eb_dual}) \\
$M_1, \rho_1$ & Transient constant and decay rate for local exp. conv. for $t \geq \tb$ (Theorem~\ref{theorem.sges}) \\
$M_\psi, \rho_\psi$ & Transient constant and decay rate for Semi-GES depending on $\psi(0)$ (Corollary~\ref{cor.sGES}) \\
$M_2, \rho_2$ & Transient constant and decay rate for GES (Theorem~\ref{theorem.ges}) \\
$\ab_1, \ab_2$ & Theoretical upper bounds on the parameter $\alpha$ for stability (Theorems~\ref{theorem.sges} and~\ref{theorem.ges}) \\
$m_{xz}$ & Strong convexity modulus of $\cLm$ with respect to primal variables (Lemma~\ref{lemma.conseq_ass_str}) \\
$m_g$ & Strong convexity modulus of the nonsmooth function $g$ (Theorem~\ref{theorem.ges})
\\[0.3cm]
\hline \hline
\end{tabular}
\end{center}

\vspace*{-1ex}
\subsection{Proof of Lemma~\ref{lemma.lyap1}}\label{proof.lemma.lyap1}
\begin{mylemma}{1}
Let Assumption~\ref{ass.zero} hold. The time derivative of $V_1$ in~\eqref{eq.lyap1} along the solutions of primal-dual gradient flow dynamics~\eqref{eq.dyn} with $\alpha > 0$ satisfies 
	\beq\non
		\dot{V}_1(t)   
		\,\leq\,   
	 	- \tfrac{\alpha}{\max(L_f,\, \mu)}
		\big(\norm{\nabla f(x(t))  \, - \,  \nabla f(\xs)}^2
	 	  \,+\, \norm{\nabla_{\yy}\cLm(x(t), \zz(t); \yy(t), \lam(t)) }^2  \,+\, \norm{\nabla_{\lam}\cLm(x(t), \zz(t); \yy(t), \lam(t)) }^2\big)
	\eeq
where $L_f$ is the Lipschitz constant of $\nabla f$.
\end{mylemma}

\vsp

\begin{proof}
Let $\xt \DefinedAs x  -  \xs$, $\zt  \DefinedAs \zz - \zs$, $\yt \DefinedAs \yy - \ys$, $\lamt \DefinedAs \lam - \lams$, and $\wt  \DefinedAs \prox_{\mu g}(\zz  +  \mu\yy) - \prox_{\mu g}(\zs  + \mu\ys)$, using the chain rule, the time derivative of $V_1$ can be written as
	\begin{align}
		\dot{V}_1
		&\,=\,
		\alpha\inner{\xd}{\xt}   \,+\, \alpha\inner{\zd}{\zt} \,+\, \inner{\yd}{\yt} \,+\, \inner{\lamd}{\lamt}
		\non
		\\
		&\,=\,
		-\, \alpha\inner{\nabla f(x) - \nabla f(\xs) \,+\, E^T(\lamt \,+\, \tfrac{1}{\alpha\mu}\lamd)}{\xt}
		\,-\, 
		\alpha\inner{\yt \,+\, \tfrac{1}{\alpha\mu}\yd \,+\, F^T(\lamt \,+\, \tfrac{1}{\alpha\mu}\lamd}{\zt} \,+\,
		\non
		\\
		& \, \phantom{=} \, 
		\phantom{-} \,
		\alpha\inner{ \zt  \,-\, \wt}{\yt}
		\,+\, \alpha\inner{ E\xt \,+\, F\zt}{\lamt}
		\non
		\\
		&\,=\, 
		-\, \alpha\inner{\nabla f(x) - \nabla f(\xs)}{\xt} \,-\, \tfrac{1}{\mu}\inner{\lamd}{E\xt \,+\, F\zt} \,-\, \tfrac{1}{\mu}\inner{\yd}{\zt} \,-\, \alpha\inner{\wt}{\yt}
		\non
		\\
		&\,=\,
		-\, \alpha\inner{\nabla f(x) - \nabla f(\xs)}{\xt} \,-\, \tfrac{\alpha}{\mu}\norm{E\xt \,+\, F\zt}^2 \,-\, \tfrac{\alpha}{\mu}(\inner{\zt - \wt}{\zt} \,+\, \inner{\wt}{\mu\yt})
		\non
		\\
		&\,=\, 
		-\, \alpha\inner{\nabla f(x) - \nabla f(\xs)}{\xt} \,-\, \tfrac{\alpha}{\mu}\norm{E\xt \,+\, F\zt}^2 \,-\, \tfrac{\alpha}{\mu}(\norm{\zt}^2  \,-\, 2\inner{\wt}{\zt} \,+\, \inner{\wt}{\zt \,+\, \mu\yt}) 
	\label{eq.lyap0_derivative}
	\end{align}
	where the second equality follows from the fact that dynamics~\eqref{eq.dyn} are zero at any solutions satisfying~\eqref{eq.kkt}, the third equality follows from the symmetry of inner products, the forth equality follows from~\eqref{eq.dyny}~and~\eqref{eq.dynlam}, and the last equality follows from the linearity of inner products. Using firm-nonexpansiveness of the proximal operator~\cite[Prop.~4.2(iv)]{baucom11}, we obtain $\norm{\wt}^2\leq\inner{\wt}{\zt+\yt}$. Substitution of this bound into~\eqref{eq.lyap0_derivative} yields
	\beq
	\ba{rcl}
	\dot{V}_1
	&\asp{\leq}& 
	-\, \alpha\inner{\nabla f(x) - \nabla f(\xs)}{\xt} \,-\, \tfrac{\alpha}{\mu}\norm{E\xt \,+\, F\zt}^2 \,-\, \tfrac{\alpha}{\mu}(\norm{\zt}^2  \,-\, 2\inner{\wt}{\zt} \,+\, \norm{\wt}^2)
	\\[\nvs]
	&\asp{\leq}&
	-\, \tfrac{\alpha}{L_f}\norm{\nabla f(x) - \nabla f(\xs)}^2 \,-\, \tfrac{\alpha}{\mu}\norm{E\xt \,+\, F\zt}^2 \,-\, \tfrac{\alpha}{\mu}\norm{\zt \,-\, \wt}^2.
	\ea
	\non
	\eeq
	where the second inequality is obtained by the cocoercivity of $\nabla f$~\cite[Cor.~18.17]{baucom11}. Rearrangement of terms completes the proof.
\end{proof}

\subsection{Proof of Lemma~\ref{lemma.gradient_dual}}\label{proof.lemma.gradient_dual}
\begin{mylemma}{2}
	The gradient of the dual function $\dual (\yy, \lam) $,
	\beq\non
		\nabla\dual(\yy, \lam)  \,=\, \tbo{\nabla_\yy\dual(\yy, \lam)}{\nabla_\lam\dual(\yy, \lam)} \,=\, \tbo{\zb (\yy, \lam) \,-\, \prox_{\mu g}(\zb (\yy, \lam) \,+\, \mu\yy)}{E\xb (\yy, \lam) \,+\, F\zb (\yy, \lam) \,-\, q}  
	\eeq
	is Lipschitz continuous with modulus $\mu$, where $(\xb, \zb)$ denotes a $(\yy, \lam)$-parameterized solution to~\eqref{eq.bar_system}.
\end{mylemma}

\vsp

\begin{proof}
We first show that quantities $Ex + Fz - q$ and $z - \prox_{\mu g}(z + \mu\yy)$ remain constant over the set of solutions to~\eqref{eq.bar_system} at $(\yy, \lam)$. Let $(\xb,  \zb)$ and $(\xbp, \zbp)$ be two different solutions to~\eqref{eq.bar_system} at $(\yy, \lam)$.  Let $\wb \DefinedAs \prox_{\mu g}(\zb + \mu\yy)$ and $\wbp\DefinedAs \prox_{\mu g}(\zbp + \mu\yy)$. Suppose, for contradiction, that either $E\xb + F\zb \neq E\xbp + F\zbp$ or $\zb - \wb  \neq \zbp - \wbp$. Since the augmented Lagrangian~\eqref{eq.aug_lagrangian} is a convex function over the primal variables, set of its minimizers is also convex, which means that $(\xh, \zh , \wh ) = ((\xb \,+\, \xbp)/2, (\zb \,+\, \zbp)/2, (\wb \,+\, \wbp)/2)$ is also a minimizer. Moreover, Lagrangian~\eqref{eq.lagrange} is a convex function over the primal variable; hence,
	\beq\label{eq.cvx_lag}
		\tfrac{1}{2}(\cL(\xb, \zb, \wb; \yy, \lam) \,+\, \cL(\xbp, \zbp, \wbp; \yy, \lam))  \,\geq\,  \cL(\xh, \zh, \wh; \yy, \lam).
	\eeq
Since $\norm{\cdot}^2$ is a strongly convex function and the arguments are not equal by the initial supposition, we have the following inequalities with at least one of them being strict
	\beq
	\ba{rcl}\label{eq.cvx_norm}
		\tfrac{1}{2}(\norm{E\xb \,+\, F\zb \,-\, q}^2 \,+\, \norm{E\xbp \,+\, F\zbp \,-\, q}^2) 
		&\asp{\geq}&
		\norm{E\xh \,+\, F\zh \,-\, b}^2
		\\[\nvs]
		\tfrac{1}{2}(\norm{\zb  \,-\, \wb}^2 \,+\, \norm{\zbp \,-\, \wbp}^2)
		&\asp{\geq}&
		\norm{\zh  \,-\, \wh}^2.
	\ea
	\eeq 
Summing~\eqref{eq.cvx_lag}~and~\eqref{eq.cvx_norm} gives
	\beq
		\dual(\yy,\ \lam) \,=\, \tfrac{1}{2}(\cLm(\xb, \zb, \wb; \yy, \lam) \,+\, \cLm(\xbp, \zbp, \wbp; \yy, \lam))  \,>\,  \cLm(\xh, \zh, \wh; \yy, \lam)
	\non
	\eeq
which is a contradiction since $\dual(\yy, \lam)$ is the minimum of $\cLm$ over all the primal variables. Consequently, the subdifferential of the dual function is a singleton, and by~\cite[Cor.~10.14]{rocwet09}, $\nabla\dual$ exists and is given by~\eqref{eq.lemma.grad_dual}.

\vsp

Next, we prove the Lipschitz continuity of $\nabla \dual$. Let $(\yy, \lam)$ and $(\yp, \lamp)$ be arbitrary points and let $(\xb, \zb)$ and  $(\xbp, \zbp)$ be any solutions to~\eqref{eq.bar_system} at $(\yy, \lam)$ and $(\yp, \lamp)$, respectively. Let $\yt \DefinedAs \yy - \yp$, $\lamt \DefinedAs \lam - \lamp$, $\xbt \DefinedAs \xb - \xbp$, $\zbt \DefinedAs \zb - \zbp$, and $\wbt \DefinedAs \prox_{\mu g}(\zb  + \mu\yy) - \prox_{\mu g}(\zbp + \mu\yp)$. Since  $(\xb,  \zb)$ and $(\xbp,  \zbp)$ minimize $\cLm$ at $(\yy,\lam)$ and $(\yt, \lamt)$, respectively, we have
	\beq
	\ba{rcl}
		0 
		&\asp{=}&
		\inner{\nabla_{x,\zz}\cLm(\xb, \zb; \yy, \lam) \,-\, \nabla_{x,\zz}\cLm(\xbp, \zbp; \yp, \lamp)}{(\xbt, \zbt)}
		\\[\nvs]
		&\asp{=}&
		\inner{\nabla f(\xb) \,-\, \nabla f(\xbp)  \,+\, \tfrac{1}{\mu}E^T(E\xbt \,+\, F\zbt \,+\, \mu\lamt)}{\xbt} \,+\, \inner{\yt \,+\, \tfrac{1}{\mu}(\zbt \,-\, \wbt) \,+\, \tfrac{1}{\mu}F^T(E\xbt \,+\, F\zbt \,+\, \mu\lamt)}{\zbt}
		\\[\nvs]
		&\asp{\geq}&
		\inner{\tfrac{1}{\mu}E^T(E\xbt \,+\, \zbt \,+\, \mu\lamt)}{\xbt} \,+\, \inner{\yt \,+\, \tfrac{1}{\mu}(\zbt \,-\, \wbt) \,+\, \tfrac{1}{\mu}F^T(E\xbt \,+\, F\zbt \,+\, \mu\lamt)}{\zbt}
		\\[\nvs]
		&\asp{=}&
		\tfrac{1}{\mu}\norm{E\xbt \,+\, F\zbt}^2 \,+\, \inner{\lamt}{E\xbt \,+\, F\zbt} \,+\, \inner{\yt \,+\, \tfrac{1}{\mu}(\zbt \,-\, \wbt)}{\zbt}
		\\[\nvs]
		&\asp{=}&
		\tfrac{1}{\mu}\norm{E\xbt \,+\, F\zbt}^2 \,+\, \inner{\lamt}{E\xbt \,+\, F\zbt} \,+\, \inner{\yt}{\zbt} \,+\, \tfrac{1}{\mu}\norm{\zbt \,-\, \wbt}^2  \,-\, \tfrac{1}{\mu}\norm{\wbt}^2 \,+\, \tfrac{1}{\mu}\inner{\wbt}{\zbt}
		\\[\nvs]
		& \asp{\geq} &
		\tfrac{1}{\mu}\norm{E\xbt \,+\, F\zbt}^2 \,+\, \inner{\lamt}{E\xbt \,+\, F\zbt} \,+\, \inner{\yt}{\zbt} \,+\, \tfrac{1}{\mu}\norm{\zbt \,-\, \wbt}^2  \,-\, \tfrac{1}{\mu}\inner{\wbt}{\zbt \,+\, \mu\yt} \,+\, \tfrac{1}{\mu}\inner{\wbt}{\zbt}
		\\[\nvs]
		& \asp{=} &
		\tfrac{1}{\mu}\norm{E\xbt \,+\, F\zbt}^2 \,+\, \inner{\lamt}{E\xbt \,+\, F\zbt} \,+\, \inner{\yt}{\zbt \,-\, \wbt} \,+\, \tfrac{1}{\mu}\norm{\zbt \,-\, \wbt}^2
		\\[\nvs]
		& \asp{=} &
		\tfrac{1}{\mu}\norm{\nabla \dual(\yy, \lam) \,-\, \nabla\dual(\yp, \lamp)}^2 \,+\, \inner{\nabla \dual(\yy, \lam) \,-\, \nabla\dual(\yp, \lamp)}{(\yt, \lamt)}
	\ea
	\non
	\eeq
where the first inequality follows from the monotonicity of $\nabla f$, the third equality follows from the symmetry of inner products, the forth equality follows from completing the square, the second inequality follows from the non-expansiveness of the proximal operator, the fifth equality follows from the linearity of inner products, and the last equality follows from~\eqref{eq.lemma.grad_dual}. The proof is completed by using the Cauchy-Schwarz inequality. 
\end{proof}
		\vspace*{-1ex}
\subsection{Proof of Lemma~\ref{lemma.eb_primal}}\label{proof.lemma.eb_primal}
\begin{mylemma}{3}
Let Assumptions~\ref{ass.zero},~\ref{ass.eb_smooth}, and~\ref{ass.eb_nonsmooth} hold. There exist positive constants $\kappa_p$ and $\delta_p$ such that the following inequalities hold when $\norm{\nabla_{x,\zz}\cLm(x,\zz;\yy,\lam)}\leq \delta_p$,
	\bseq
	\begin{eqnarray}
\kappa_p\dist((x,\zz),\cP(\yy,\lam))  &\asp{\leq}&  \norm{\nabla_{x,\zz}\cLm(x,\zz;\yy,\lam)} \label{eq.lemma.eb_primal.eb.proof}
\\
(\kappa_p/2)\dist^2((x,\zz),\cP(\yy,\lam)) &\asp{\leq}& \cLm(x,\zz;\yy,\lam) \,-\, \dual(\yy,\lam)\label{eq.lemma.eb_primal.qg.proof}
\\
\cLm(x,\zz;\yy,\lam) \,-\, \dual(\yy,\lam) &\asp{\leq}& ({L_{xz}}/{2\kappa_p})\norm{\nabla_{x,\zz}\cLm(x,\zz;\yy,\lam)} ^2 \label{eq.lemma.eb_primal.pl.proof}
\end{eqnarray}
\eseq
where $L_{x\zz}$ is the Lipschitz constant of~$\nabla_{x,\zz}\cLm$.
\end{mylemma}

\vsp

\begin{proof}
The proof is based on the Hoffman error-bound condition associated with generalized gradient map of composite objective functions~\cite{luotse92,tse10}. We consider minimizing the augmented Lagrangian~\eqref{eq.aug_lagrangian} with respect to primal variables $(x,\zz,\ww)$ and denote the set of all minimizers at a given dual pair $(\yy,\lam)$ by $\cPw(\yy,\lam)$. Due to~\eqref{eq.wopt}, we have 
\beq\label{eq.dist_w}
\cPw(\yy,\lam)  \,=\,  \set{(x,\zz,\prox_{\mu g}(\zz \,+\, \mu\yy))\,|\,(x,\zz)\in\cP(\yy,\lam)}.
\eeq
Under Assumptions~\ref{ass.eb_smooth} and~\ref{ass.eb_nonsmooth}, the error bound conditions~\cite[Lemma 7]{tseyun09} and~\cite[Theorem 2]{tse10} (see~\cite{zhoso17} for a recent overview of related results) imply the existence of positive constants $\kappa_p$ and $\delta$ such that the distance to $\cPw(\yy,\lam)$ at any $(\yy,\lam)$ is upper bounded by the magnitude of the generalized gradient map associated with the augmented Lagrangian, i.e., the following inequality holds 
\beq\non
\ba{rcl}
\kappa_p\dist((x,\zz,\ww),\cPw(\yy,\lam))  &\asp{\leq}& \norm{\cG_{\cLm}(x,\zz,\ww;\yy,\lam)}
\ea
\eeq 
when $\norm{\cG_{\cLm}(x,\zz,\ww;\yy,\lam)}\leq \delta_p$. Here, the generalized gradient map $\cG_{\cLm}$ is given by
\beq\non
\ba{rcl}
\cG_{\cLm}(x,\zz,\ww;\yy,\lam) &\asp{=}&\dfrac{1}{\mu}\thbo{x  \,-\, \prox_{0}(x \,-\, \mu\nabla_x\cLm(x,\zz,\ww;\yy,\lam))}{\zz  \,-\, \prox_{0}(\zz \,-\, \mu\nabla_\zz\cLm(x,\zz,\ww;\yy,\lam))}{\ww  \,-\, \prox_{\mu g}(\ww \,-\, \mu\nabla_\ww\left[ \cLm(x,\zz,\ww;\yy,\lam)  \,-\, g(\ww)\right] )}
\\[0.6cm]
&\asp{=}&
\thbo{\nabla_x\cLm(x,\zz,\ww;\yy,\lam)}{\nabla_\zz\cLm(x,\zz,\ww;\yy,\lam)}{(1/\mu)(\ww  \,-\, \prox_{\mu g}(\zz \,+\, \mu\yy )}
\ea
\eeq
where the second equality is obtained using the fact that $\prox$ operator associated with zero is identity map. Now, since the third entry of $\cG_{\cLm}$ is zero at $w = \prox_{\mu g}(\zz \,+\, \mu\yy)$, identity~\eqref{eq.dist_w} together with the definition of proximal augmented Lagrangian~\eqref{eq.pal} implies that
\beq\non
\kappa_p\dist((x,\zz,\wb(\zz,\yy)),\cPw(\yy,\lam))  \,=\, \kappa_p\dist((x,\zz),\cP(\yy,\lam))  \,\leq\, \norm{\cG_{\cLm}(x,\zz,\wb(\zz,\yy);\yy,\lam)}   \,=\, \norm{\nabla_{x,\zz}\cLm(x,\zz;\yy,\lam)}.
\eeq
Moreover, since the proximal augmented Lagrangian is a smooth convex function in primal variables (see~\cite{dinhudhijovCDC18} for an explicit expression of $L_{xz}$), the equivalence between the error bound, PL, and quadratic growth conditions~\cite[Theorem 2]{karnutsch16} yields~\eqref{eq.lemma.eb_primal.qg.proof} and~\eqref{eq.lemma.eb_primal.pl.proof}.
\end{proof}

		\vspace*{-1ex}
\subsection{Proof of Lemma~\ref{lemma.eb_dual}}\label{proof.lemma.eb_dual}
\begin{mylemma}{4}
Let Assumptions~\ref{ass.zero},~\ref{ass.eb_smooth}, and~\ref{ass.eb_nonsmooth} hold. There exist positive constants $\kappa_d$ and $\delta_d$ such that the following inequality holds when $\norm{\nabla\dual(\yy,\lam)}\leq \delta_d$,
	\bseq
	\begin{eqnarray}
\kappa_d\dist((\yy,\lam),\cDs)  &\asp{\leq}&\norm{\nabla\dual(\yy,\lam)}  \label{eq.lemma.eb_dual.eb.proof}
\\
(\kappa_d/2)\dist^2((\yy,\lam),\cDs) &\asp{\leq}&  \dopt \,-\, \dual(\yy,\lam) \label{eq.lemma.eb_dual.qg.proof}
\\
\dopt \,-\, \dual(\yy,\lam) &\asp{\leq}& ({\mu}/{2\kappa_d})\norm{\nabla\dual(\yy,\lam)} ^2 \label{eq.lemma.eb_dual.pl.proof}
\end{eqnarray}
\eseq
\end{mylemma}

\vsp

\begin{proof}
The proof follows from~\cite[Lemma 2.3-(c)]{honluo17}. For completeness, we verify the conditions in~\cite[Lemma 2.3]{honluo17}: The lifted problem~\eqref{eq.lifted} can be written as
\beq\non
\ba{rcl}
	\minimize\limits_{\xt}
	& \asp{} &
	\ft_1(x) \,+\, \ft_2(\zz) \,+\, \ft_3(\ww)
	\\[\nvs]
	\subjectto 
	& \asp{} &
	\tbo{E}{0}x\,+\, \tbo{F}{I}\zz + \tbo{0}{-I}\ww\,=\, \tbo{q}{0}
\ea
\eeq
where $\ft_1(x) \DefinedAs f(x)$, $\ft_2(\zz) = 0$, and $\ft_3(\ww) = g(\ww)$. Condition (a) and (e) in~\cite[Lemma 2.3]{honluo17} is verified by Assumption~\ref{ass.zero}, condition (d) by Assumption~\ref{ass.eb_nonsmooth}, and conditions (b) and (c) by Assumption~\ref{ass.eb_smooth}. Moreover, since the dual function is concave and has a Lipschitz continuous gradient with modulus $\mu$ (see Lemma~\ref{lemma.gradient_dual}), the equivalence between the error bound, PL, and quadratic growth conditions~\cite[Theorem 2]{karnutsch16} yields~\eqref{eq.lemma.eb_dual.qg.proof} and~\eqref{eq.lemma.eb_dual.pl.proof}. For the inclusion of group lasso penalty function in Assumption~\ref{ass.eb_nonsmooth}, see the remark after Lemma 2.3 in~\cite{honluo17}.
\end{proof}

		\vspace*{-1ex}
\subsection{Proof of Lemma~\ref{lemma.lyap2_upper}}\label{proof.lemma.lyap2_upper}
\begin{mylemma}{5}
	Lyapunov function $V_2$ in \eqref{eq.lyap2} satisfies
	\beq\non
	V_2(x, \zz; \yy, \lam)   \,\leq\, c_1\dist^2((x,\zz,\yy,\lam),\PPs) 
	\eeq 
	where $c_1 = (L_{xz}/2+1)\max(1,\mu)$ and $L_{x\zz}$ is the Lipschitz constant of~$\nabla_{x,\zz}\cLm$.
\end{mylemma}

\vsp

\begin{proof}
The dual gap can be bounded by using Lipschitz continuity of $\nabla \dual$ as
\beq\label{eq.upper_dualgap}
\dopt \,-\, \dual(\yy,\lam) \,\leq\, (\mu/2)\dist^2((\yy,\lam),\cDs).
\eeq
As for the quadratic upper bound on the primal gap, let $(\xb,\zb) \in\cP(\yy,\lam)$. Adding and subtracting $\dopt$ from the primal gap yields
\beq\label{eq.dual_diff}
\cLm(x,\zz;\yy,\lam)   \,\mp\, \dopt \,-\,  \cLm(\xb,\zb;\yy,\lam)  \,=\, \cLm(x,\zz;\yy,\lam)  \,-\, \cLm(\xs,\zs;\ys,\lams)   \,+\, \dopt  \,-\, \dual(\yy,\lam)
\eeq
where $(\xs,\zs,\ys,\lams)$ is an arbitrary point in $\PPs$. The second difference term in~\eqref{eq.dual_diff}, i.e., the dual gap, is bounded by~\eqref{eq.upper_dualgap}, while the first difference can be bounded by using saddle inequality~\eqref{eq.saddle} as
\beq\non
\cLm(x,\zz;\yy,\lam)   \,-\, \cLm(\xs,\zs;\ys,\lams) \,\leq\, \cLm(x,\zz;\yy,\lam)   \,-\, \cLm(\xs,\zs;\yy,\lam).
\eeq 	
Since $\nabla_{x,\zz}\cLm$ is smooth, the quadratic upper bound yields
\beq\label{eq.inner_prod_up}
\cLm(x,\zz;\yy,\lam)   \,-\, \cLm(\xs,\zs;\yy,\lam)  \,\leq\, \inner{\nabla_{x,\zz}\cLm(\xs,\zs;\yy,\lam)}{(x,\zz) \,-\, (\xs,\zs)}  \,+\, \tfrac{L_{xz}}{2}\norm{(x,\zz) \,-\, (\xs,\zs)}^2.
\eeq 
The inner product can be upper bounded using Fenchel-Young inequality as follows:  
\begin{align}
\inner{\nabla_{x,\zz}\cLm(\xs,\zs;\yy,\lam)}{(x,\zz) - (\xs,\zs)} 
&\stackrel{(i) }{ \,\leq\, }\tfrac{1}{2}(\norm{\nabla_{x,\zz}\cLm(\xs,\zs;\yy,\lam)}^2 \,+\,\norm{(x,\zz) - (\xs,\zs)}^2)
\non
\\
& \stackrel{(ii) }{ \,\leq\, }L_{xz}(\cLm(x,\zz;\yy,\lam)   - \dual(\yy,\lam)) \,+\, \tfrac{1}{2}\norm{(x,\zz) - (\xs,\zs)}^2
\non
\\
& \stackrel{(iii) }{ \,\leq\, }L_{xz}(\dopt   - d(\yy,\lam)) \,+\, \tfrac{1}{2}\norm{(x,\zz) - (\xs,\zs)}^2
\non
\\
& \stackrel{(iv) }{ \,\leq\, }\tfrac{\mu L_{xz}}{2}\dist^2((\yy,\lam),\cDs) \,+\, \tfrac{1}{2}\norm{(x,\zz) - (\xs,\zs)}^2.
\label{eq.upper_inner}
\end{align}
Here, $(ii)$ is given by the Lipschitz continuity of $\nabla_{x,\zz}\cLm$ and convexity of $\cLm$ in primal variables, which yields
\beq\non
\hspace{-2ex}\tfrac{1}{2L_{xz}}\norm{\nabla_{x,\zz}\cLm(\xs,\zs;\yy,\lam)}^2 \leq \cLm(\xs,\zs;\yy,\lam)  - \dual(\yy,\lam)
\eeq
whereas $(iii)$ is given by saddle inequality~\eqref{eq.saddle}
\beq\non
\begin{split}
	\cLm(\xs,\zs;\yy,\lam)  &\,-\, \cLm(\xb,\zb;\yy,\lam)  \,\leq\, \cLm(\xs,\zs;\ys,\lams) \,-\, \cLm(\xb,\zb;\yy,\lam).
\end{split}
\eeq
Lastly, $(iv)$ follows from~\eqref{eq.upper_dualgap}. Substituting~\eqref{eq.upper_inner} in~\eqref{eq.inner_prod_up}~gives
\begin{align*}
&\cLm(x,\zz;\yy,\lam)   \,-\, \cLm(\xs,\zs;\yy,\lam)  \,\leq\, \tfrac{\mu L_{xz}}{2}\norm{(\yy,\lam) \,-\, (\ys,\lams)}^2  \,+\, \tfrac{1 \,+\, L_{x\zz}}{2}\norm{(x,\zz) \,-\, (\xs,\zs)}^2  .
\end{align*}
The result follows from the fact that $(\xs,\zs,\ys,\lams)$ is an arbitrary solution.
\end{proof}

		\vspace*{-1ex}
\subsection{Proof of Lemma~\ref{lemma.lyap2}}\label{proof.lemma.lyap2}
\begin{mylemma}{6}
	Let Assumptions~\ref{ass.zero},~\ref{ass.eb_smooth}, and~\ref{ass.eb_nonsmooth} hold and let $\tb\geq t_0$ be such that $\norm{\nabla_{x,\zz}\cLm(x(\tb),\zz(\tb);\yy(\tb),\lam(\tb))}\leq \delta_p$ and $\norm{\nabla\dual(\yy(\tb),\lam(\tb))}\leq \delta_d$ for constants $\delta_p$ and $\delta_d$ given in Lemmas~\ref{lemma.eb_primal} and~\ref{lemma.eb_dual}, respectively. The time derivative of $V_2$ along the solutions of~\eqref{eq.dyn} with a positive time scale $\alpha \in (0,\ab_1)$ satisfies
	\beq\non
		\dot{V}_2(t)   \,\leq\,   -\, \rho_1 V_2(t) \quad \forall t\geq \tb
	\eeq
	where $\ab_1 = 0.5\kappa_p^2\left(\sigu^2([E\ F]) +  4\right)^{-1}$ and $\rho_1 = \min(0.5,\alpha)/\max\left( L_{xz}/(2\kappa_p), \mu/(2\kappa_d)\right)$.
\end{mylemma}

\vsp 

\begin{proof}
The time derivative of $V_2$ along the solutions of~\eqref{eq.dyn} can be obtained by using the chain rule as 
	\begin{align}
		\dot{V}_2  
		&\,=\,
		\inner{\nabla_x\cLm}{\xd} \,+\, \inner{\nabla_\zz\cLm}{\zd} \,+\, \inner{\nabla_\yy\cLm \,-\, 2\nabla_\yy\dual}{\yd} \,+\, \inner{\nabla_\lam\cLm \,-\, 2\nabla_\lam\dual}{\lamd} 
		\\[\nvs]
		&\,=\,
		-\, \norm{\nabla_{x,\zz}\cLm}^2 \,+\, \alpha\inner{\nabla_\yy\cLm \,-\, 2\nabla_\yy\dual}{\nabla_\yy\cLm} \,+\, \alpha\inner{\nabla_\lam\cLm \,-\, 2\nabla_\lam\dual}{\nabla_\lam\cLm}
		\\[\nvs]
		&\,=\,
		-\, \norm{\nabla_{x,\zz}\cLm}^2\,+\, \alpha\norm{\nabla_\yy\cLm \,-\, \nabla_\yy\dual}^2 \,+\, \alpha\norm{\nabla_\lam\cLm \,-\, \nabla_\lam\dual}^2  \,-\, \alpha\norm{\nabla_\yy \dual}^2\ \,-\, \alpha\norm{\nabla_\lam \dual}^2 
		\label{eq.v2dot_bound}  
	\end{align}
where the last equality is obtained via completing the square. Let $(\xb, \zb)$ be an arbitrary point in $\cP(\yy, \lam)$. Using the gradient expression~\eqref{eq.lemma.grad_dual}, we bound the first positive term \blue{in~\eqref{eq.v2dot_bound} using the triangle inequality and the} nonexpansiveness of the proximal operator \blue{as follows},
	\begin{align}
		\norm{\nabla_\yy\cLm \,-\, \nabla_\yy\dual}^2
		&=
		\norm{\zz \,-\, \prox_{\mu g}(\zz \,+\, \mu\yy) \,-\, \left(\zb \,-\, \prox_{\mu g}(\zb \,+\, \mu\yy)\right)}^2\non
		\\
		&\leq
		2\left( \norm{\zz \,-\, \zb}^2 \,+\, \norm{\prox_{\mu g}(\zz \,+\, \mu\yy) \,-\, \prox_{\mu g}(\zb \,+\, \mu\yy)}^2\right) \non
		\\
		&\leq
		4\norm{\zz \,-\, \zb}^2\label{eq.proof.v2dot1}
	\end{align}
where the first line is obtained by using the triangle inequality and the third line follows from the non-expansiveness of the proximal operator. As for the second positive term, we have 
	\begin{align}
		\norm{\nabla_\lam\cLm \,-\, \nabla_\lam\dual}^2
		&=
		\norm{Ex \,+\, F\zz \,-\, (E\xb \,+\, F\zb)}^2\non
		\\
		&\leq
		\sigu^2([E \ F])\norm{(x, \zz) \,-\, (\xb, \zb)}^2.\label{eq.proof.v2dot2}
	\end{align}
Substituting~\eqref{eq.proof.v2dot1} and~\eqref{eq.proof.v2dot2} back into $\dot{V}_2$ and using the fact that $(\xb, \zb)$ is an arbitrary point in $\cP(\yy,\lam)$ yields
	\beq\label{eq.lyap2_time_der}
		\dot{V}_2
		 \,\leq\, 
		-\, \norm{\nabla_{x,\zz}\cLm}^2  \,+\, \alpha(4 \,+\, \sigu^2([E \ F])) \dist^2((x, \zz), \cP(\yy, \lam))\,-\, \alpha\norm{\nabla \dual}^2.
	\eeq
Also, using~\eqref{eq.lemma.eb_primal.eb}, we obtain 
	\beq\non
		\dot{V}_2
		 \,\leq\, 
		-\, (1  \,-\, \alpha\kappa_p^{-2}(4 \,+\, \sigu^2([E \ F]))) \norm{\nabla_{x,\zz}\cLm}^2 \,-\, \alpha\norm{\nabla \dual}^2.
	\eeq
Choosing $\alpha$ as in the lemma results in
	\beq\non
		\dot{V}_2
		 \,\leq\, 
		-\, \min(0.5,\alpha) (\norm{\nabla_{x,\zz}\cLm}^2 \,+\, \norm{\nabla \dual}^2).
	\eeq
which together with PL inequalities~\eqref{eq.lemma.eb_primal.pl} and~\eqref{eq.lemma.eb_dual.pl} concludes the proof.
\end{proof}

		\vspace*{-1ex}
\subsection{Proof of Lemma~\ref{lemma.upper_lam}}\label{proof.lemma.upper_lam}
\begin{mylemma}{7}
	Let Assumptions~\ref{ass.zero} and~\ref{ass.lip} hold and let $\pas = (\xas, \zas, \yas, \lamas) \DefinedAs \lim_{t\to\infty}\pp(t)$. Then,
	\beq\non
		\norm{\lam(t) \,-\, \lamas}^2 \,\leq\, c_2(\norm{\xb(t) \,-\, \xas}^2 \,+\, \norm{\nabla_\lam\dual(\yy(t), \lam(t))}^2),\quad \forall t \,\geq\, t_0
	\eeq
	where $(\xb(t), \zb(t))$ is an arbitrary point in $\cP(\yy(t), \lam(t))$, $c_2 = \max(2L_f^2/\sigl^2(E),\,  1/\mu^{2})$, and $L_f$ is the Lipschitz constant of $\nabla f$. 
\end{mylemma}

\vsp

\begin{proof}
We utilize~\eqref{eq.bar_system1} to upper bound $\norm{\lam - \lams}$. Since $\lamd \in \cR([E\ F])$ by Assumption~\ref{ass.zero}, the fundamental theorem of calculus yields that $\lam(t) - \lam(t_0)\in\cR([E\ F])$ for any $t\geq t_0$. Moreover, Theorem~\ref{theorem.gas} ensures that the solutions of~\eqref{eq.dyn} converge to a point in $\PPs$, i.e. $\pas \in \PPs$. Hence, using~\eqref{eq.bar_system1}, we obtain
\beq\non
\norm{\nabla f(\xb) \,-\, \nabla f(\xas)}^2 
 \,=\, 
\norm{E^T(\lam  \,-\,  \lamas  \,+\, (1/\mu)\nabla_\lam\dual(\yy,\lam))}^2.
\eeq
Under Assumption~\ref{ass.lip}, both $\lam(t) - \lamas$ and $ \nabla_\lam\dual(\yy,\lam)$ do not have any component in $\cN(E^T)$, which implies
\beq\non
\norm{\nabla f(\xb) \,-\, \nabla f(\xas)}^2 
 \,\geq\, 
\sigl^2(E)\norm{\lam  \,-\,  \lamas  \,+\,  \tfrac{1}{\mu}\nabla_\lam\dual(\yy,\lam)}^2.
\eeq
The following basic inequality~\cite{denyin16},
	\beq\label{eq.ineq1}
		\norm{u \,+\, v}^2 \,\geq\, \tfrac{1}{1 \,+\, \zeta}\norm{u}^2 \,-\, \tfrac{1}{\zeta}\norm{v}^2,
		~~\forall u,v \,\in\,\reals{n},
		~~ \forall\zeta \,>\, 0 
	\eeq
	implies
\beq\non
\norm{\nabla f(\xb) \,-\, \nabla f(\xas)}^2 \,+\,    \tfrac{\sigl^2(E)}{\mu^2}\norm{\nabla_\lam\dual(\yy,\lam)}^2 
 \,\geq\, 
\tfrac{\sigl^2(E)}{2}\norm{\lam  \,-\,  \lamas}^2,
\eeq
which together with the Lipschitz continuity of $\nabla f$ completes the proof.
\end{proof}

		\vspace*{-1ex}
\subsection{Proof of Lemma~\ref{lemma.conseq_ass_str}}\label{proof.lemma.conseq_ass_str}
\begin{mylemma}{8}
Let Assumptions~\ref{ass.zero} and~\ref{ass.str_cvx} hold and let $\mu m_g \leq 1$.
\begin{enumerate}
\item[(a)] The proximal augmented Lagrangian is strongly convex in primal variables $(x,z)$ with modulus $m_{xz}$, see~\eqref{eq.mxz} for an explicit expression of  $m_{xz}$.
\item[(b)] There is a unique solution to problem~\eqref{eq.intro}, i.e., $\cPs = \set{(\xs,\zs)}$, while $\cDs$ may not be a singleton.
\item[(c)] The time derivative of quadratic Lyapunov function $V_1$ in~\eqref{eq.lyap1} along the solutions of~\eqref{eq.dyn} with any $\alpha>0$ satisfies
\beq\non
\dot{V}_1(t) \,\leq\, -\,\alpha m_{xz}\norm{(x(t), \zz(t)) \,-\, (\xs,\zs) }^2, \quad t \,\geq\, t_0.
\eeq 
\item[(d)] The time derivative of nonquadratic Lyapunov function $V_2$ in~\eqref{eq.lyap2} along the solutions of dynamics~\eqref{eq.dyn} with any $\alpha\in(0,\ab_2]$ satisfies
\beq\non
	\dot{V}_2(t)   \,\leq\,   -\, \min(0.5,\alpha)\left(\norm{(x(t),\zz(t)) \,-\, (\xb(t),\zb(t))}^2  \,+\, \norm{\nabla\dual(\yy(t), \lam(t))}^2\right), \quad t \,\geq\, t_0
\eeq
where $\set{(\xb(t), \zb(t))} = \cP(\yy(t), \lam(t))$ and $\ab_2 = 0.5m_{x\zz}^2\big/\left(\sigu^2([E\ F]) +  4\right)$.
\end{enumerate}
\end{mylemma}

\vsp

\begin{proof}
We exploit the quadratic term in~\eqref{eq.pal} to induce strong convexity along directions in which $f$ and $g$ lack it.
Let $I^c$ and $J^c$ denote the complementary sets of $I$ and $J$ defined in Assumption~\ref{ass.str_cvx}, i.e., $I^c \DefinedAs \set{1,\cdots,k} \setminus I$ and $J^c \DefinedAs \set{1,\cdots,\ell}\setminus J$. Also, let $x_I$ denote the collection of $x$ variables, $E_{I}$ the row-concatenation of matrices, and $f_I$ the sum of functions that are associated with the blocks indexed by $i \in I$, i.e., the smooth but not strongly convex blocks. Similarly, we define tuple $(\zz_J, g_J, f_J)$ for nonsmooth and not strongly convex blocks. Additionally, we denote the strong convexity constant of $f_{I^c}$ by $m_f$ and of $g_{I^c}$ by $m_g$. We also define $m_{fg}$ as $m_{fg} = \min(m_f, m_g)$ if both $m_f$ and $m_g$ are nonzero. In the case that one of $m_f$ and $m_g$ is zero, then $m_{fg}$ is equal to the other. If both $m_f$ and $m_g$ are zero, so is $m_{fg}$. If $m_{fg}$ is zero, then $[E~F]$ is a full-column rank matrix by Assumption~\ref{ass.str_cvx}.

\vsp

\noindent\textbf{Strong Convexity of $\cLm$ in primal variables.} For arbitrary points $(x, \zz)$ and $(\xp, \zp)$, let $\xt  \DefinedAs x  -  \xp$ and $\zt  \DefinedAs \zz  -  \zp$. The strong convexity of $f_{I^c}$ and the contractive mapping $\prox_{\mu g_{J^c}}$~\cite[Proposition 23.13]{baucom11} yields
\begin{align}
		\inner{\nabla_{x,\zz}& \cLm(x,\zz; \yy,  \lam)   \,-\,   \nabla_{x,\zz} \cLm(\xp, \zp; \yy, \lam)}{(\xt, \zt)}\non
		\\
		&\,=\, \inner{\nabla f(x) - \nabla f(\xp)}{\xt} \,+\, \tfrac{1}{\mu}\norm{E\xt \,+\, F\zt}^2  \,+\,  \norm{\zt}^2  \,-\, \inner{\prox_{\mu g}(\zz \,+\, \mu\yy) \,-\, \prox_{\mu g}(\zp \,+\, \mu\yy)}{\zt}\non 
		 \\
		 &\,\geq\, 
		m_f\norm{\xt_{I^c}}^2 \,+\, \tfrac{1}{\mu}\norm{E\xt \,+\, F\zt}^2  \,+\, \tfrac{m_g}{\mu m_g  \,+\, 1}\norm{\zt_{J^c}}^2.
		\label{eq.str_cvx_v1}
\end{align}
There are two cases:
	\begin{itemize}
	\item If $m_{fg} = 0$, then $m_f = m_g = 0$ and Assumption~\ref{ass.str_cvx} guarantees that $[E~F]$ is full-column rank, which leads to
	\beq\label{eq.str_temp1}
	m_f\norm{\xt_{I^c}}^2 \,+\, \tfrac{1}{\mu}\norm{E\xt \,+\, F\zt}^2  \,+\, \tfrac{m_g}{2}\norm{\zt_{J^c}}^2
	 \,=\, 
	\tfrac{1}{\mu}\norm{E\xt \,+\, F\zt}^2  \,\geq\, \tfrac{\sigl^2([E\ F])}{\mu}\norm{(\xt, \zt)}^2.
	\eeq
	\item If $m_{fg} \neq 0$, then the inequality~\eqref{eq.ineq1} together with the definition of $m_{fg}$ for any $\zeta>0$ gives
	\begin{align}
	\hspace{-2ex}m_f\norm{\xt_{I^c}}^2  +  \tfrac{1}{\mu}\norm{E\xt + F\zt}^2  + \tfrac{m_g}{2}\norm{\zt_{J^c}}^2
	&\geq
	\tfrac{m_{fg}}{2}\norm{(\xt_{I^c}, \zt_{J^c})}^2 + \tfrac{1}{\mu (1 +  \zeta)}\norm{E_{I}x_{I} + F_{I}\zz_J}^2   - \tfrac{1}{\mu\zeta}\norm{E_{I^c}\xt_{I^c}  +  F_{J^c}\zt_{J^c}}^2
	\non
	\\
	&\geq
	\big(\tfrac{m_{fg}}{2} -  \tfrac{\sigu^2([E_{I^c}\ F_{J^c}])}{\mu\zeta} \big)\norm{(\xt_{I^c}, \zt_{J^c})}^2 +  \tfrac{\sigl^2([E_{I}\ F_{J}])}{\mu  +  \mu\zeta}\norm{(\xt_{I}, \zt_{J})}^2
	\non
	\\
	&\geq
	\tfrac{m_{fg}\sigl^2([E_I\ F_J])}{m_{fg}\mu \,+\, 4\sigu^2([E_{I^c}\ F_{J^c}])}\norm{(\xt, \zt)}^2
	\label{eq.str_temp2}
	\end{align}
	where the second line is obtained by using matrix norm and the third line is obtained by setting $\zeta =  {4\sigu^2([E_{I^c}\ F_{J^c}])}\big/(m_{fg}\mu)$.
	\end{itemize}

Recalling that $\mu m_g \leq 1$ and using~\eqref{eq.str_temp1} and~\eqref{eq.str_temp2} to lower bound~\eqref{eq.str_cvx_v1} yields the following strong convexity constant,
\beq
m_{xz}  \,\DefinedAs\, \left\{
\ba{ll}\label{eq.mxz}
\dfrac{\sigl^2([E\ F])}{\mu}, &m_{fg}  \,=\, 0
\\[0.4cm]
\dfrac{m_{fg}\sigl^2([E_I\ F_J])}{m_{fg}\mu \,+\, 4\sigu^2([E_{I^c}\ F_{J^c}])}, &m_{fg} \,>\,0.
\ea
\right.
\eeq

\vsp

\noindent\textbf{Uniqueness of the solution.} Let $(\xs, \zs,  \ys, \lams)$ and $(\xss, \zss, \yss, \lamss)$ be arbitrary solutions to~\eqref{eq.kkt} such that $(\xs, \zs) \neq (\xss,  \zss)$. The strong convexity of $f_{I^c}$ and $g_{J^c}$ yields
	\beq
	\ba{rcl}
		f(\xss) 
		& \asp{\geq} &
		f(\xs)  \,+\, \nabla f(\xs)^T(\xss \,-\, \xs) \,+\, \tfrac{m_f}{2} \norm{\xss_{I^c} \,-\, \xs_{I^c}}^2
		\\[\nvs]
		g(\zss) 
		& \asp{\geq} &
		g(\zs)  \,+\, r^T(\zss \,-\, \zs) \,+\, \tfrac{m_g}{2} \norm{\zss_{J^c} \,-\, \zs_{J^c}}^2.
	\ea
	\non
	\eeq 
Adding up these two inequalities results in
	\beq\label{eq.strong_convexity_IJ}
		-\tbo{\nabla f(\xs)}{r}^T\tbo{\xss \,-\, \xs}{\zss \,-\, \zs}
	 \,\geq\,  
  \tfrac{m_f}{2} \norm{\xss_{I^c} \,-\, \xs_{I^c}}^2 
  \,+\,  
  \tfrac{m_g}{2} \norm{\zss_{J^c} \,-\, \zs_{J^c}}^2 
  \,\geq\,
  0.
	\eeq
Since $(\xss - \xs,  \zss - \zs)$ is a feasible direction, the first order optimality condition yields that for any $r  \in \partial g(\zs)$, 
	\beq\non
		\tbo{\nabla f(\xs)}{r}^T\tbo{\xss \,-\, \xs}{\zss \,-\, \zs} \,\geq\,  0
	\eeq
which in conjunction with~\eqref{eq.strong_convexity_IJ} implies that $(\xs_{I^c}, \zs_{J^c}) = (\xss_{I^c}, \zss_{J^c})$. Thus, the initial assumption $(\xs, \zs) \neq (\xss,  \zss)$ yields that $(\xs_{I}, \zs_{J}) \neq (\xss_{I}, \zss_{J})$. On the other hand, optimality condition~\eqref{eq.opt3} requires $E(\xs - \xss) + F(\zs - \zss) = 0$, which together with $(\xs_{I^c}, \zs_{J^c}) = (\xss_{I^c}, \zss_{J^c})$ implies that $E_I(\xs_I - \xss_I) - F_J(\zs_J - \zss_J)  = 0$, but this is a contradiction because $(\xs_{I}, \zs_{J}) \neq (\xss_{I}, \zss_{J})$ and $[E_I,~F_J]$ is full-column rank by Assumption~\ref{ass.str_cvx}. In conclusion, $(\xs, \zs) = (\xss, \zss)$.

\vsp

\noindent\textbf{Time derivative of $V_1$.}
Let $\xt \DefinedAs x  -  \xs$, $\zt  \DefinedAs \zz - \zs$, $\yt \DefinedAs \yy - \ys$, $\lamt \DefinedAs \lam - \lams$, and $\wt  \DefinedAs \prox_{\mu g}(\zz  +  \mu\yy) - \prox_{\mu g}(\zs  + \mu\ys)$. Starting with equation~\eqref{eq.lyap0_derivative} derived in Section~\ref{proof.lemma.gradient_dual}, we partition the right hand side as the sum of terms related to strongly convex blocks and the rest:
		\beq
		\ba{l}
			\dot{V}_1
			 =  
			- \alpha\inner{\nabla f_I}{\xt_I}  -  \alpha\inner{\nabla f_{I^c}}{\xt_{I^c}} - \tfrac{\alpha}{\mu}\norm{E\xt + F\zt}^2 - \tfrac{\alpha}{\mu}(\norm{\zt}^2  - 2\inner{\wt}{\zt} \,+\, \inner{\wt_J}{\zt_{J} + \mu\yt_{J}} + \inner{\wt_{J^c}}{\zt_{J^c} + \mu\yt_{J^c}})
			\\[0.2cm]
			\leq
			-\, \alpha\inner{\nabla f_I}{\xt_I}  \,-\,  \alpha\inner{\nabla f_{I^c}}{\xt_{I^c}} \,-\, \tfrac{\alpha}{\mu}\norm{E\xt \,+\, F\zt}^2 \,-\, \tfrac{\alpha}{\mu}(\norm{\zt}^2  \,-\, 2\inner{\wt}{\zt} \,+\, \norm{\wt_J}^2 \,+\, (\mu m_g+1)\norm{\wt_{J^c}}^2)
			\\[0.2cm]
			\leq
			-\, \alpha\inner{\nabla f_{I^c}}{\xt_{I^c}} \,-\, \tfrac{\alpha}{\mu}\norm{E\xt \,+\, F\zt}^2 \,-\, \tfrac{\alpha}{\mu}(\norm{\zt \,-\, \wt}^2 \,+\, \mu m_g\norm{\wt_{J^c}}^2)
			\\[0.2cm]
			\leq
			-\, \alpha\, m_f\norm{\xt_{I^c}}^2 \,-\, \tfrac{\alpha}{\mu}\norm{E\xt \,+\, F\zt}^2 \,-\, \tfrac{\alpha}{\mu}(\norm{\zt \,-\, \wt}^2 \,+\, \mu m_g\norm{\wt_{J^c}}^2)
			\\[0.2cm]
			=
				-\, \alpha\, m_f\norm{\xt_{I^c}}^2 - \tfrac{\alpha}{\mu}\norm{E\xt + F\zt}^2 - \tfrac{\alpha}{\mu}(\norm{\zt_J - \wt_J}^2 + (1 - \mu m_g)\norm{\zt_{J^c} - \wt_{J^c}}^2 +  \mu m_g(\norm{\zt_{J^c} - \wt_{J^c}}^2 +\norm{\wt_{J^c}}^2))
				\\[\nvs]
			\leq
			-\, \alpha\, m_f\norm{\xt_{I^c}}^2 \,-\, \tfrac{\alpha}{\mu}\norm{E\xt \,+\, F\zt}^2 \,-\, \tfrac{\alpha}{\mu}(\norm{\zt_J \,-\, \wt_J}^2 \,+\, (1 \,-\, \mu m_g)\norm{\zt_{J^c} \,-\, \wt_{J^c}}^2 \,+\,  \tfrac{\mu m_g}{2}\norm{\zt_{J^c}}^2)
			\\[0.2cm]
			\leq
			-\, \alpha\, m_f\norm{\xt_{I^c}}^2 \,-\, \tfrac{\alpha}{\mu}\norm{E\xt \,+\, F\zt}^2 \,-\, \tfrac{\alpha m_g}{2}\norm{\zt_{J^c}}^2
			\\[0.2cm]
			\leq
			-\alpha m_{xz}\norm{(\xt, \zt)}^2
		\ea		
		\non
		\eeq
		where the first inequality follows from the nonexpansiveness of $\prox_{\mu g_{J}}$ and contraction of $\prox_{\mu g_{J^c}}$, the second inequality from the monotonicity of $\nabla f_I$, the third inequality from the strong monotonicity of $\nabla f_{I^c}$, the forth inequality from the triangle inequality, the fifth inequality from the removal of negative terms, the sixth inequality from definition of $m_{x\zz}$ together with inequalities~\eqref{eq.str_temp1} and~\eqref{eq.str_temp2}.
	
\vsp

\noindent\textbf{Time Derivative of $V_2$.}
The time derivative of $V_2$ along the solutions of~\eqref{eq.dyn} is already obtained in~\eqref{eq.lyap2_time_der} for the proof of Lemma~\ref{lemma.lyap2}:
	\beq
\dot{V}_2
		 \,\leq\, 
		-\, \norm{\nabla_{x,\zz}\cLm}^2  \,+\, \alpha(4 \,+\, \sigu^2([E \ F])) \norm{(x,z) \,-\, (\xb,\zb)}^2 \,-\, \alpha\norm{\nabla \dual}^2.
	\non
	\eeq
The strong convexity of the proximal augmented Lagrangian in primal variables yields
	\beq\non
		\dot{V}_2
		 \,\leq\, 
		-\,\left( m_{xz}^2  \,-\, \alpha(4 \,+\, \sigu^2([E \ F])) \right)\norm{(x,\zz) \,-\, (\xb,\zb)}^2 \,-\, \alpha\norm{\nabla \dual}^2. 
	\eeq
Choosing $\alpha$ as in the lemma concludes the proof.
\end{proof}
\end{document}